\title{ALGEBRA OF FAMILIES OF ALTERNATING KNOTS AND LINKS}
\author{E. PI\~NA \\
Department of Physics\\ Universidad Aut\'onoma Metropolitana - Iztapalapa, \\
P. O. Box 55 534 \ Mexico, D. F., 09340 Mexico \\
e-mail: pge@xanum.uam.mx}

\documentclass[12pt]{article}
\usepackage{graphicx}
\usepackage{psfrag}
\usepackage{amsfonts}
\begin{document}
\date {In the centenary of H. Poincar\'e (1854-1912)}
\maketitle

\abstract{Families of alternating knots (links) and tangles are studied using as building block the conway defined as the twisting of two strands.
The regular representation of knots assumes the projection has the minimal number of overpassings, and the minimal number of conways. The continued
fraction associated to rational knots is represented by gaussian brackets and products of 2-dimensional matrices. This gives birth to an algebra
of rational knots and tangles which is easily generalized to alternating knots. A collection of 65 families of prime alternating knots
with one to six conways is found. Eleven families with six conways show peculiar behavior not present in families with a lower or equal number of conways.}

\newpage

\section{INTRODUCTION}
This paper deals exclusively with alternating knots and alternating links, namely the string forming the knot (link) overpasses and underpasses
alternately at the crossing points. Fig. 1 shows examples of both an alternating and a non-alternating knot. The word knot will encompass both knots
and links except when distinction is relevant.

In 1970 J. Conway \cite{co} published a seminal work on knots where the simplest prime knots are tabulated. For each knot J. Conway provides in his
table of knots a natural number $C \ \epsilon \ \mathbb{N} \cup 0$, which will be here referred as Conway's number, equal to the value of the Alexander
polynomial \cite{al} for $x = - 1$, namely, equal
to the sum of the absolute values of the coefficients of the Alexander polynomial, which plays a relevant role in the theory of knots. The coefficients
of the Alexander polynomial appear next to the corresponding knot in the table of knots drawings in Appendix C of D. Rolfsen's \cite{ro} book
on knots. Rolfsen calls to Conway's number the determinant of the knot.

\begin{figure}
\quad \quad \quad \scalebox{0.3}{\includegraphics{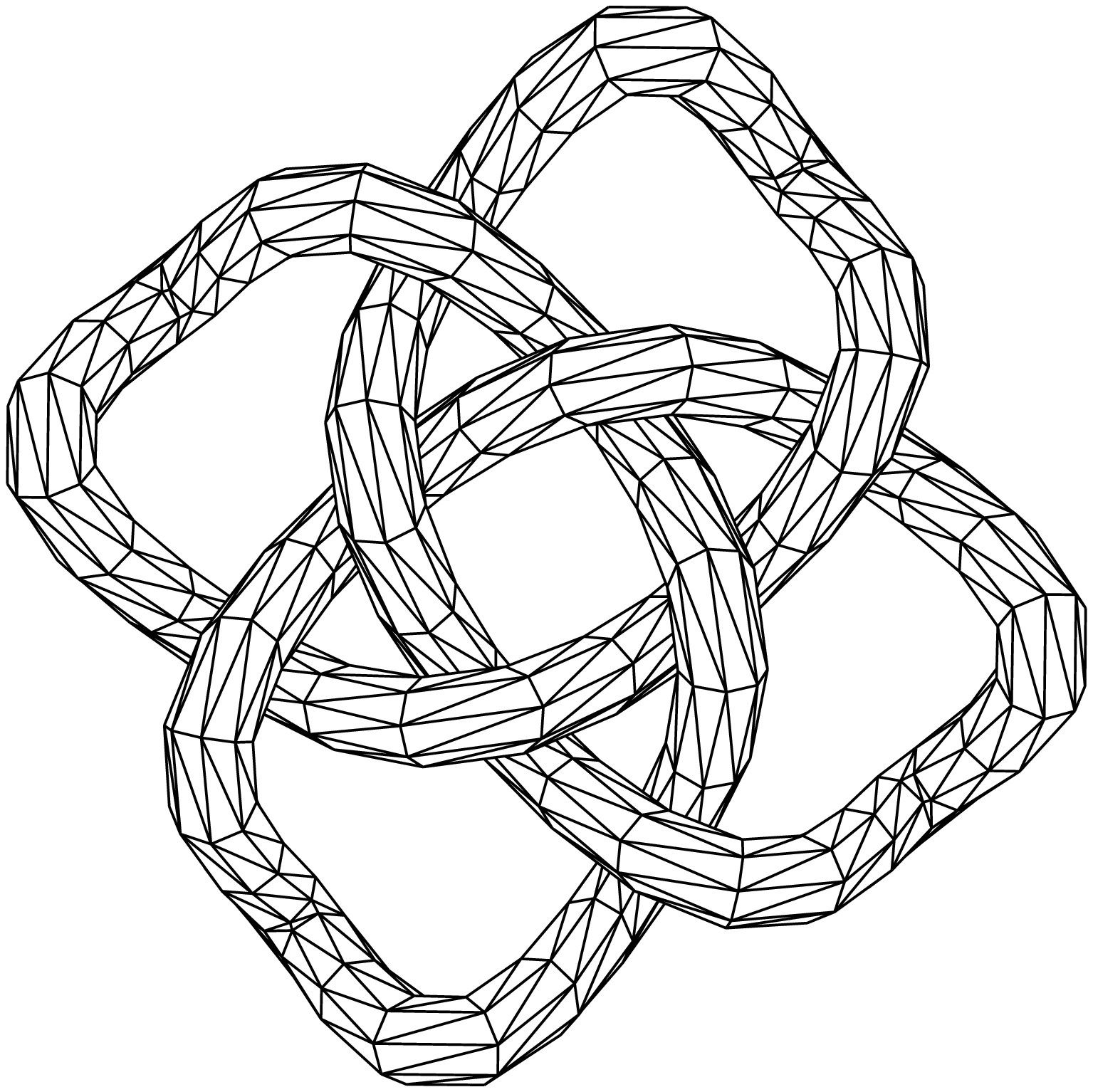}} \quad \quad  \scalebox{0.3}{\includegraphics{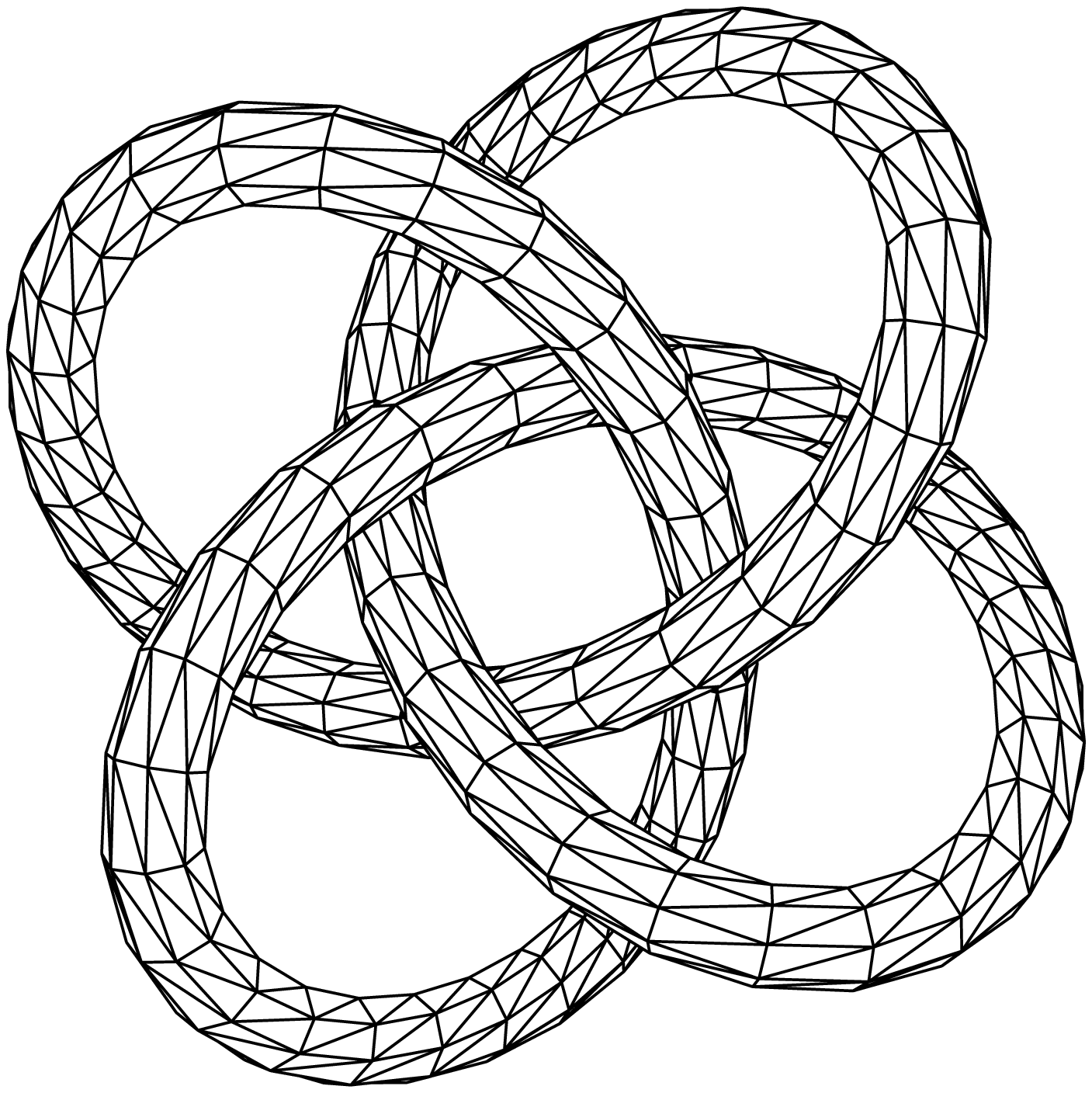}}

\quad \quad \caption{Alternating knot. \quad \quad Non-alternating knot.}
\end{figure}

This number, for many knots, is also equal to the sum of the absolute value of the coefficients of the Jones polynomial \cite{jo}. Tables of the
coefficients of the Alexander and Jones polynomials were tabulated for the simplest prime knots in appendix II of the book by
K. Murasugi \cite{mu}. Note that the previous assertion is not true for the non-alternating knot $8_{19}$ in this table.

Conway's number $C$ also appears in the group theory of the complement of a knot when one has the Wirtinger presentation \cite{ro} by means
of two generators and an
algebraic relation among the generators and its inverses with $2 C$ factors equal to the identity. In many cases, the relation is presented as
$C$ factors
of a palindromic product equal to $C$ factors, where the generators have been interchanged. For example, for the figure eight knot represented in Fig. 2,
with Conway's number equal to 5, the group of the complement of the knot is presented \cite{ro} as
$\{x_1, x_3 | \, x_1 x_3 x_1^{-1} x_3 x_1 = x_3 x_1 x_3^{-1} x_1 x_3 \}$,
where $x_1$ and $x_3$ are two generators of the group. An analogous presence of Conway's number is found in the Wirtinger presentation of many groups
of the complement of knots with two generators and one relation.

The same Conway's number was found multiplying the number of crossings of any knot equal to the derivative for $x = 2$ of the characteristic polynomial
of the adjacency matrix associated to the graph of alternating knots oriented as in Fig. 2 \cite{pi}.

One of the main objectives of this paper is an alternative definition of Conway's number independent of the previous numbers.
This opens the requirement of a formal proof of the identification of these numbers, which is not developed in this paper.

\begin{figure}

\hfil\scalebox{0.5}{\includegraphics{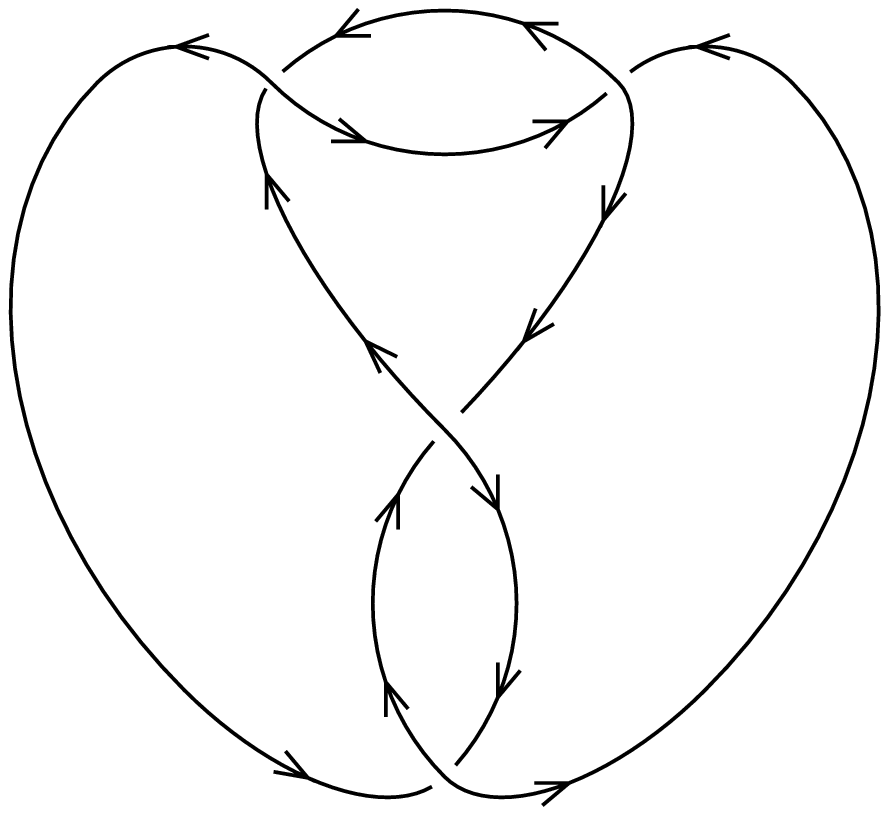}}\hfil

\caption{Figure eight knot with four superpositions showing the orientation given to the strands (edges). Its Conway number is 5, its conways are 2 and 2,
related by $5/2 = 2 + 1/2$.}

\

\hfil\scalebox{0.7}{\includegraphics{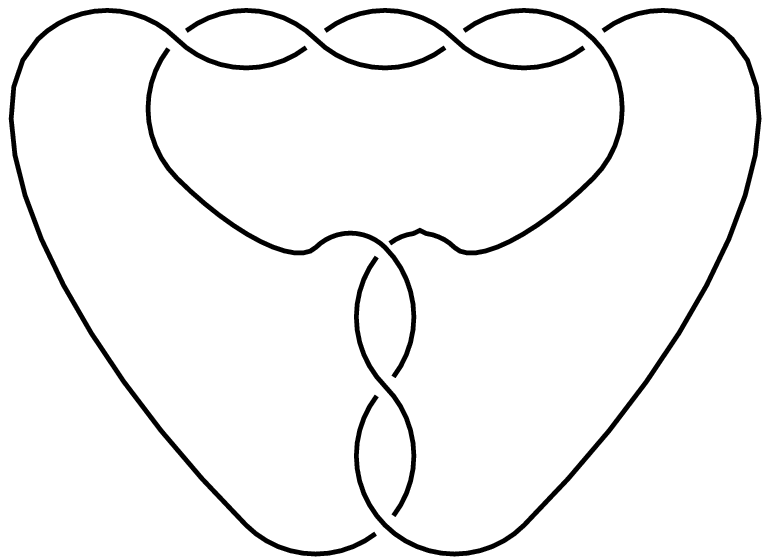}}\hfil

\

\caption{The knot $7_3$ formed by two conways equal to 3 and 4. Its Conway's number is 13. One has the relations $13/4 = 3 + 1/4$
and $13/3 = 4 + 1/3$. According to this drawing, the conway with 4 overpassings is right handed and the conway with 3 crossings is left handed.}

\

\psfrag{p}{\LARGE{4}}
\psfrag{q}{\LARGE{3}}
\hfil\scalebox{0.5}{\includegraphics{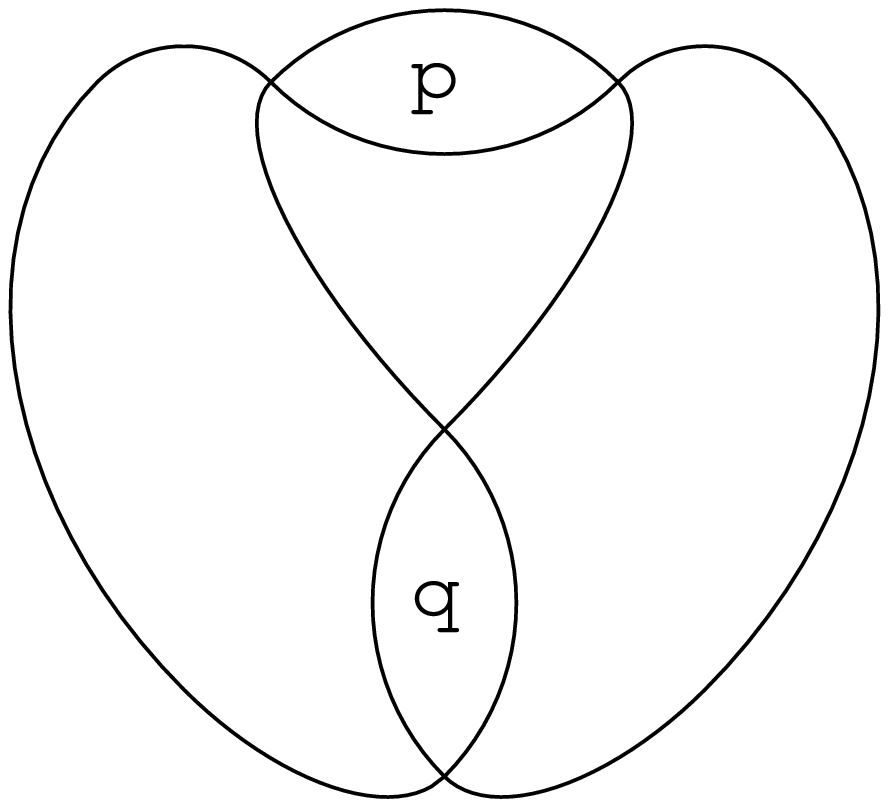}}\hfil

\

\caption{The same knot as in Fig. 3 with compact notation.}

\

\end{figure}

The rational knot characterization has been documented by Kauffman and Lambropoulou \cite{kl}. Perhaps the best definition of rational knot is by
the induction generated by the following Figs. 6-9.

For most of the knots tabulated by Conway \cite{co}, he found a subclass of alternating knots called \textit{rational knots}. He associated to these knots
a ratio of two natural numbers, where the numerator is the Conway's number of the knot and the denominator is Conway's number of another related
knot. When the ratio of these two numbers is developed in a continued fraction he found a set of integers, associated to the knot, called in the present paper
{\sl conways}, corresponding to the number of successive crossings, which may be separated from the rest of the knot by four strands. See Figs. 2, 3 and
4 to illustrate the concept of a conway. In this paper a conway is the geometric part of the knot and conway is the number of crossings in this conway.
The rational knots were discovered by Conway by providing the relation to the ratio of those two Conway's numbers.

As a number, the conway is the largest number allowed by the topology of the knot. This digression is necessary to take into account the {\sl flype}
transformation of a knot discovered by Tait \cite{ta} that is capable of reducing the number of conways of a knot to a minimum. Two knots are equivalent
(isotopic) by a flype. A flype is the modification of the representation of a knot, or a family of knots, when any 2-tangle, represented by the square,
is connected in the knot to a two conways equal to $a_1$ and $a_2$ ($a_1, a_2 \ \epsilon \ \mathbb{N} \cup 0$) as shown in Fig. 5. When rotating
the square tangle in the direction of handedness of the two conways, one conway loses one crossing by untwisting it, while the other conway wins one
crossing by twisting it
$$
a_1 \longrightarrow a_1 + 1\, \quad a_2 \longrightarrow a_2 - 1 \, .
$$
By succesive rotations (flypes) one of the conways may disappear and the other conway adds the numbers of the two conways and becomes just one conway.
$$
a_1 \longrightarrow a_1 + a_2 \, \quad a_2 \longrightarrow 0 \, .
$$
One assumes that this operation has always been performed until a minimal number of conways is obtained for any knot and for any family of knots.
\begin{figure}

\psfrag{a}{\LARGE{$a_1$}}
\psfrag{b}{\LARGE{$a_2$}}
\hfil\scalebox{0.7}{\includegraphics{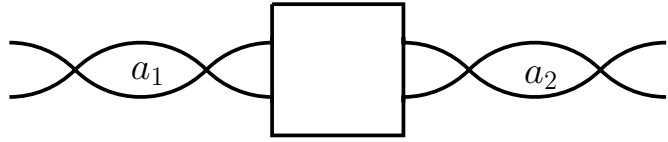}}\hfil

\

\caption{Flypes are needed to form a single conway equal to $a_1 + a_2$}
\end{figure}

\subsection*{Definition 1.} A family of alternating knots is formed by conways.

\subsection*{Definition 2.} A conway is the number of $a_j$ succesive crossings in the projection of a knot. Represent the twisting of two strands
with a definite handedness left or right $a_j$ times. In general the conway is a positive number, represented algebraically by the letter $a_j$
with an integer subindex $j$ which labels the conway. In the graphic representations of the knot the conway is represented by an "eye" connected
at the commissures with two strands on each side and the number of crossings $a_j$ in the apple of the eye as indicated in Fig. 4, representing the
knot in Fig. 3. When the conway is equal to 1 the handedness of the conway is not defined, is neither right nor left handed. When the conway
is equal to 0 the two strands are separated with no twisting, defining the orientation (although without crossing.)

Conways formed by two or more crossings are of two different orientations either left or right-handed, which are similar to the handedness
of a mechanical screw. A corkscrew is generally right-handed. The two handednesses are interchanged upon reflection in a mirror. Fig. 2 is a knot
formed by two conways of two crossings each. The upper conway is right-handed in the horizontal position and the lower conway is
left-handed in the vertical position. Something similar happens for the knot in Fig. 3; the conway with 4 crossings is right-handed,
the conway with 3 crossings is left-handed. The handedness of the conways of an alternating knot is determined by the handedness of
any of its conways. Its mirror image has the opposite handedness for all of its conways. To deal with the characterization of the handedness
of a particular knot, or family of knots, the orientation of the strands of the knot was chosen as in Fig. 2: the strands leaving an overpassing
 go in opposite directions, the strands going through an underpassing point at the crossing point. With this orientation of the strands in
Figs. 2 and 3, the "eyes" of the right-handed conway have a counter-clockwise orientation; the "eyes" of the left-handed conway have a
clockwise orientation.

An equivalent orientation which seems more efficient to me is obtained by the coloring of the knots and families of knots, using white and other colors
in the form called of the checkerboard. White regions are opposed at each crossing. Colored regions are opposed at each crossing. We choose white
as the exterior color. Then white conways are chosen to be left-handed and colored conways are then right-handed. Figs. 6-9 show the checkerboard coloring
of rational knots. Figs. 27-52 for the 44 families of 6 conways are another examples of white and checkerboard coloration (with two or three colors
to represent the factorization of a knot in tangles).

\subsection*{Definition 3.} The seed of a family of knots is a knot or a link corresponding to the family when all the conways are equal to 1.
The family is
formed substituting each crossing of the seed by a conway. A seed corresponds in general to several families because at each crossing of the
seed where no handedness is defined, one introduces a conway with two different possibilities of right and left handedness that in many cases
produces different families. The strands forming the seed are inherited to all the families corresponding to the seed, which besides have the
strands forming the conways added to the seed. If the seed has $n$ crossings, the number of strands of the seed is equal to $2 n$.

\begin{figure}

\psfrag{a}{\LARGE{$a_1$}}
\psfrag{b}{\LARGE{$a_2$}}
\psfrag{c}{\LARGE{$a_3$}}
\psfrag{d}{\LARGE{$a_4$}}
\psfrag{e}{\LARGE{$a_5$}}

\hfil\scalebox{0.5}{\includegraphics{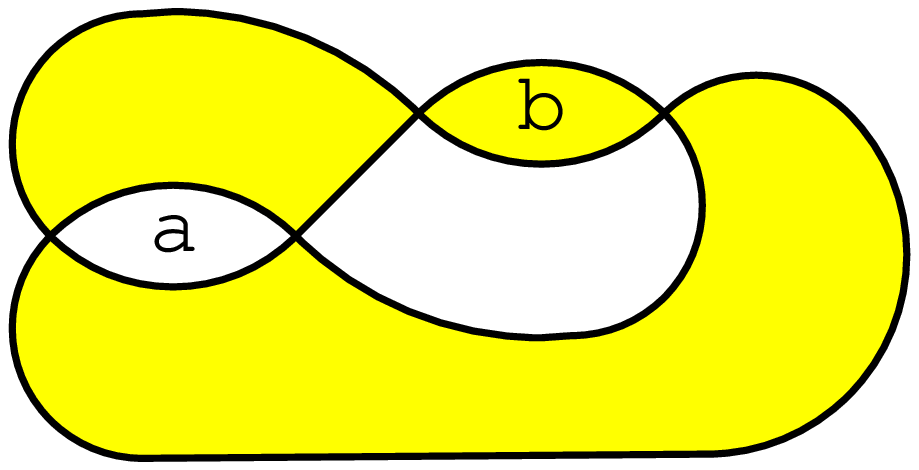}}\hfil

\caption{Family of rational knots with 2 conways}

\

\hfil\scalebox{0.5}{\includegraphics{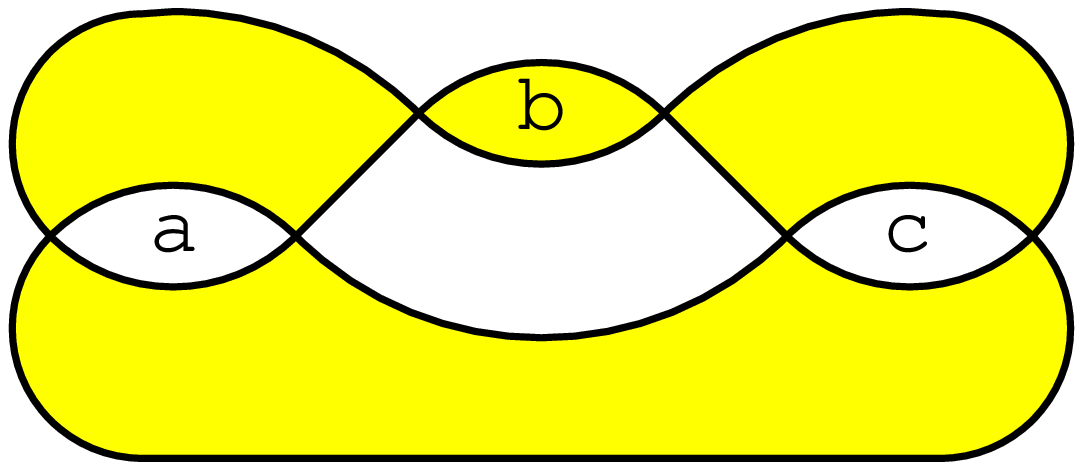}}\hfil

\caption{Family of rational knots with 3 conways}

\

\hfil\scalebox{0.5}{\includegraphics{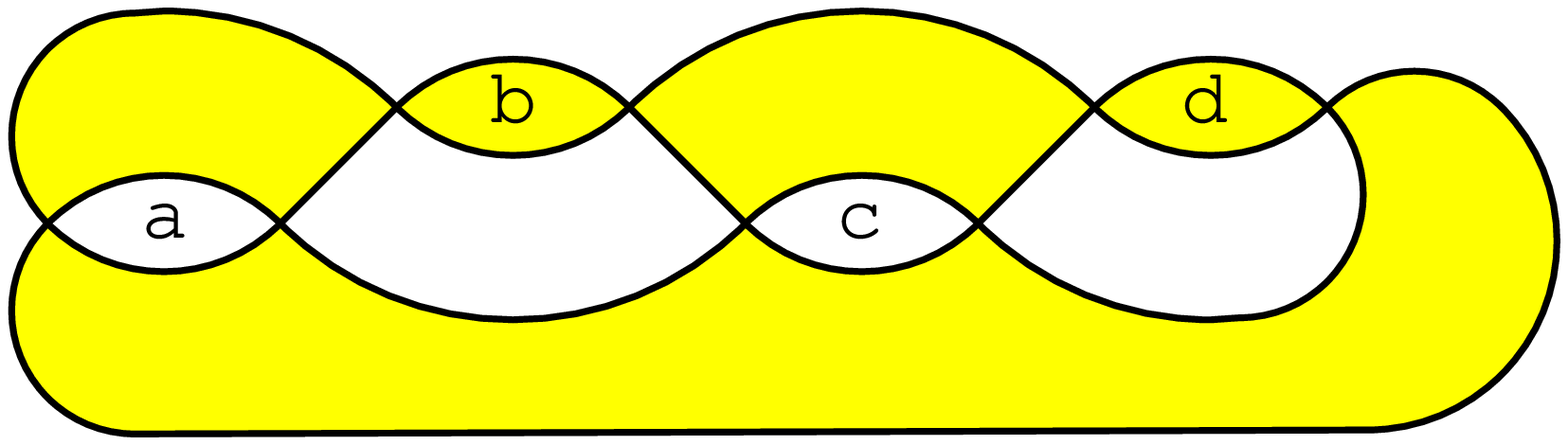}}\hfil

\caption{Family of rational knots with 4 conways}

\

\hfil\scalebox{0.5}{\includegraphics{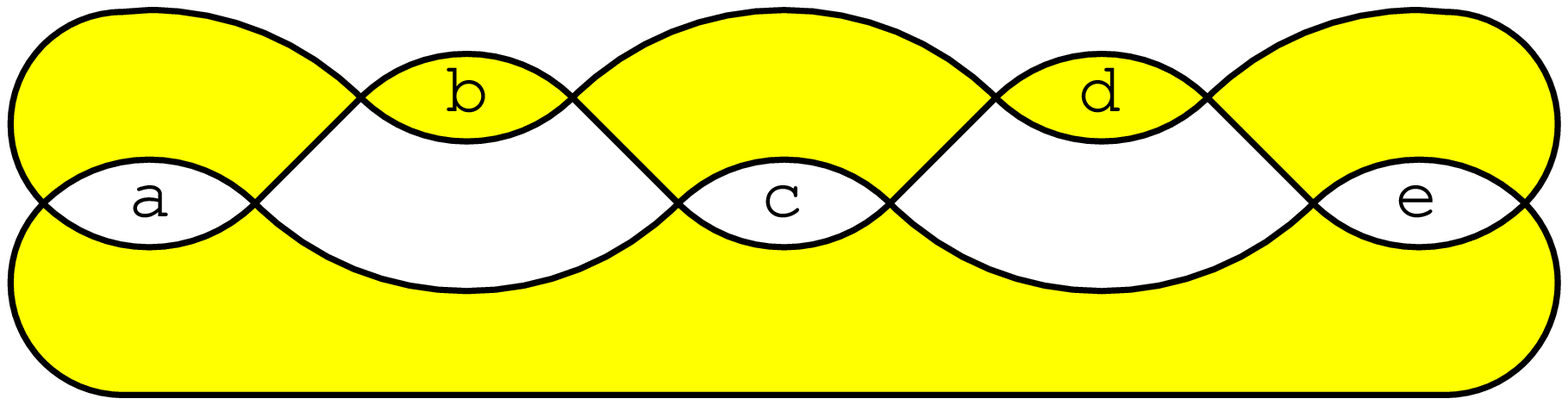}}\hfil

\caption{Family of rational knots with 5 conways}

\end{figure}

Conway \cite{co} introduced in his table of prime knots the conways (numbers) of each knot. For rational knots the conways (numbers) are written
successively by him
with no separation. For the non rational knots he separated the numbers by semicolons or used other set of symbols. The association of conways to
knots is reproduced in the tables of knots in D. Rolfsen's book \cite{ro} and in the Appendix of the C. Adams' book \cite{ad}. In both, Conway's paper \cite{co} is cited to explain the meaning of the conways (numbers) of the knot.

 This paper is an extension of a previous study \cite{pi} of Conway's functions of families of alternating knots with up to five conways. In the present paper all the
families of knots with up to six conways are drawn and with each figure the corresponding Conway's function of each family of knots is associated, which is defined in the next section. The principal tool to be presented and used is a linear algebra related to the the separation
of the knot in {\sl 2-tangles}. The simplest family of 2-tangles is one conway. Any family of alternating knots and any family of alternating 2-tangles
is the union of conways. The concept of tangle is a fundamental concept of the theory of knots (see for example definitions and properties of tangles in
books \cite{mu} and \cite{ad}). The Kauffman and Lambropoulou work \cite{kl} on rational knots and tangles was strongly influential on this and the previous work \cite{pi}. Their work based on the continued fraction decomposition of two integers is here generalized into a linear algebra of alternating
(not just rational) knots and tangles. Nevertheless the first axiom of our algebra is: The algebra here defined applies without changes to the rational case studied and proved by Kauffman and Lambropoulou \cite{kl} which is a subset of the set of alternating knots.

\subsection*{Axiom 1} All the algebraic properties of the families of rational knots and 2-tangles are particular cases of the algebraic properties of families of alternating knots and 2-tangles.

\section{CONWAY'S FUNCTION OF A FAMILY OF ALTERNATING KNOTS}
In this work to each family of alternating knots formed by $n$ conways, equal to the positive numbers $a_1, a_2, ... , a_n$, a function
of these conways is associated. This function of the family, called Conway's function (C-function), computed for the particular choice of the values of the conways of a knot of the
family, is equal to the Conway's number used by him in his table \cite{co}. Rolfsen \cite{ro} calls it the determinant of the knot since it is the determinant
of the symmetric part of a Seifert matrix.

\subsection*{Axiom 2.} The C-function of a family of knots is defined on the natural numbers $\mathbb{N} \cup 0$, function of conways. It is a polynomial, formed by  monomials each a sum of products of conways, linear in
all the conways, with the possibility of monomials not containing some conways and sometimes a term equal to 1. The coefficient of each monomial is $+1$. In general it is formed by two terms,
one having a particular conway as a common factor and other term without this conway. This separation is possible for any selection of the particular
family and for any choice of the conway of that family.

\subsection*{Axiom 3.} For two or more separated knots the C-function is 0.

\subsection*{Axiom 4.} For composite knots the C-function is equal to the product of the C-functions of the composite knots. A particular definition of composite knots is presented when defining 1-tangles.

\subsection*{Axiom 5.} The C-function of a knot formed by just one conway is equal to the number of crossings forming the conway. Important particular cases are
the case that equals 1, for the unknotted circle (with or without a simple twist), and the case that equals 0, corresponding to two unlinked circles.

\subsection*{Axiom 6.} The number of monomials of the C-function is equal to the Conway's number of the corresponding seed.

\subsection*{Axiom 7.} For a particular family of knots all the number of factors (conways) in each monomial is an even number or an odd number. Number 1,
corresponding to zero $a_j$-factors is considered an even number of factors.

\subsection*{Axiom 8.} The number of unknotted twists do not appear in the C-function. Adding twists to a family of knots does not change the C-function.

To construct the families of $n$ conways corresponding to one seed of $n$ crossings one places one conway with two different possible orientations at each crossing of the seed. Suppress the cases when two conways could be replaced by a single conway, the cases which are only different by a relabeling of the conways, and the cases which are equivalent by a flype.

The C-function of a given family of $n$ conways is constructed by induction as follows. Assume that the C-function of all the families of $n-1$ conways or $n-2$ conways is known. Starting from the graphic representation of the family of $n$ conways, successively the number of crossings of a particular conway of the family is decreased to zero, corresponding to a family with $n-1$ or $n-2$ conways where the chosen conway disappears joining the regions with the same color that the color of the conway. These two regions which were separated by a crossing in the seed now form a single region. The resulting family of $n-1$ or $n-2$ conways has a C-function such that all the monomials are also monomials of the considered family of $n$ conways. In most cases the C-function of the family of $n$ conways is obtained by the superposition of the monomials obtained by this procedure. There exists a term, not produced in this way, corresponding to the case when a term of the C-function is equal to the product of all the conways: $\Pi_j a_j$. This case is identified since the total number of terms in the C-function is equal to the Conway's number of the seed. Using this method the C-function was computed in \cite{pi} for families of knots from one to five conways.

Sometimes a C-function of an assumed family of $n$ conways depends on two conways as the sum of them. This case corresponds to a family with two conways
of the same color connected by the commissures, corresponding to a single conway, therefore to a family with less than $n$ conways which should be
supressed from the full collection of families of $n$ conways. Frequently, to recognize we are in this case one should use flypes in order to have two
contiguous conways.

The number of families increases with the number of conways. For one or two conways we have only 1 family. For three conways we have 2 families. For
four conways we have 5 families. For five conways there are 12 families. For six conways I find 44 families.

The seed of the family of $n$ conways has $2 n$ strands, which in the family connect these conways in all the different possible handednesses.

The number of monomials in the C-functions of rational knots is equal to the Conway's number of the corresponding seed, which is a Fibonacci number.

Although I do not know if gaussian brackets have been used in knot theory, it is obvious to me that they should be used for rational knots and tangles as I do.
Gaussian brackets, discovered by Euler \cite{eu} and Gauss \cite{ga} and used frequently in geometric optics \cite{gh} are defined by
\begin{equation}
\begin{array}{l}
G[ \ ] = 1 \\
G[a_1] = a_1 \\
G[a_1, a_2] = a_2 G[a_1] + G[ \ ] = 1 + a_1 a_2 \\
G[a_1, a_2, ... , a_n] = a_n G[a_1, a_2, ... ,a_{n-1}] + G[a_1, ... , a_{n-2}]
\end{array}\, ,
\end{equation}

\subsection*{Axiom 9.} The families of rational knots have as C-function the gaussian bracket of the conways of the family.

The basis to this axiom is that Gauss finds \cite{ga} that the value of a continued fraction is the ratio of two gaussian brackets
\begin{equation}
\frac{G[a_1, a_2, ... , a_n]}{G[a_2, a_3, ... , a_n]} = a_1 + \frac{1}{a_2 + \frac{1}{a_3 + ... +\frac{1}{a_n}}}\, .
\end{equation}
Therefore the Conway's function of the circle is 1, The Conway's function of the torus knots formed by one conway is $a_1$ and the Conway's function
of the family formed by two conways is $a_1 a_2 + 1$. Other Conway's functions of families of rational knots are:
\begin{equation}
\begin{array}{l}
G[a_1, a_2, a_3] = a_3 G[a_1, a_2] + G[a_1] = a_1 a_2 a_3 + a_1 + a_3 \\
G[a_1, a_2, a_3, a_4] = a_4 G[a_1, a_2, a_3] + G[a_1, a_2] =\\
(a_1 a_2 + 1)(a_3 a_4 + 1) + a_1 a_4  \\
G[a_1, a_2, a_3, a_4, a_5] = a_5 G[a_1, a_2, a_3, a_4] + G[a_1, a_2, a_3] =\\
a_1 + a_5 + a_1 (a_2 + a_4) a_5 + (a_1 a_2 + 1) a_3 (a_4 a_5 + 1) \\
G[a_1, a_2, a_3, a_4, a_5, a_6] = a_6 G[a_1, a_2, a_3, a_4, a_5] + G[a_1, a_2, a_3, a_4] = \\
(a_1 a_2 a_3 + a_1 + a_3)(a_4 a_5 a_6 + a_4 + a_6) + (a_1 a_2 + 1)(a_5 a_6 + 1)
\end{array} \, .
\end{equation}

For 2, 3, 4, 5, conways the rational families of knots were represented in figures 6-9. In these figures the handedness of the conways was distinguished by different color. One notes that the successive conways in a rational family of knots have different handedness.

\subsection*{Theorem 1.} Figures 6-9 have C-functions equal to gaussian brackets.

\noindent {\bf Proof:} The proof of the correspondence of the gaussian bracket with the C-function of the rational knots appearing in Figs. 6-9 comes from doing 0 the $a_1$ conway in those figures: this amounts to opening both commissures of the conway $a_1$ in the figures. Then the $a_2$ conway becomes a simple twist which disappears also from the family. Fig. 6 becomes a circle with C-function equal to $G[ \ ] = 1$. Fig. 7 becomes a torus family of one conway equal to $a_3$ and C-function equal to $G[a_3] = a_3$. Fig. 8 becomes the rational knot of Fig. 6 with conways $a_3$ and $a_4$ and C-function $G[a_3, a_4] = a_3 a_4 + 1$. Fig. 9 becomes the rational knot of Fig. 7 with C-function $G[a_3, a_4, a_5] = a_3 a_4 a_5 + a_3 + a_5$. This is exactly the definition of the gaussian brackets. $\square$

For 1, 2, 3, 4, 5, 6 conways the families of rational knots were represented again by Figs. 13, 14, 20, first of 26 and first of 42 in the appendix, respectively.

In geometric optics the gaussian brackets form part of a linear algebra with their main property based on the factorization with matrices.

\subsection*{Theorem 2.}
The product of the following matrices is expressed in terms of gaussian brackets as
$$
\left( \begin{array}{cc}
a_1 & 1 \\
1 & 0
\end{array} \right) \left( \begin{array}{cc}
a_2 & 1 \\
1 & 0
\end{array} \right) \left( \begin{array}{cc}
a_3 & 1 \\
1 & 0
\end{array} \right) . . . \left( \begin{array}{cc}
a_n & 1 \\
1 & 0
\end{array} \right) =
$$
\begin{equation}
\left( \begin{array}{cc}
G[a_1, a_2, a_3, ..., a_n] & G[a_1, a_2, ... a_{n-1}] \\
G[a_2, a_3, ..., a_n] & G[a_2, a_3, ... a_{n-1}]
\end{array} \right) \, .
\end{equation}

\noindent {\bf Proof:} The proof of this property follows by induction from the gaussian bracket definition (1). $\square$

This property is used in this paper to express the Conway's functions associated to rational knots (gaussian brackets) as a product of a linear
algebra of vectors and matrices whose entries are gaussian brackets corresponding to rational tangles and other kinds of non rational tangles. Different
factorizations corresponding to different dissections in tangles to compose the same knot.

The C-function of a family of rational knots of $n$ conways appears in one component of the previous product of matrices and is obtained by the projection
$$
G[a_1, a_2, a_3, ..., a_n] =
$$
\begin{equation}
\left( \begin{array}{cc}
1 & 0
\end{array} \right) \left( \begin{array}{cc}
a_1 & 1 \\
1 & 0
\end{array} \right) \left( \begin{array}{cc}
a_2 & 1 \\
1 & 0
\end{array} \right) \left( \begin{array}{cc}
a_3 & 1 \\
1 & 0
\end{array} \right) . . . \left( \begin{array}{cc}
a_n & 1 \\
1 & 0
\end{array} \right) \left( \begin{array}{c}
1 \\
0
\end{array} \right)\, .
\end{equation}

Making $a_1 = 0$ in the previous product one has a linear algebra representation of the proof of Theorem 1
\begin{equation}
\left( \begin{array}{cc}
1 & 0
\end{array} \right) \left( \begin{array}{cc}
0 & 1 \\
1 & 0
\end{array} \right) \left( \begin{array}{cc}
a_2 & 1 \\
1 & 0
\end{array} \right) = \left( \begin{array}{cc}
1 & 0
\end{array} \right) \, . \quad \square
\end{equation}

\subsection*{Theorem 3.}
The C-function of the family of rational knots has a palindromic symmetry.

\noindent {\bf Proof.} The proof follows from the fact that one has the representation (5) of the C-function as a product of matrices which are symmetric matrices and the elementary
property of the transpose matrix of a product of matrices. $\square$

This property is well known for the decomposition of a ratio of two numbers in continued fractions and for gaussian brackets.

\subsection*{Corollary 1}
From theorems 2 and 3, making $a_n = 0$ the last conway produces a similar property at the other end of the rational knot
\begin{equation}
\left( \begin{array}{cc}
a_{n-1} & 1 \\
1 & 0
\end{array} \right) \left( \begin{array}{cc}
0 & 1 \\
1 & 0
\end{array} \right) \left( \begin{array}{c}
1 \\
0
\end{array} \right) = \left( \begin{array}{c}
1 \\
0
\end{array} \right)
\end{equation}

\subsection*{Theorem 4.}
Making $a_j = 0$ a conway in the interior of a rational knot the two neighbor conways join in a single conway that equals $a_{j-1} + a_{j+1}$.

\noindent {\bf Proof.} One has two proofs: a topological proof that comes from figures 7-9. When one opens the two commissures of the eye whose $a_j$ is going
to be zero in those figures the two neighbor conways are placed one after the other forming a single conway. A corresponding algebraic
proof is given by the remarkable equation
\begin{equation} \label{parabolica}
\left( \begin{array}{cc}
a_{j-1} & 1 \\
1 & 0
\end{array} \right) \left( \begin{array}{cc}
0 & 1 \\
1 & 0
\end{array} \right) \left( \begin{array}{cc}
a_{j+1} & 1 \\
1 & 0
\end{array} \right) = \left( \begin{array}{cc}
a_{j-1} + a_{j+1} & 1 \\
1 & 0
\end{array} \right) \, . \quad \square
\end{equation}
These theorems give information on the families of rational knots. Equation (\ref{parabolica}) introduces algebraically the essential property that
two contiguous conways (by the commissures) are equivalent to a single conway. If the matrices in (\ref{parabolica}) represent the conways, their
joining to form a single conway is glued by using the matrix
\begin{equation}
M = \left( \begin{array}{cc}
0 & 1 \\
1 & 0
\end{array} \right)\, ,
\end{equation}
which acts as a metric matrix. This matrix play an essential role in the following algebra of tangles that are represented by products of vectors
and matrices and uses this matrix to define the interior product between them.

The advantage of linear algebra with respect to continued fractions is that division by zero is not defined {\sl a priori} for continued fractions, whereas doing some $a_j$ equal to zero in the gaussian bracket or in a matrix is well defined.

\begin{figure}

\psfrag{a}{\LARGE{$a_1$}}
\psfrag{b}{\LARGE{$a_2$}}

\hfil\scalebox{0.45}{\includegraphics{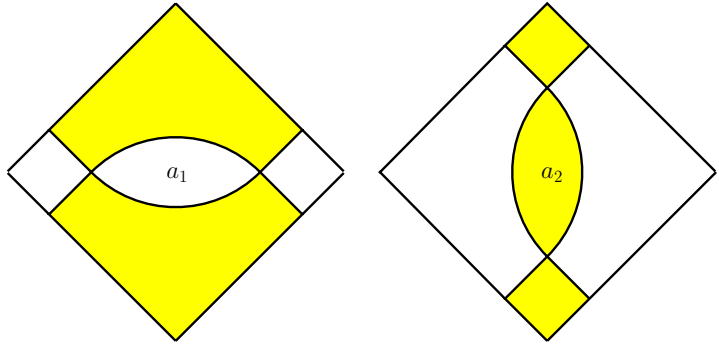}}\hfil

\caption{Rational tangles formed by one conway. Its covariant and contravariant components are $(1 , a_1)$ and $(a_2, 1)$}
\end{figure}

\begin{figure}

\psfrag{a}{\LARGE{$a_1$}}
\psfrag{b}{\LARGE{$a_2$}}
\psfrag{d}{\LARGE{$a_3$}}
\psfrag{c}{\LARGE{$a_4$}}

\hfil\scalebox{0.35}{\includegraphics{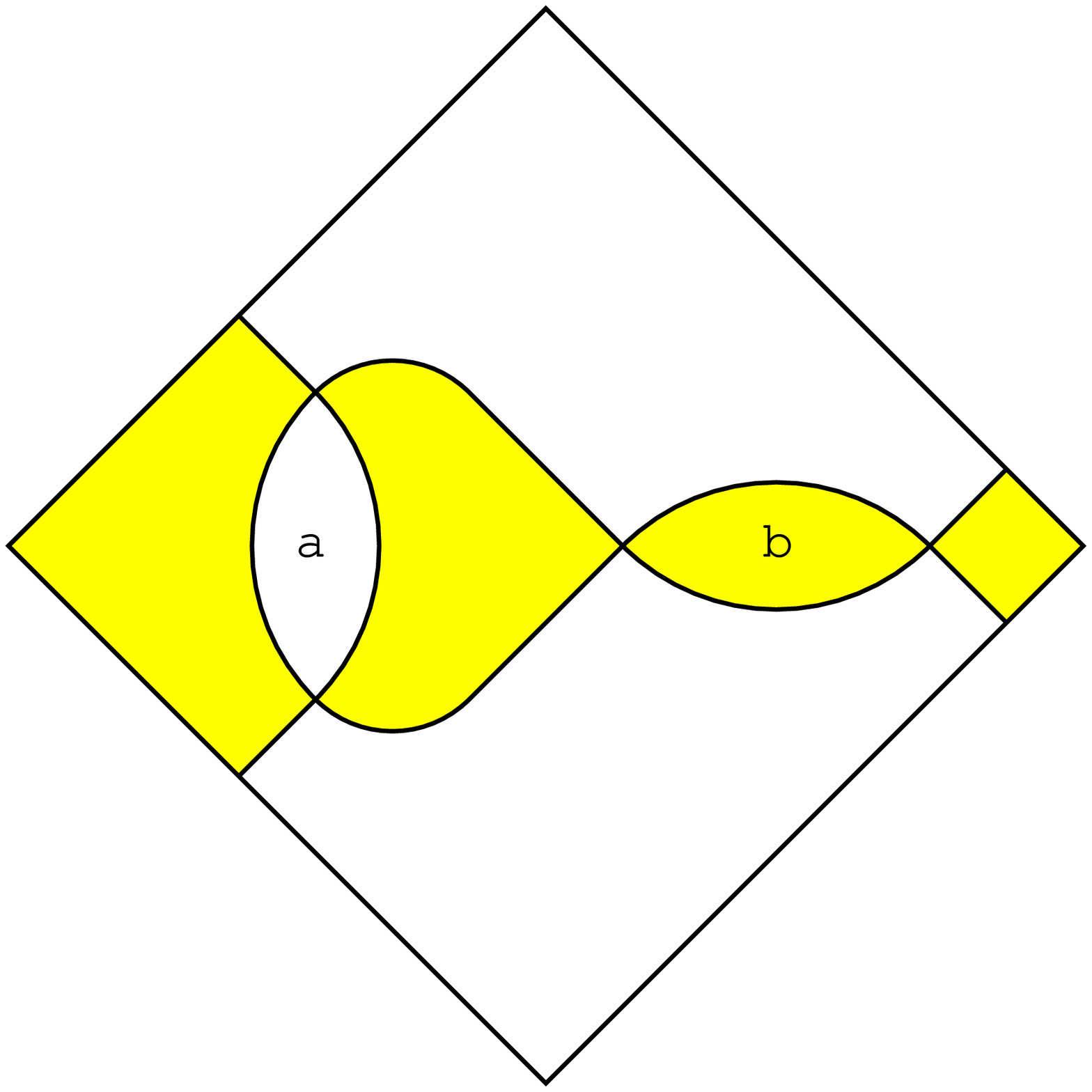}} \quad \scalebox{0.35}{\includegraphics{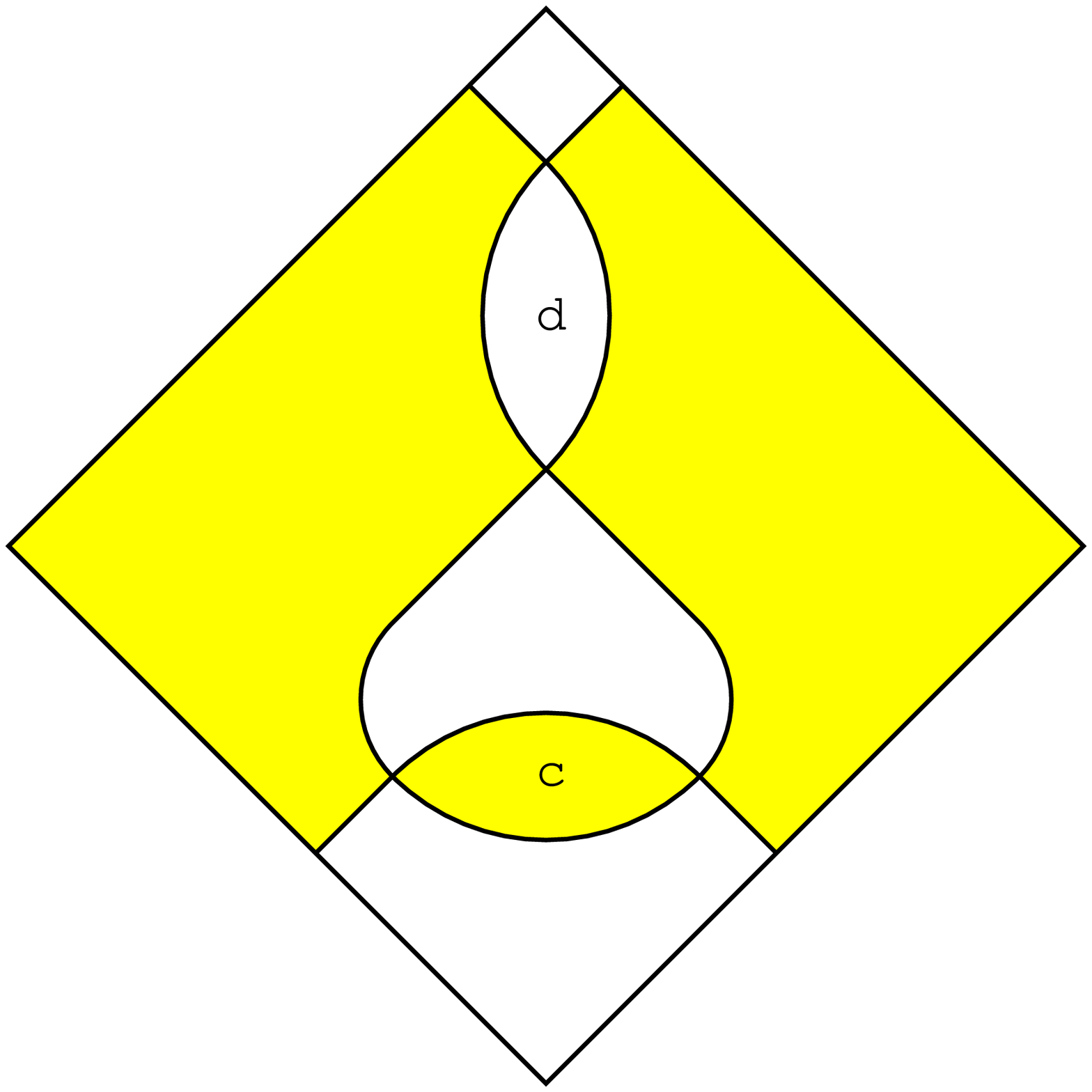}}\hfil

\caption{Rational tangles formed by two conways. Its covariant and contravariant components are $(1 + a_1 a_2, a_1)$ and $(a_4, 1 + a_3 a_4)$.}

\end{figure}

\section{ALGEBRA OF FAMILIES OF RATIONAL TANGLES}
We introduce our concept of \textit{family of tangles}. We cut a knot family in the interior of a strand connecting two conways. The two cutted extremes
are fixed to the boundary of a 3-sphere containing the knot. This define a 1-tangle. If we connect two 1-tangles joining the cut end we find a
composite knot. The opposite operation, to separate a knot by cutting only two strands is possible for composite knots. When a knot can not be separated
in two as in composite knots, it is called a prime knot. A composite knot can be decomposed in two or more knots.

In what follows, we are mainly interested in families of prime knots. Composite knots are considered only if it is unavoidable.

From a family of alternating knots we construct families of alternating 2-tangles by cutting two different strands connecting
conways, corresponding to the seed of the
family. If a family is separated in two parts by cutting four strands, we call a 2-tangle to each part. These strands are fixed to four points in the
boundary of a ball containing the knot. This fixedness precludes the unknotting. Each 2-tangle gives two possible knots by connecting the four strands in
the two possible forms allowed by the alternating restriction.

In this section an inductive approach is used to obtain an algebra of families of tangles constructed by using the experience on families of rational knots.

Returning to Figs. 6 and 8 of rational knots, cut these two families at the four strands in the middle of the two families. For figure 6 one obtains
two similar tangles formed by one conway with different handedness (recognized by the different color). These two tangles are recognized in Fig. 10. The
factorization (5) for this case gives the equation for the C-function of a rational knot
\begin{equation}
G[a_1, a_2] = \left( \begin{array}{cc}
1 & 0
\end{array} \right) \left( \begin{array}{cc}
a_1 & 1 \\
1 & 0
\end{array} \right) \left( \begin{array}{cc}
a_2 & 1 \\
1 & 0
\end{array} \right) \left( \begin{array}{c}
1 \\
0
\end{array} \right) = \left( \begin{array}{cc}
a_1 & 1
\end{array} \right) \left( \begin{array}{c}
a_2 \\
1
\end{array} \right)\, ,
\end{equation}
as the product of two similar vectors. Notice that the two parts have different handedness, that the gluing of conways in (8) was using the
metric matrix $M$ defined in (9) and that the square of matrix $M$ is the unit matrix. The previous equation may be rewritten as
\begin{equation}
G[a_1, a_2] = \left( \begin{array}{cc}
a_1 & 1
\end{array} \right) M M\left( \begin{array}{c}
a_2 \\
1
\end{array} \right) = \left( \begin{array}{cc}
1 & a_1
\end{array} \right) M \left( \begin{array}{c}
a_2 \\
1
\end{array} \right)\, .
\end{equation}
To the 2-tangle formed by one white conway one associates the (covariant) vector $(1, a)$, where $a$ is the number of crossings or overpassings forming it.
To the corresponding colored conway one associates the (contravariant) vector $(a, 1)$. The two components transform one into the other by the metric
matrix M. The interior product between this two vectors (white and colored) is by means of the $M$ metric matrix.

A similar dissection in the middle of Fig. 8 corresponding to the family of rational knots formed by four conways leads to the factorization of
the C-function of the family in the interior product of two vectors corresponding to the tangle in Fig. 11, where different coloration distinguishes
the covariant  $(a_1 a_2 + 1, a_1)$ and contravariant components $(a_4, a_3 a_4 + 1)$ of the tangles factoring the C-function
\begin{equation}
G[a_1, a_2, a_3, a_4] = \left( \begin{array}{cc}
a_1 a_2 + 1 & a_1
\end{array} \right) M \left( \begin{array}{c}
a_4 \\
a_3 a_4 + 1
\end{array} \right)\, .
\end{equation}

Noting the factorizations
$$
\left( \begin{array}{cc}
a_1 a_2 + 1 & a_1
\end{array} \right) = \left( \begin{array}{cc}
1 & a_1
\end{array} \right) M \left( \begin{array}{cc}
a_2 & 1 \\
1 & 0
\end{array} \right) \, ,
$$
and
$$
 \left( \begin{array}{c}
a_4 \\
a_3 a_4 + 1
\end{array} \right) = \left( \begin{array}{cc}
0 & 1 \\
1 & a_3
\end{array} \right) M \left( \begin{array}{c}
a_4 \\
1
\end{array} \right)\, .
$$
one deduces other factorizations for the C-function
$$
G[a_1, a_2, a_3, a_4] =
$$
$$
\left( \begin{array}{cc}
1 & a_1
\end{array} \right) M \left( \begin{array}{cc}
a_2 & 1 \\
1 & 0
\end{array} \right) M \left( \begin{array}{cc}
0 & 1 \\
1 & a_3
\end{array} \right) M \left( \begin{array}{c}
a_4 \\
1
\end{array} \right) =
$$
$$
\left( \begin{array}{cc}
1 & a_1
\end{array} \right) M \left( \begin{array}{c}
a_2 a_3 a_4 + a_2 + a_4\\
a_3 a_4 + 1
\end{array} \right) =
$$
$$
\left( \begin{array}{cc}
a_1 a_2 + 1 & a_1 a_2 a_3 + a_1 + a_3
\end{array} \right) M \left( \begin{array}{c}
a_4 \\
1
\end{array} \right)\, .
$$

These equivalent factorizations are used to transform the expression (5) by adding identity factors $M^2$ into
$$
G[a_1, a_2, a_3, ..., a_{2 k}] =
$$
\begin{equation}
\left( \begin{array}{cc}
0 & 1
\end{array} \right) M \left( \begin{array}{cc}
0 & 1 \\
1 & a_1
\end{array} \right) M \left( \begin{array}{cc}
a_2 & 1 \\
1 & 0
\end{array} \right) M . . . \left( \begin{array}{cc}
a_{2 k} & 1 \\
1 & 0
\end{array} \right) M \left( \begin{array}{c}
0 \\
1
\end{array} \right)\, .
\end{equation}
and
$$
G[a_1, a_2, a_3, ..., a_{2k+1}] =
$$
\begin{equation}
\left( \begin{array}{cc}
0 & 1
\end{array} \right) M \left( \begin{array}{cc}
0 & 1 \\
1 & a_1
\end{array} \right) M \left( \begin{array}{cc}
a_2 & 1 \\
1 & 0
\end{array} \right) M . . . \left( \begin{array}{cc}
0 & 1 \\
1 & a_{2k+1}
\end{array} \right) M \left( \begin{array}{c}
1 \\
0
\end{array} \right)\, .
\end{equation}
In these factorizations, all the matrices corresponding to white conways are of the form of the matrix of the conway $a_1$. All the matrices
corresponding to colored conways, the form of the matrix of the conway $a_2$.

\subsection*{Theorem 5.}
The C-function of a rational knot is expressed as the interior product of two vectors associated to rational tangles which components are gaussian brackets
\begin{equation}
G[a_1, a_2, a_3, ..., a_n] = \left( \begin{array}{cc}
G[a_1, ..., a_{j-1}] & G[a_1, ..., a_j]
\end{array} \right) M \left( \begin{array}{c}
G[a_{j+1}, ..., a_n] \\
G[a_{j+2}, ..., a_n]
\end{array} \right)
\end{equation}

\noindent {\bf Proof:} The proof results directly from factorizations (4) and (5). $\quad \square$

Until now, the 2-tangles considered in this paper are rational 2-tangles such as those studied by Kauffman and Lambropoulou \cite{kl}. To a rational
2-tangle one associates two components of C-functions, which for rational 2-tangles are called numerator and denominator \cite{kl}, since they give
the continued fraction decomposition of the ratio of C-functions in the conways $a_j$, as explained before. The cut strands of the 2-tangle, if
reconnected, give two different knots whose C-functions are the numerator and denominator.

On the contrary, the presentation of rational 2-tangles in this paper is related to the decomposition of a rational knot in two 2-tangles,
such that the vectors of these 2-tangles whose components are the so called numerator and denominator have been composed by an interior product defined
by the matrix $M$ to produce the C-function of the dissected knot.

These properties will be generalized for alternating knots and alternating 2-tangles in the next section according to the following axioms.

\subsection*{Axiom 10.} Every alternating 2-tangle has associated two C-functions equal to the C-functions of the knots formed when the 2-tangle
is closed by joining the free edges in two allowed forms.

\subsection*{Axiom 11.} The C-function of all the families of alternating knots is expressed as the internal product of two alternating 2-tangles that are obtained by a dissection of the knot in these two tangles. In most cases this product is not unique because in many families of knots, different splitting are allowed.

\begin{figure}

\psfrag{a}{\LARGE{$a_1$}}
\psfrag{b}{\LARGE{$a_3$}}
\psfrag{d}{\LARGE{$a_3$}}
\psfrag{c}{\LARGE{$a_4$}}

\hfil\scalebox{0.35}{\includegraphics{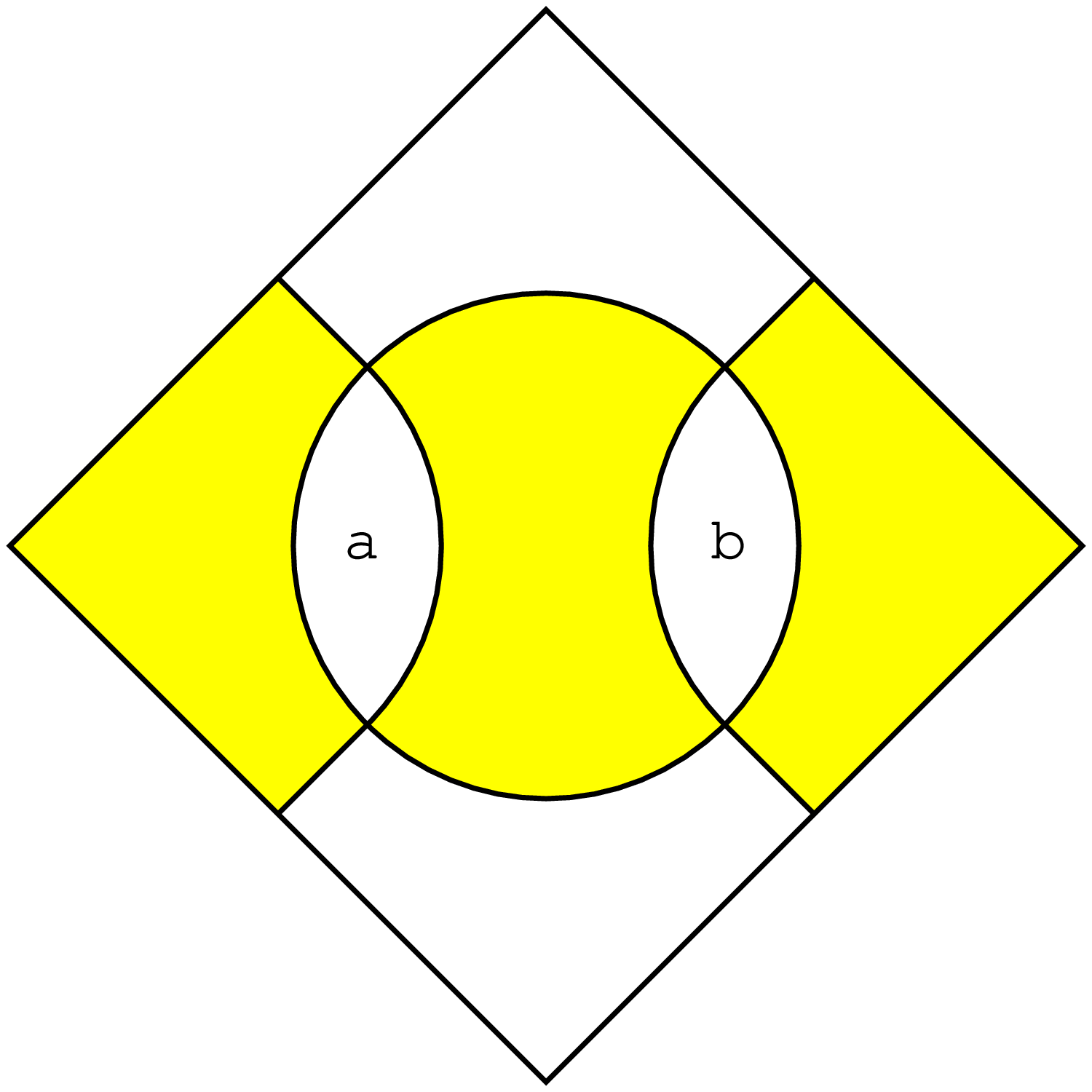}} \quad \scalebox{0.35}{\includegraphics{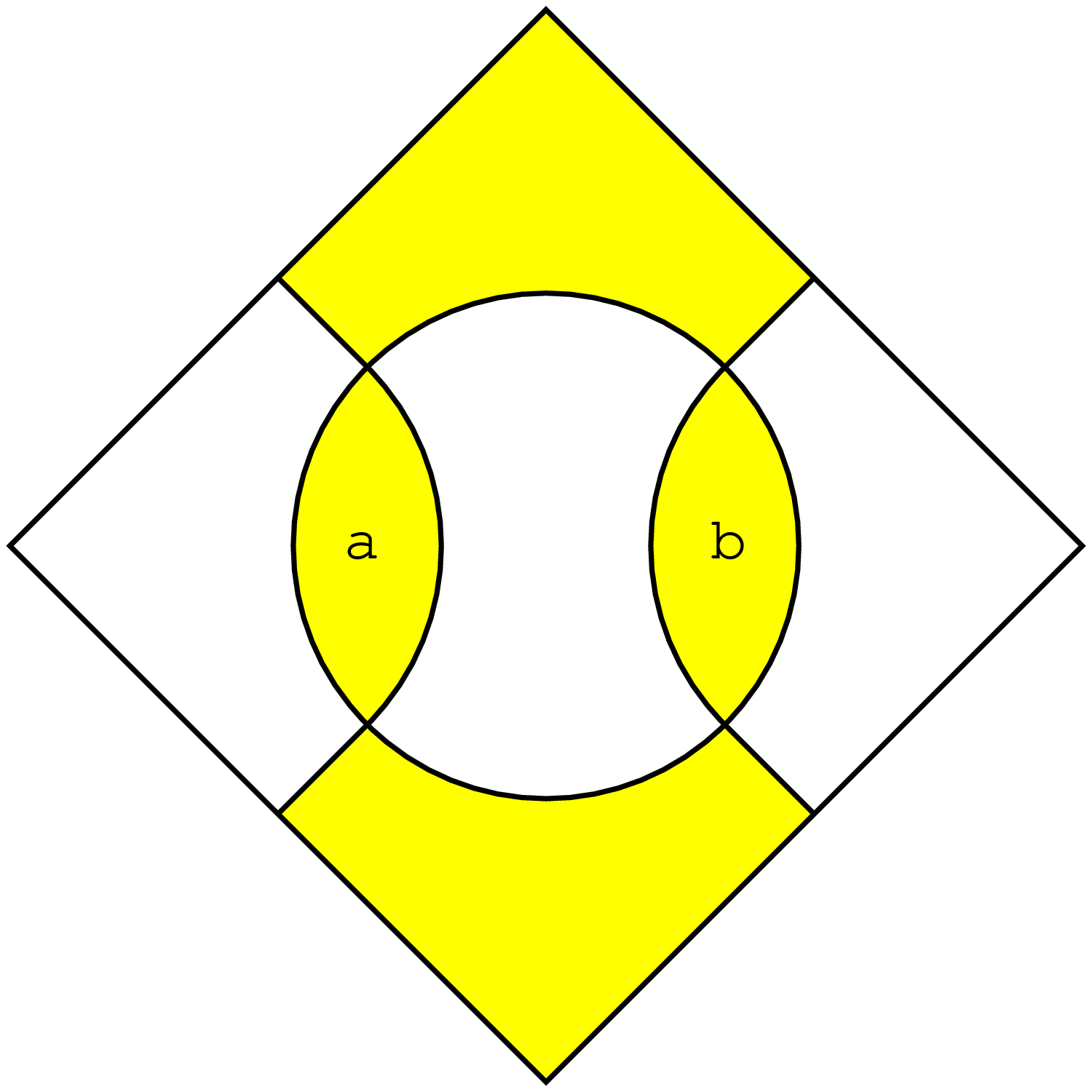}}\hfil

\caption{Non rational tangles formed by two conways. Its covariant and contravariant components are $(a_1 + a_3, a_1 a_3)$ and $(a_1 a_3, a_1 + a_3)$.}

\end{figure}

\subsection*{Definition 4} If $(A_1, A_2)$ are the Conway's functions of one tangle, and $(B_1, B_2)$ of the second tangle, the internal product
of both tangles is defined by
\begin{equation}
A_1 B_2 + A_2 B_1 = \left( A_1, A_2 \right) \left( \begin{array}{cc}
0 & 1 \\
1 & 0
\end{array} \right) \left( \begin{array}{c}
B_1 \\
B_2
\end{array} \right)\, .
\end{equation}

In the rest of this section some rational knots are factorized in two 2-tangles, one of them not rational. These alternating not rational 2-tangles
have two associated components according to Axiom 9. In the next section two alternating 2-tangles are joined to form alternating knots, and
the C-function of a family of alternating knots is computed.

The matrices in product (4) generally do not commute, except just in the trivial case when two contiguous conways have the same value. However
the projections in (5) give the possibility of doing flypes at the left end of the factorization of the rational knots as stated by the property
$$
\left( \begin{array}{cc}
1 & 0
\end{array} \right) M \left( \begin{array}{cc}
0 & 1 \\
1 & a_1
\end{array} \right) M \left( \begin{array}{cc}
a_2 & 1 \\
1 & 0
\end{array} \right) =
$$
$$
\left( \begin{array}{cc}
1 & a_1
\end{array} \right) M \left( \begin{array}{cc}
a_2 & 1 \\
1 & 0
\end{array} \right) = \left( \begin{array}{cc}
0 & 1
\end{array} \right) M \left( \begin{array}{cc}
1 & a_1 \\
a_1 & 0
\end{array} \right) M \left( \begin{array}{cc}
a_2 & 1 \\
1 & 0
\end{array} \right) =
$$
\begin{equation}
\left( \begin{array}{cc}
0 & 1
\end{array} \right) M \left( \begin{array}{cc}
a_2 & 1 \\
1 & 0
\end{array} \right) M \left( \begin{array}{cc}
1 & a_1 \\
a_1 & 0
\end{array} \right) = \left( \begin{array}{cc}
a_2 & 1
\end{array} \right) M \left( \begin{array}{cc}
1 & a_1 \\
a_1 & 0
\end{array} \right) \, .
\end{equation}

Commutation of matrices from the second line to the third line is a particular case of the commutation of matrices in the following lemma.

\subsection*{Lemma 1.}
Symmetric matrices with a zero entry in the same place of the diagonal commute with respect to the interior product defined by $M$.

{\bf Proof:} It results from explicit computation
\begin{equation}
\left( \begin{array}{cc}
A_1 & B_1 \\
B_1 & 0
\end{array} \right) M \left( \begin{array}{cc}
A_2 & B_2 \\
B_2 & 0
\end{array} \right) = \left( \begin{array}{cc}
A_2 B_1 + A_1 B_2 & A_2 B_2 \\
A_2 B_2 & 0
\end{array} \right) \, .
\end{equation}
The other case is proved by transforming these matrices by $M$. $\quad \square$

\subsection*{Remark 1.} Note that this family of matrices defines an abelian group in the field of rational numbers. The unit element is matrix $M$.

A trivial case of lemma 1 is found in equation (8).

The rational knot of three conways that appear in Fig. 7 could be factorized as
\begin{equation}
G[a_1, a_2, a_3] = \left( \begin{array}{cc}
1 + a_1 a_2 & a_1
\end{array} \right) M \left( \begin{array}{c}
1 \\
a_3
\end{array} \right)\, ,
\end{equation}
in terms of the two 2-tangles at the left of Figs. 10 and 11. Using the identity (17) one has
\begin{equation}
G[a_1, a_2, a_3] = \left( \begin{array}{cc}
a_2 & 1
\end{array} \right) M \left( \begin{array}{cc}
1 & a_1 \\
a_1 & 0
\end{array} \right) M \left( \begin{array}{c}
1 \\
a_3
\end{array} \right)\, ,
\end{equation}
which gives the C-function of the family of rational knots with three conways as
\begin{equation}
G[a_1, a_2, a_3] = \left( \begin{array}{cc}
a_2 & 1
\end{array} \right) M \left( \begin{array}{c}
a_1 + a_3 \\
a_1 a_3
\end{array} \right)\, ,
\end{equation}
which is symmetric in conways $a_1$ and $a_3$. This last factorization includes the tangle with C-function $(a_2, 1)$ to the right in Fig. 10 and a
2-tangle that is not rational, with components $(a_1 + a_3, a_1 a_3)$, with the same handedness for their conways $a_1$ and $a_3$. See the left 2-tangle
in Fig. 12. The components of this family of 2-tangles are the sum of the conways $a_1$ and $a_3$, equivalent to one conway, and the product of conways
$a_1 a_3$ corresponding to a composite knot. This family of 2-tangles is important when it appears as the product of tangles producing a family of prime
knots, as for example in equation (21).

\section{ALGEBRA OF FAMILIES OF ALTERNATING KNOTS AND TANGLES}
The algebra previously developed in this paper is generalized to alternating tangles. One notes that the C-function of families of knots can be computed
with the interior product of the vectors associated to two families of tangles whose union forms the family of knots.

\subsection*{Axiom 12.}
For alternating knots and tangles, the connection of tangles is allowed just if colored regions coincide and white regions coincide.

This prescription takes into account the alternating nature of the knots and tangles considered in this paper.

Fig. 10 shows the two different colored tangles of one conway. Fig. 11 shows two colored tangles of two conways
which obey the previous axiom, where the two conways of Fig. 10 were connected according to this axiom.

One has two families of knots of three conways, with the same seed, the trefoil knot. The rational one has one conway with different handedness and two
with the same handedness shown in Fig. 14; and the other family of knots with the three conways with the same handedness as shown in Fig. 15.
Its C-function computed in \cite{pi} is equal to the product of the two vectors associated with 2-tangles of the same handedness in Figs. 10 and 12, namely
\begin{equation}
a_1 a_2 + a_2 a_3 + a_3 a_1 = \left( \begin{array}{cc}
1 & a_1
\end{array} \right) M \left( \begin{array}{c}
a_2 + a_3 \\
a_2 a_3
\end{array} \right)\, .
\end{equation}

The C-functions of the five families of knots formed by four conways, see Figs. 16-20, are easily computed by the interior product of different
combinations of the two tangles in Figs. 11 and 12. The correct factors to be used depends on the coloration. Different factorizations of these
families give the expression for the vectors corresponding to 2-tangles formed by three tangles.

The C-functions of the twelve families of knots formed by five conways were computed by the interior product of one 2-tangle of two conways with one
2-tangle of three conways. In some cases, the C-functions of the 2-tangles involved are determined by closing the 2-tangle in the two allowed forms
to determine the components of the tangle. Two useful properties for this computation are expressed as axioms 13 and 14 below.

\subsection*{Axiom 13.}
The dual vector of a 2-tangle is obtained by interchanging the two components. This corresponds to an interchange of the coloration of the tangle, and
of the handedness of all the conways forming the tangle. Some authors distinguish covariant and
contravariant components of the vector, related by the dual property, or the metric matrix.

\subsection*{Axiom 14.}
Assuming  two families of alternating knots or tangles are different in the handedness of one conway (with conway equal to $a_j$). If the C-function of one family is expressed as $a_j A + B$ (with $A$ and $B$ C-functions), then the C-function of the other family is $A + a_j B$.

By definition, any family of knots is the composition of two 2-tangles: any conway of the family and the rest of the family which is a family of
2-tangles. This last is, by definition, formed by conways.

I have used two or three colors in Figs. 26-52 to show the decomposition of the families of knots in families of tangles with a corresponding
factorization of the C-functions. For the families in Figs. 27-29 and 32-45 one factorization with two 2-tangles formed by three conways was illustrated with two colors. The components of the corresponding 2-tangles are easily recognized from the expression of the C-function in the figure caption.

\subsection*{Definition 5.} Two families of 3-tangles are formed when a family of knots is dissected by cutting three different strands of the seed.

Four families of knots formed by six conways in Figs. 30 and 31 can not be factorized in two 2-tangles of 3 conways, but factorizes into two 2-tangles of two and four conways, that are represented with three colors, or alternatively they are factorized in two 3-tangles of three conways. Two of these families are separated in two 3-tangles identified by colors green and red.

The seven last families represented in Figs. 46-52 can be separated in the product of one conway with a 2-tangle of five conways, or as shown in the figure with colors green and red, they are separated in two 3-tangles of three conways.

A linear algebra for the families that separates in two 3-tangles has been developed for the simple cases that appear in families of alternating knots of six conways. The algebra is now in five dimensions. The metric matrix to express the C-function of these families as an interior product of 5-dimensional vectors is
\begin{equation}
\left( \begin{array}{ccccc}
0 & 0 & 0 & 0 & 1 \\
0 & 0 & 1 & 1 & 0 \\
0 & 1 & 0 & 1 & 0 \\
0 & 1 & 1 & 0 & 0 \\
1 & 0 & 0 & 0 & 0
\end{array} \right)
\end{equation}

In the symmetric 3-tangles two vectors appear
$$
(a_1 a_2 a_3, a_2 a_3, a_3 a_1, a_1 a_2, a_1 + a_2 + a_3) \, , \quad (a_1 a_2 + a_2 a_3 + a_3 a_1, a_1, a_2, a_3, 1)\, .
$$

In the non symmetric tangles (with $a_1$ non symmetric) they appear the vectors
$$
(a_2 a_3, a_1 a_2 a_3, a_3, a_2, 1 + a_1 a_2 + a_1 a_3) \, , \quad (a_2 + a_1 a_2 a_3 + a_3, 1, a_1 a_2, a_3 a_1, a_1)\, .
$$

A table with the same figures of the 65 families of alternating families of knots and explicit factorizations of the C-functions will be released
simultaneously with this paper \cite{tb}. One factorization is associated to two colors. The corresponding factors of this factorization provide
a table of the vectors associated to 2-tangles formed by three conways.

\


\begin{figure}

\section*{Appendix}

{\bf FAMILIES OF PRIME ALTERNATING KNOTS FROM ONE TO SIX CONWAYS}

\psfrag{p}{\LARGE{$a_1$}}
\psfrag{q}{\LARGE{$a_2$}}
\hfil\scalebox{0.5}{\includegraphics{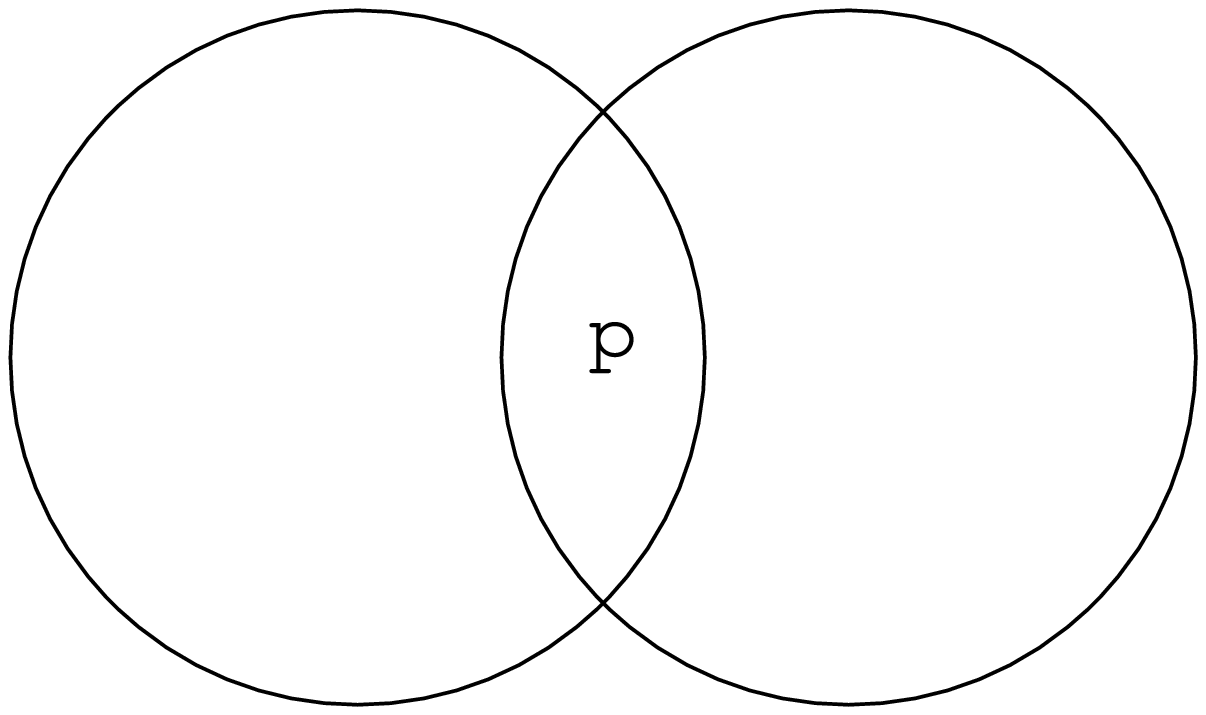}}\hfil\scalebox{0.5}{\includegraphics{nudo4-2.ps}}\hfil

\

\caption{Rational knots with one or two conways. Its Conway's functions are $a_1$ and $1 + a_1 a_2$.}
\end{figure}

\begin{figure}

FAMILIES OF KNOTS WITH THREE CONWAYS. ITS SEED IS THE TREFOIL KNOT. TWO CASES

\psfrag{a}{\LARGE{$a_1$}}
\psfrag{b}{\LARGE{$a_2$}}
\psfrag{c}{\LARGE{$a_3$}}
\hfil\scalebox{0.5}{\includegraphics{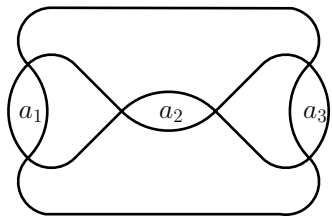}}\hfil

\caption{Rational knot formed by three conways, two of them with one handedness, the other ($a_2$) with the opposite handedness. Its Conway's function is $a_1 a_2 a_3 + a_1 + a_3$.}

\end{figure}

\begin{figure}

\psfrag{a}{\LARGE{$a_1$}}
\psfrag{b}{\LARGE{$a_2$}}
\psfrag{c}{\LARGE{$a_3$}}
\hfil\scalebox{0.5}{\includegraphics{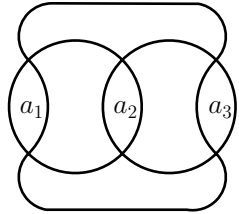}}\hfil

\caption{Knot formed by three conways with the same handedness. Its Conway's function is $a_1 a_2 + a_2 a_3 + a_3 a_1 \, .$ The first member of this family in the knot table is $8_5$. In the links table the first examples are $6_3^1$ and $7_4^2$.}

\end{figure}

\begin{figure}
FAMILIES OF KNOTS OF FOUR CONWAYS WITH THE LINK CALLED THE SOLOMON KNOT AS SEED. TWO CASES

\psfrag{a}{\LARGE{$a_1$}}
\psfrag{b}{\LARGE{$a_2$}}
\psfrag{c}{\LARGE{$a_3$}}
\psfrag{d}{\LARGE{$a_4$}}
\hfil\scalebox{0.5}{\includegraphics{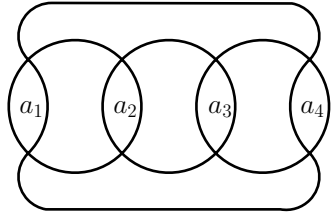}}\hfil

\caption{Family of knots formed by four conways with the same handedness. Its C-function is $a_1 a_2 a_3 + a_2 a_3 a_4 + a_3 a_4 a_1 + a_4 a_1 a_2$.}

\end{figure}

\begin{figure}

\psfrag{a}{\LARGE{$a_1$}}
\psfrag{b}{\LARGE{$a_2$}}
\psfrag{c}{\LARGE{$a_3$}}
\psfrag{d}{\LARGE{$a_4$}}
\hfil\scalebox{0.5}{\includegraphics{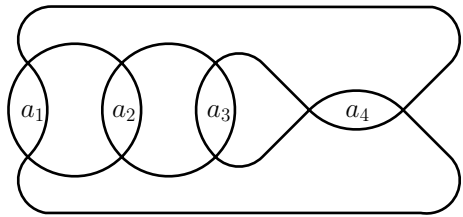}}\hfil

\caption{Knot formed by three conways with the same handedness and one olf the opposite handedness. Its C-function is
$a_1 a_2 a_3 a_4 + a_1 a_2 + a_2 a_3 + a_3 a_1 \, .$}

\end{figure}

\begin{figure}

KNOTS OF FOUR CONWAYS WITH THE KNOT $4_1$ AS SEED. THREE CASES

\psfrag{a}{\LARGE{$a_1$}}
\psfrag{b}{\LARGE{$a_2$}}
\psfrag{c}{\LARGE{$a_3$}}
\psfrag{d}{\LARGE{$a_4$}}
\hfil\scalebox{0.5}{\includegraphics{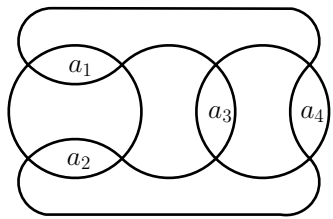}}\hfil

\caption{Knot formed by four conways in two pairs of opposite handedness. Its C-function is $a_1 a_2 a_3 a_4 + (a_1 + a_2)(a_3 + a_4) \, .$}
\end{figure}

\begin{figure}

\psfrag{a}{\LARGE{$a_1$}}
\psfrag{b}{\LARGE{$a_2$}}
\psfrag{c}{\LARGE{$a_3$}}
\psfrag{d}{\LARGE{$a_4$}}
\hfil\scalebox{0.5}{\includegraphics{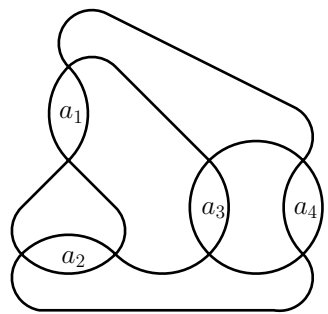}}\hfil

\caption{Knot formed by two conways with the same handedness connected to a T-tangle of two conways of opposite handedness. Its C-function is
$(a_1 a_2 + 1)(a_3 + a_4) + a_2 a_3 a_4\, .$}
\end{figure}

\begin{figure}

\

\

\psfrag{a}{\LARGE{$a_1$}}
\psfrag{b}{\LARGE{$a_2$}}
\psfrag{c}{\LARGE{$a_3$}}
\psfrag{d}{\LARGE{$a_4$}}
\hfil\scalebox{0.5}{\includegraphics{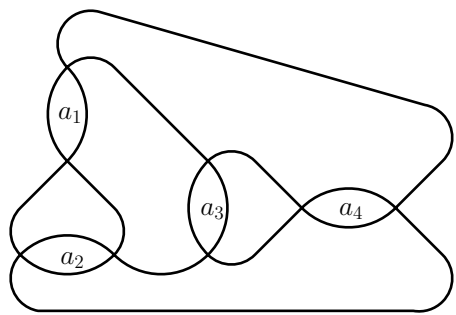}}\hfil

\caption{Rational knot formed by four conways. Its C-function is $(a_1 a_2 + 1)(a_3 a_4 + 1) + a_2 a_3\, $. An example is the knot $8_{12}$.}

\end{figure}

\begin{figure}
FAMILIES OF KNOTS FORMED BY FIVE CONWAYS WITH THE KNOT $5_1$ AS SEED. TWO CASES
\psfrag{a}{\LARGE{$a_1$}}
\psfrag{b}{\LARGE{$a_2$}}
\psfrag{c}{\LARGE{$a_3$}}
\psfrag{d}{\LARGE{$a_4$}}
\psfrag{e}{\LARGE{$a_5$}}

\hfil\scalebox{0.45}{\includegraphics{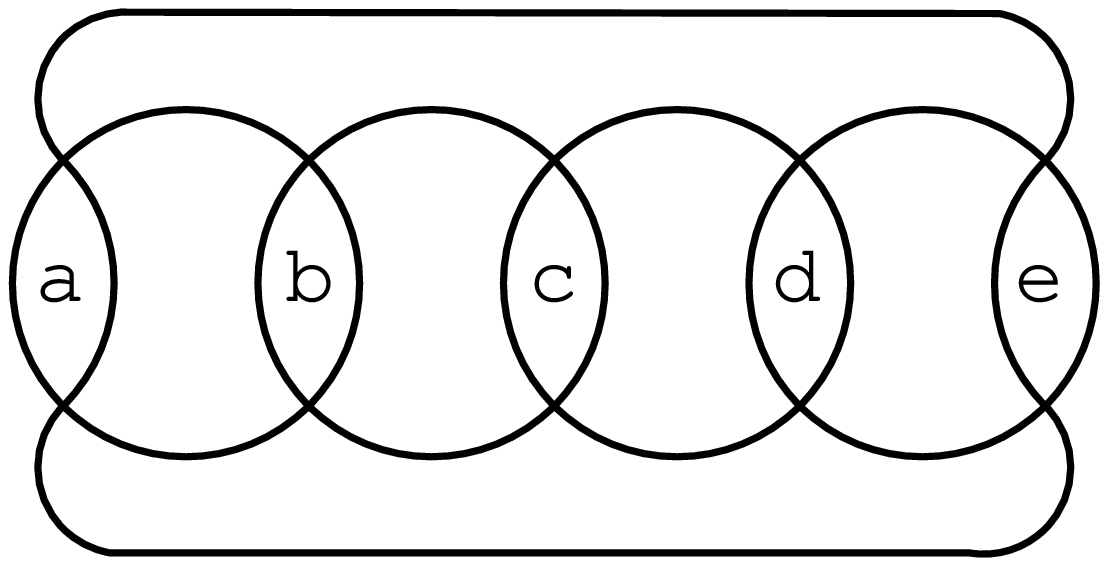}} \quad \quad \scalebox{0.45}{\includegraphics{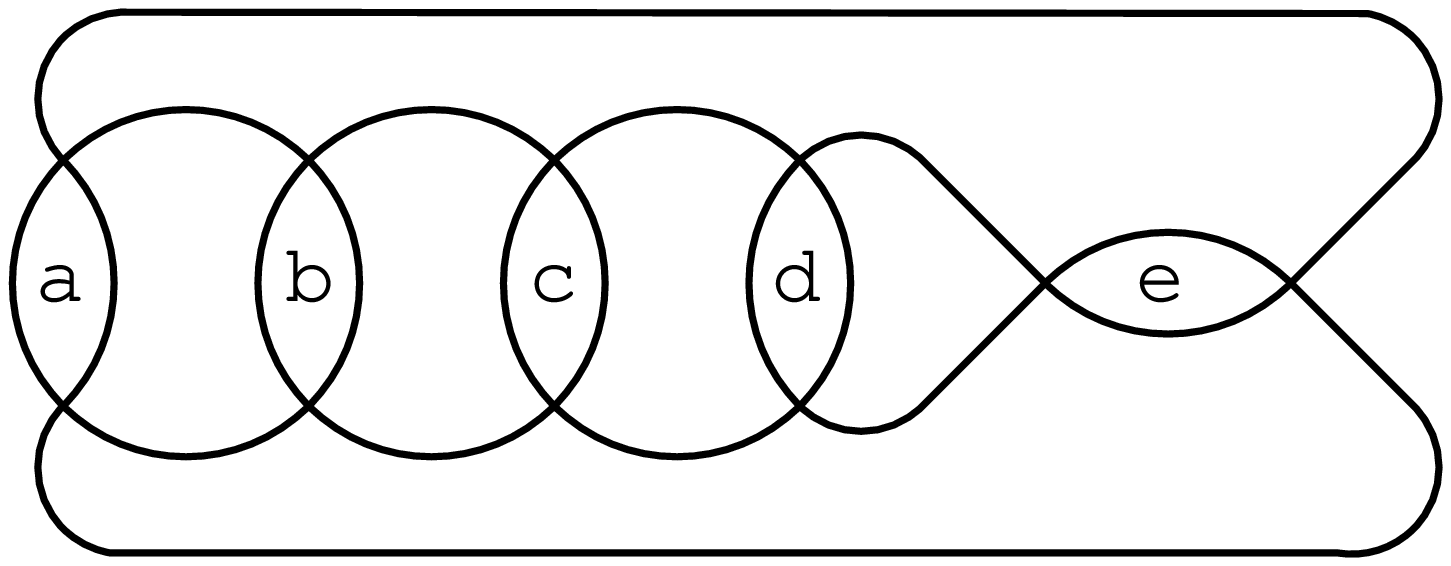}}\hfil

\caption{Families of knots formed by five conways, five or four with the same handedness with C-functions with five terms.
$ a_1 a_2 a_3 a_4 + a_2 a_3 a_4 a_5 + a_3 a_4 a_5 a_1 + a_4 a_5 a_1 a_2 + a_5 a_1 a_2 a_3 $ and
$a_1 a_2 a_3 + a_2 a_3 a_4 + a_3 a_4 a_1 + a_4 a_1 a_2 + a_1 a_2 a_3 a_4 a_5$, respectively.}

\end{figure}

\begin{figure}
FAMILIES OF FIVE CONWAYS ASOCIATED TO THE KNOT $5_2$ AS SEED, FOUR CASES
\psfrag{a}{\LARGE{$a_1$}}
\psfrag{b}{\LARGE{$a_2$}}
\psfrag{c}{\LARGE{$a_3$}}
\psfrag{d}{\LARGE{$a_4$}}
\psfrag{e}{\LARGE{$a_5$}}

\scalebox{0.45}{\includegraphics{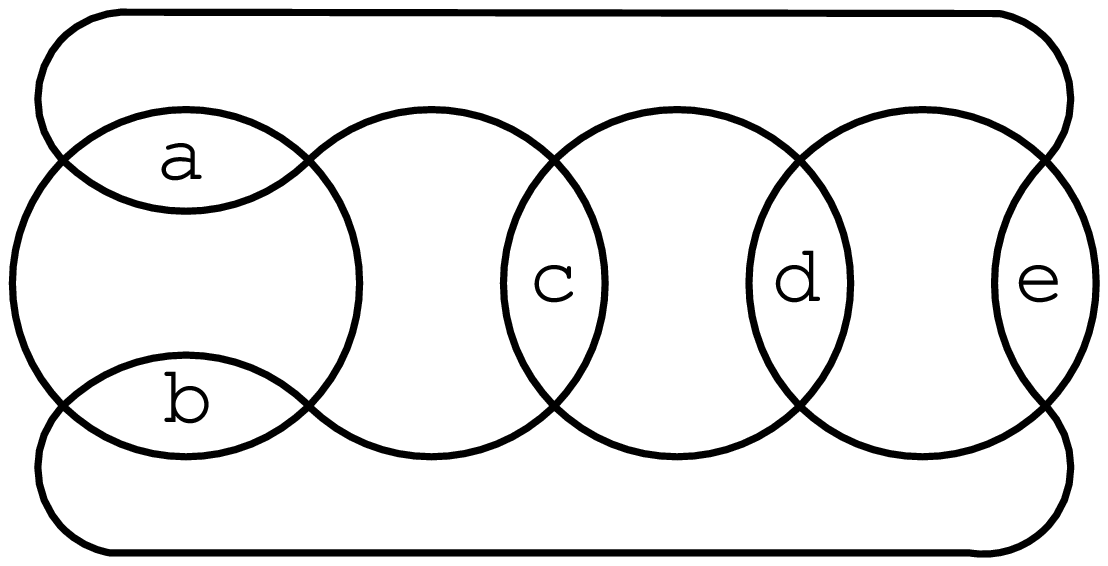}}\quad \scalebox{0.45}{\includegraphics{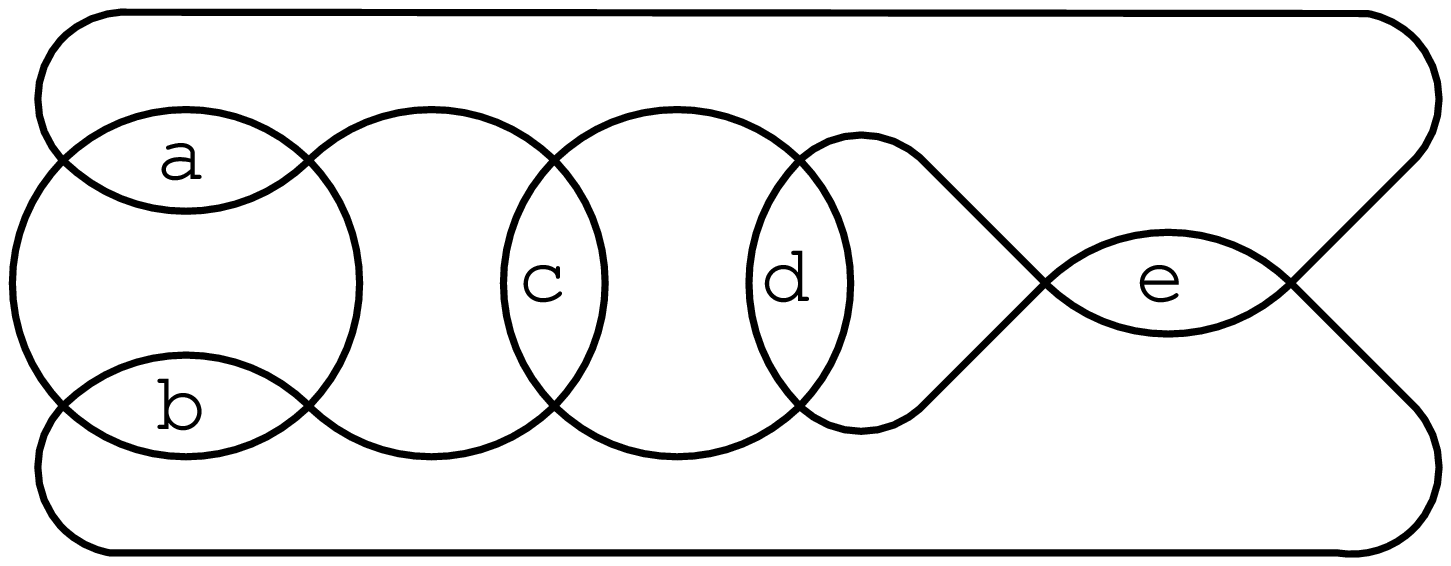}}

\caption{Families with five conways, two with one handedness the other three with the opposite handedness, or two pairs of different handedness. Its
C-functions of seven terms are
$(a_1 + a_2)(a_3 a_4 + a_4 a_5 + a_5 a_3) + a_1 a_2 a_3 a_4 a_5$ and
$a_1 a_2 a_3 a_4 + (a_1 + a_2)(a_3 + a_4 + a_3 a_4 a_5)$.}

\end{figure}

\begin{figure}
\psfrag{a}{\LARGE{$a_1$}}
\psfrag{b}{\LARGE{$a_2$}}
\psfrag{c}{\LARGE{$a_3$}}
\psfrag{d}{\LARGE{$a_4$}}
\psfrag{e}{\LARGE{$a_5$}}

\hfil\scalebox{0.45}{\includegraphics{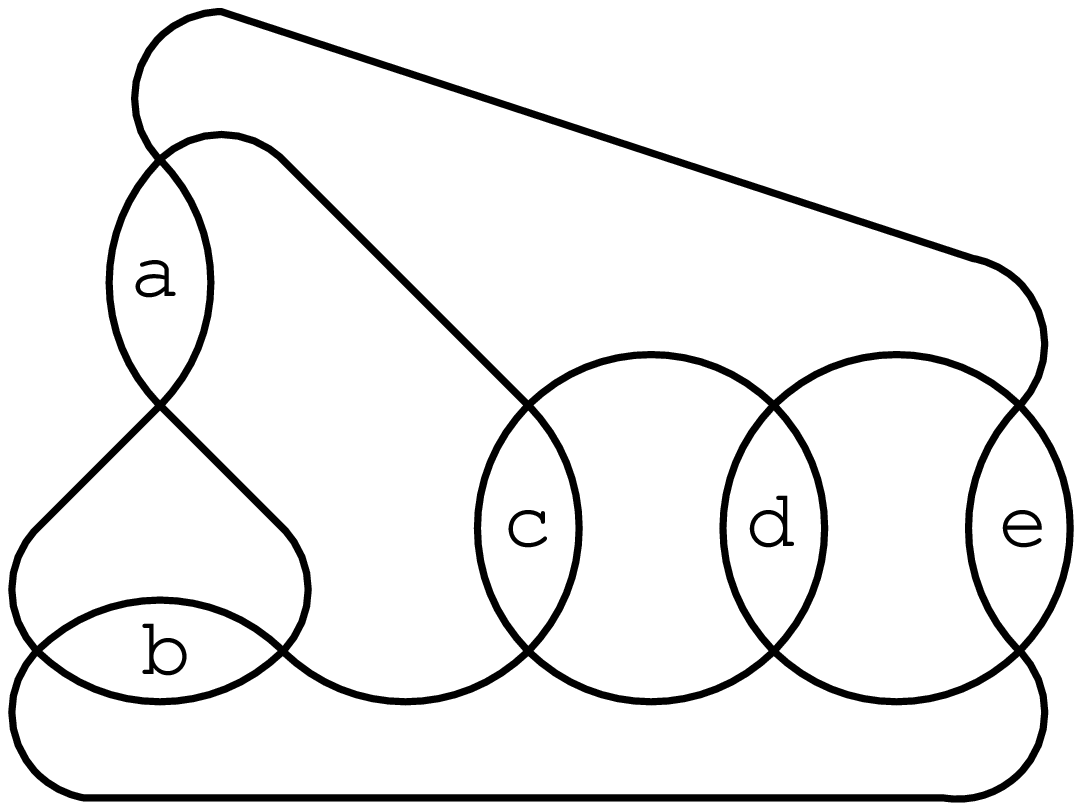}} \quad \quad \scalebox{0.45}{\includegraphics{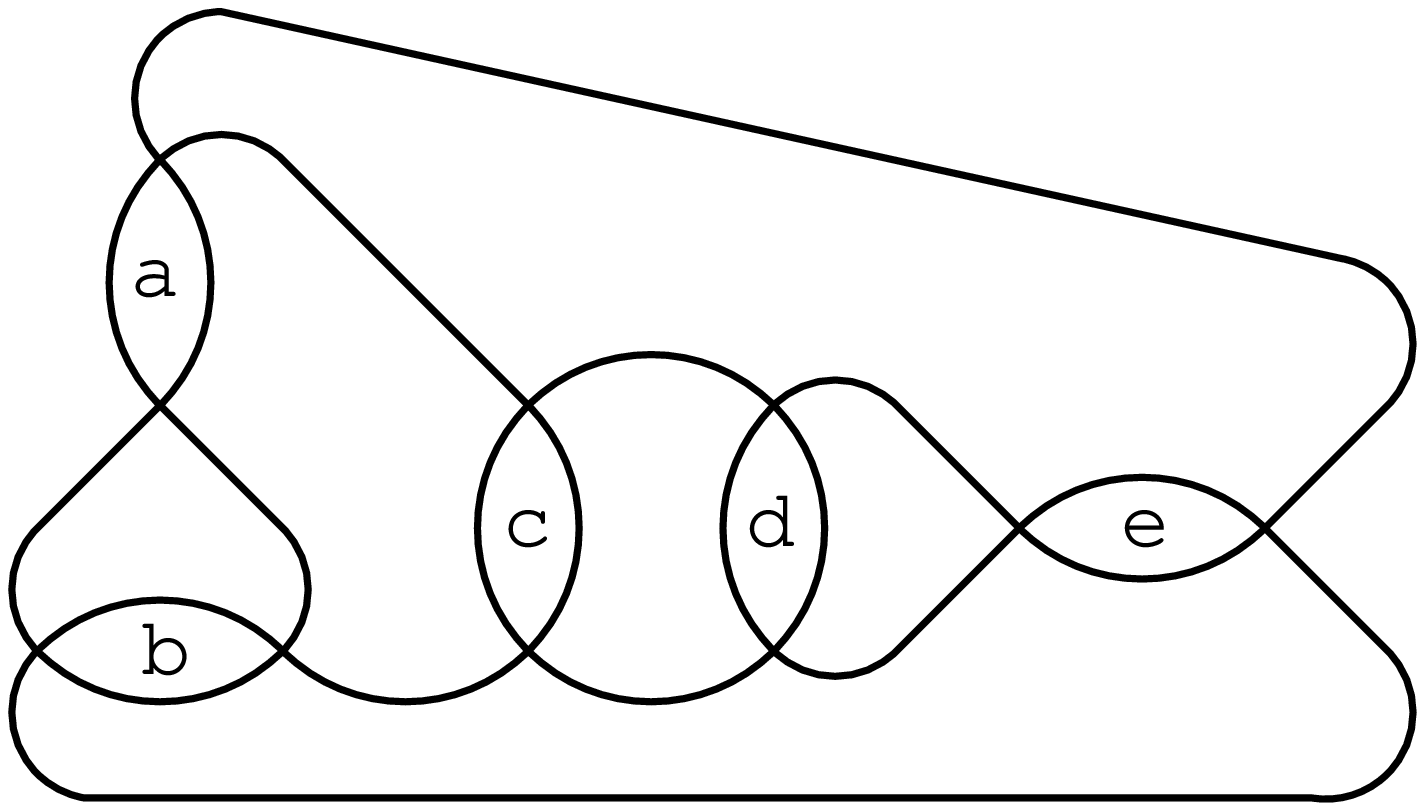}}\hfil

\caption{Continuation of families with five conways. One has now one or two T-tangles and two or three conways with the same handedness with C-functions of seven
terms $(a_1 a_2 + 1)(a_3 a_4 + a_4 a_5 + a_5 a_3) + a_2 a_3 a_4 a_5$ and
$a_2 a_3 a_4 + (1 + a_1 a_2)(a_3 + a_4 + a_3 a_4 a_5) $.}

\end{figure}

\begin{figure}
FAMILIES FORMED BY FIVE CONWAYS WITH THE WHITEHEAD LINK $5_1^2$ AS SEED, SIX CASES
\psfrag{a}{\LARGE{$a_1$}}
\psfrag{d}{\LARGE{$a_2$}}
\psfrag{c}{\LARGE{$a_3$}}
\psfrag{b}{\LARGE{$a_4$}}
\psfrag{e}{\LARGE{$a_5$}}

\hfil \scalebox{0.45}{\includegraphics{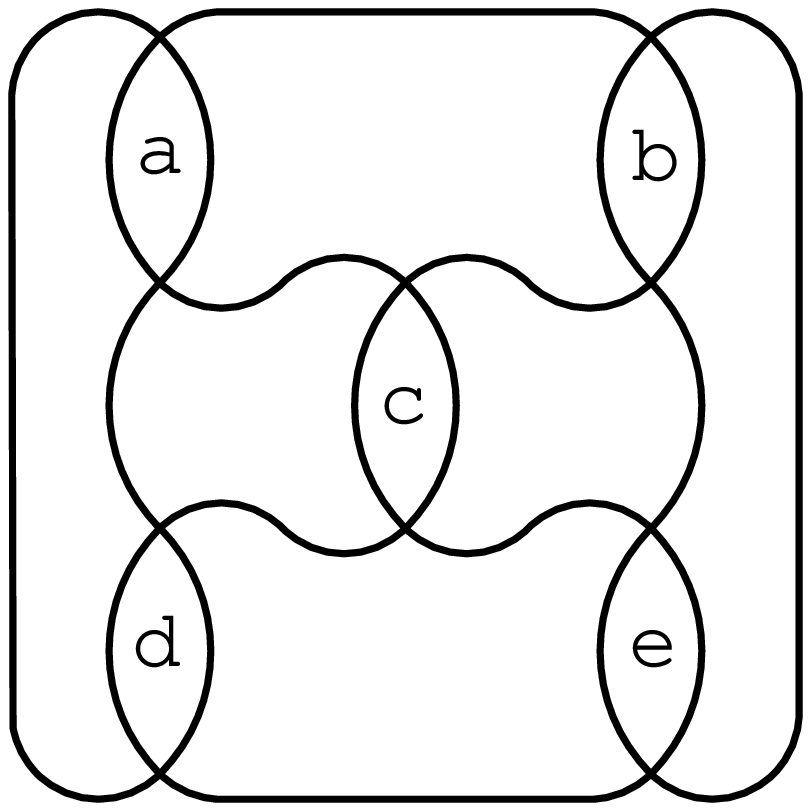}} \quad \quad \scalebox{0.45}{\includegraphics{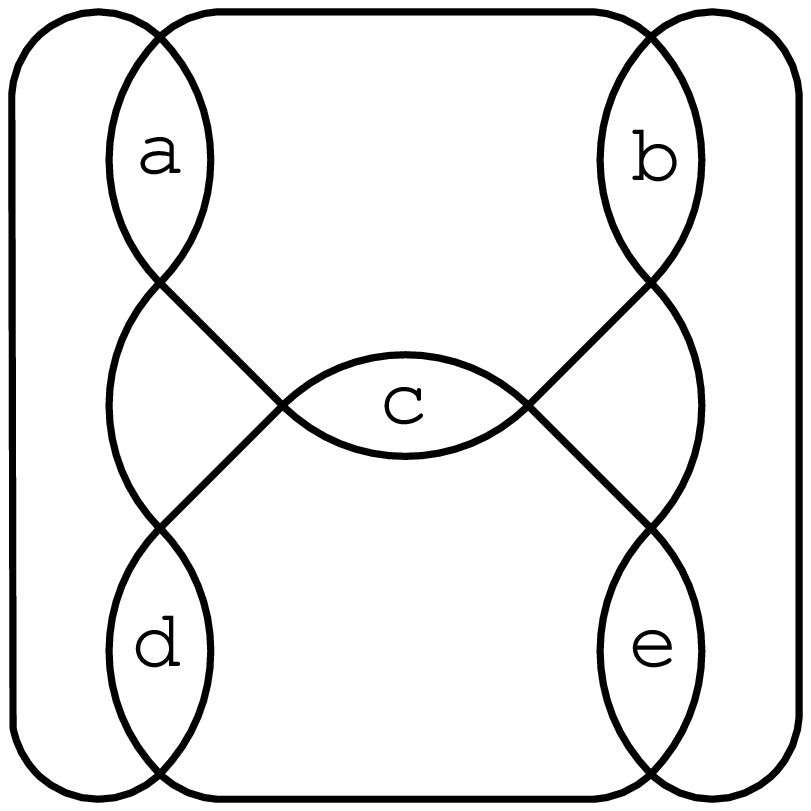}} \hfil

\caption{Families with five conways. Its C-functions with eight monomials are
$(a_1 + a_2)(a_4 + a_5) + a_1 a_2 a_3 (a_4 + a_5) + (a_1 + a_2)a_3 a_4 a_5$ and
$(a_1 + a_2) a_3 (a_4 + a_5) + a_1 a_2 (a_4 + a_5) + (a_1 + a_2) a_4 a_5$.}

\end{figure}

\begin{figure}
\psfrag{a}{\LARGE{$a_1$}}
\psfrag{d}{\LARGE{$a_2$}}
\psfrag{c}{\LARGE{$a_3$}}
\psfrag{b}{\LARGE{$a_4$}}
\psfrag{e}{\LARGE{$a_5$}}

\hfil\scalebox{0.45}{\includegraphics{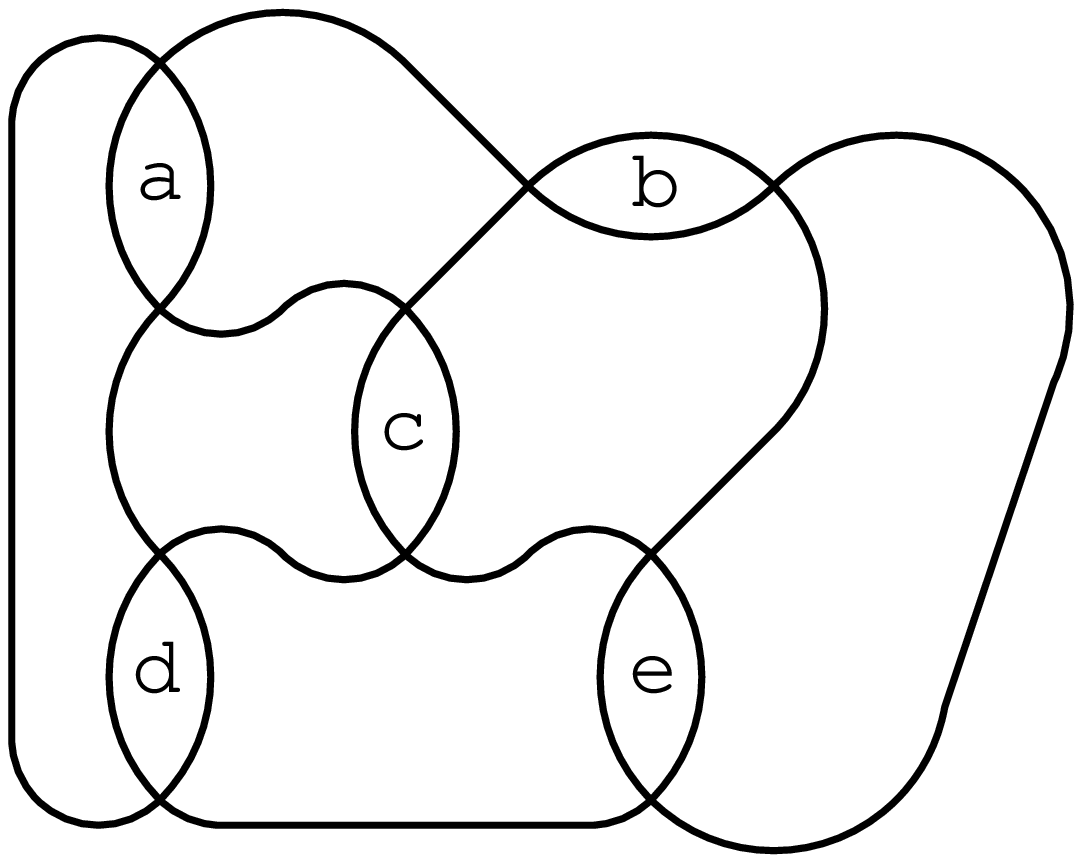}} \quad \quad \quad \quad\scalebox{0.45}{\includegraphics{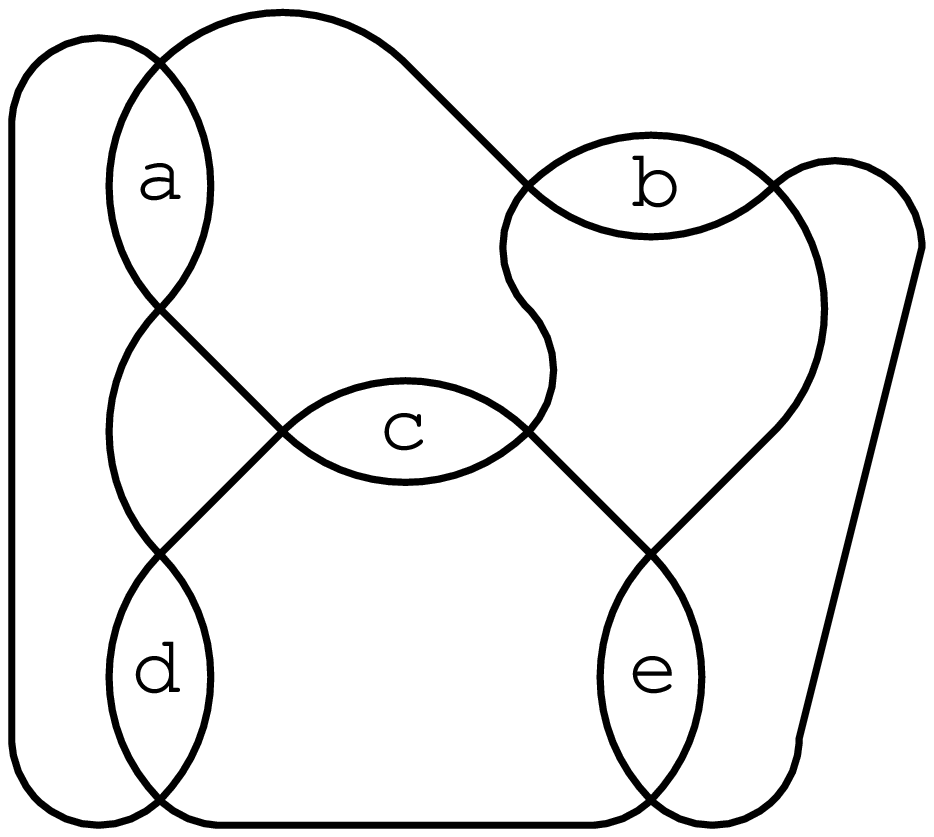}}\hfil

\caption{Continuation of the families having five conways. Its C-functions have eight monomials
$(a_1 a_2 a_3 + a_1 + a_2)(1 + a_4 a_5) + (a_1 + a_2) a_3 a_5$ and
$(a_1 + a_2)(a_3 a_4 a_5 + a_3 + a_5) + a_1 a_2 (1 + a_4 a_5)\, .$}

\end{figure}

\begin{figure}
\psfrag{a}{\LARGE{$a_1$}}
\psfrag{d}{\LARGE{$a_2$}}
\psfrag{c}{\LARGE{$a_3$}}
\psfrag{b}{\LARGE{$a_4$}}
\psfrag{e}{\LARGE{$a_5$}}

\hfil\scalebox{0.45}{\includegraphics{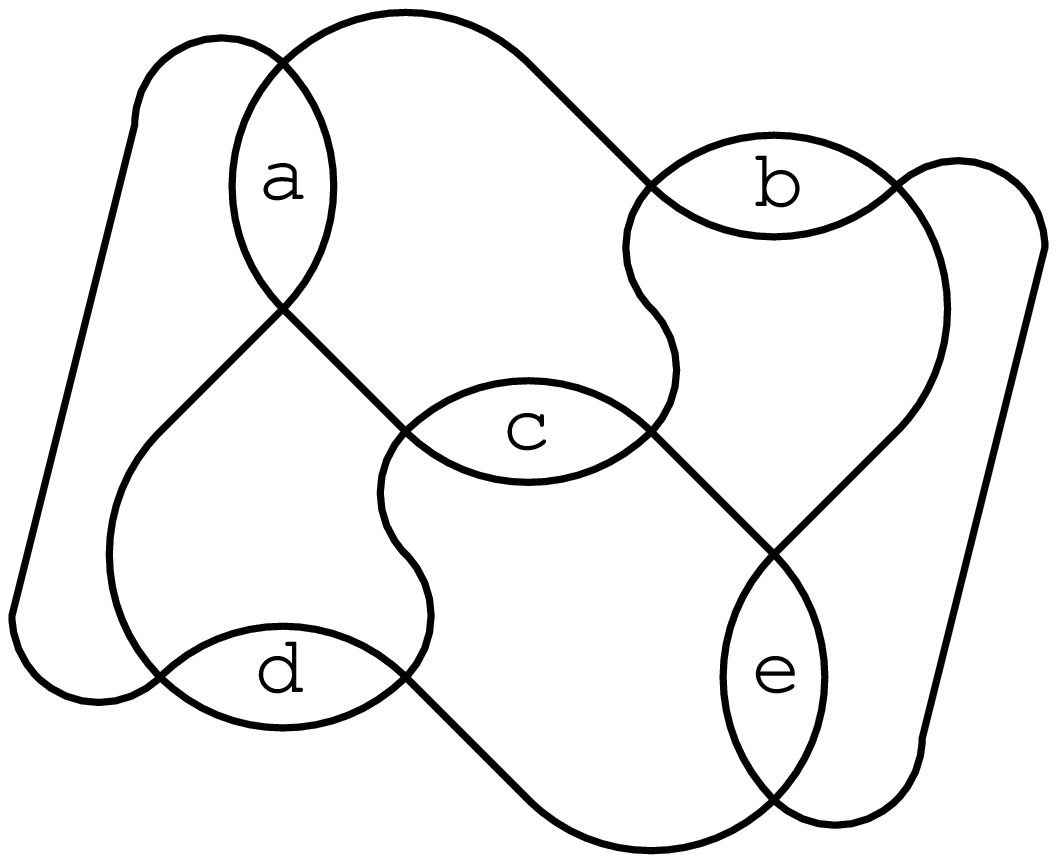}} \quad \quad \scalebox{0.45}{\includegraphics{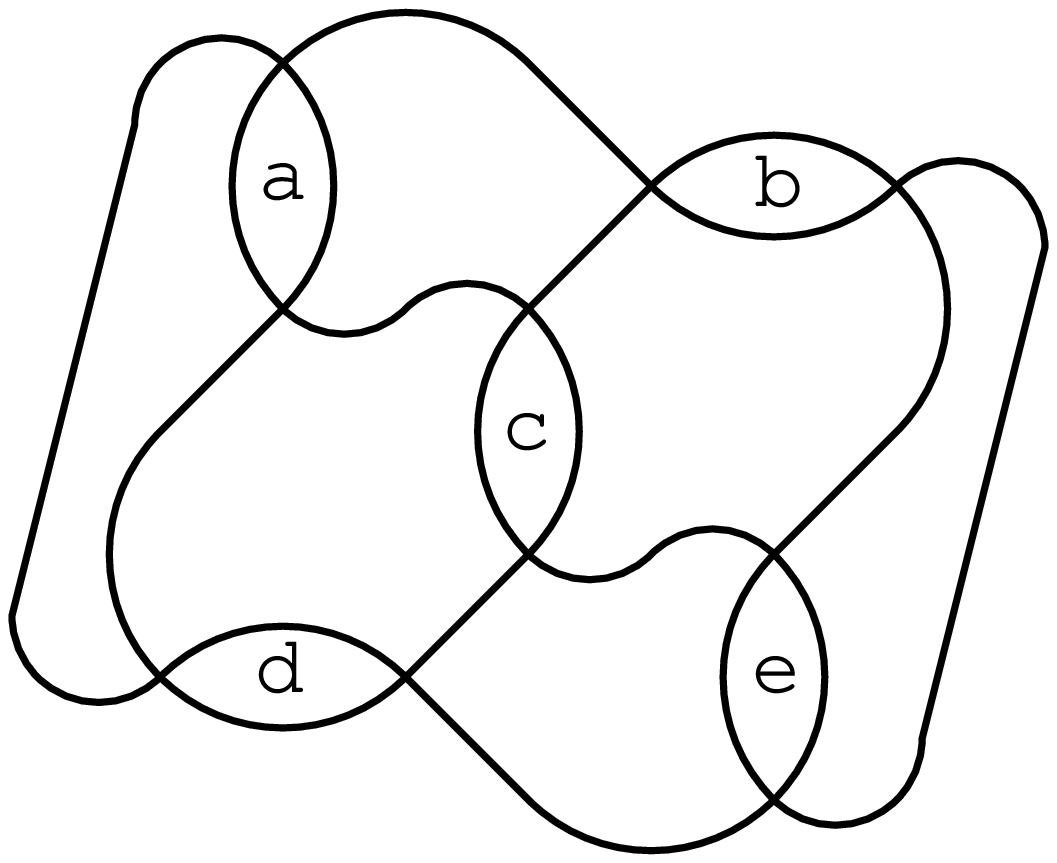}}\hfil

\caption{Termination of the families with five conways and C-functions with eight terms
$a_1 + a_5 + a_1 (a_2 + a_4) a_5 + (a_1 a_2 + 1) a_3 (a_4 a_5 + 1)$, which is rational, and
$1 + a_1 (a_2 + a_3) + (a_3 + a_4) a_5 + a_1 (a_2 a_3 + a_3 a_4 + a_4 a_2) a_5\, .$}

\end{figure}

\begin{figure}

TWO CASES OF FAMILIES OF KNOTS FORMED BY SIX CONWAYS AND THE LINK $6_2^1$ AS SEED

\psfrag{a}{\LARGE{$a_1$}}
\psfrag{b}{\LARGE{$a_2$}}
\psfrag{c}{\LARGE{$a_3$}}
\psfrag{d}{\LARGE{$a_4$}}
\psfrag{e}{\LARGE{$a_5$}}
\psfrag{f}{\LARGE{$a_6$}}
\scalebox{0.40}{\includegraphics{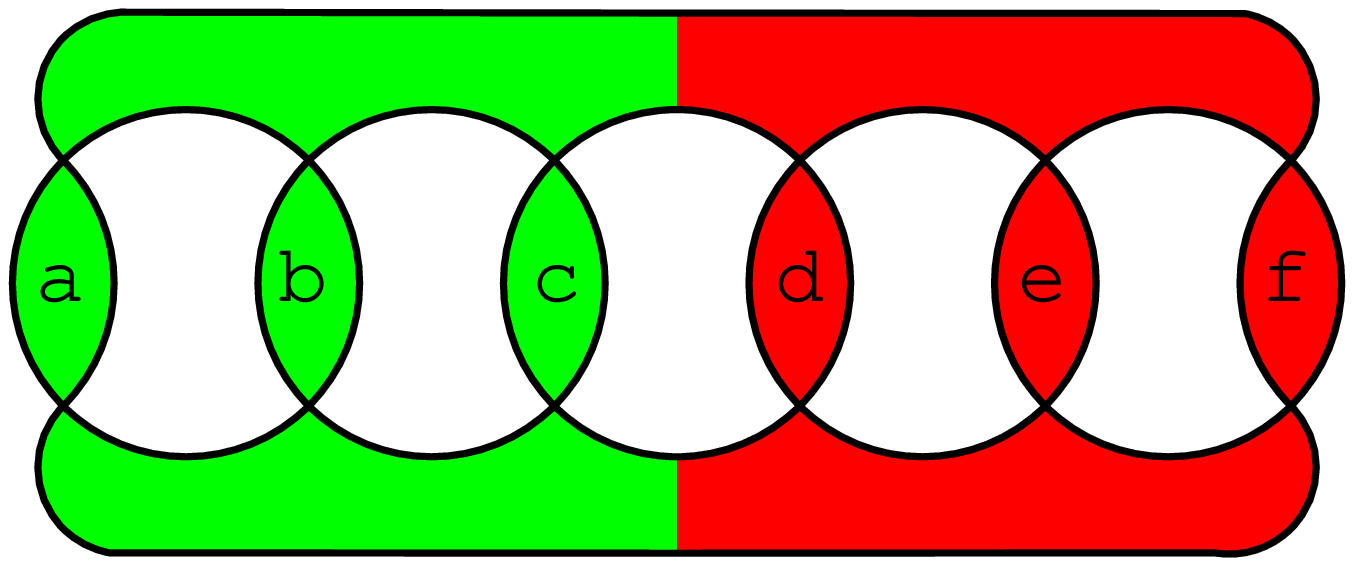}} \scalebox{0.40}{\includegraphics{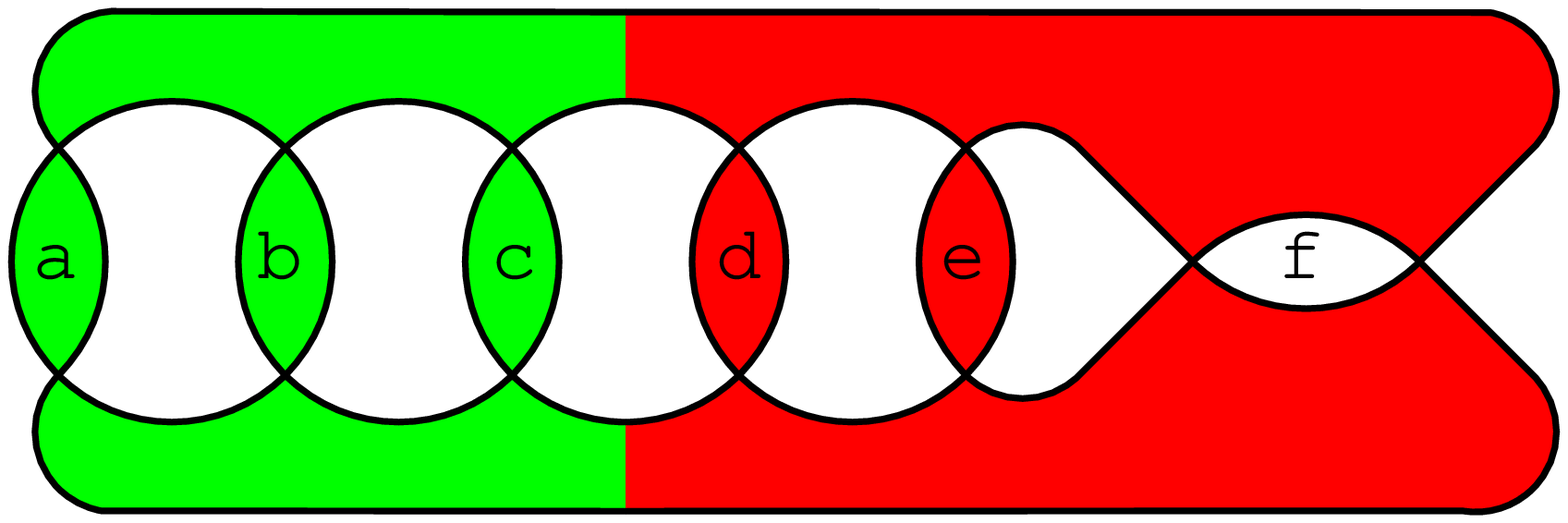}}

\caption{Families of knots formed by six conways. Its C-function are factored in the product of two tangles formed each by three conways
$(a_1 a_2 + a_2 a_3 + a_3 a_1) a_4 a_5 a_6 + a_1 a_2 a_3 (a_4 a_5 + a_5 a_6 + a_6 a_4)$
and $(a_1 a_2 + a_2 a_3 + a_3 a_1) a_4 a_5 + a_1 a_2 a_3 (a_4 + a_5 + a_4 a_5 a_6)$}
\end{figure}

\begin{figure}

FOUR CASES OF FAMILIES OF KNOTS FORMED BY SIX CONWAYS WITH THE KNOT $6_1$ AS SEED

\psfrag{a}{\LARGE{$a_1$}}
\psfrag{b}{\LARGE{$a_2$}}
\psfrag{c}{\LARGE{$a_3$}}
\psfrag{d}{\LARGE{$a_4$}}
\psfrag{e}{\LARGE{$a_5$}}
\psfrag{f}{\LARGE{$a_6$}}

\scalebox{0.40}{\includegraphics{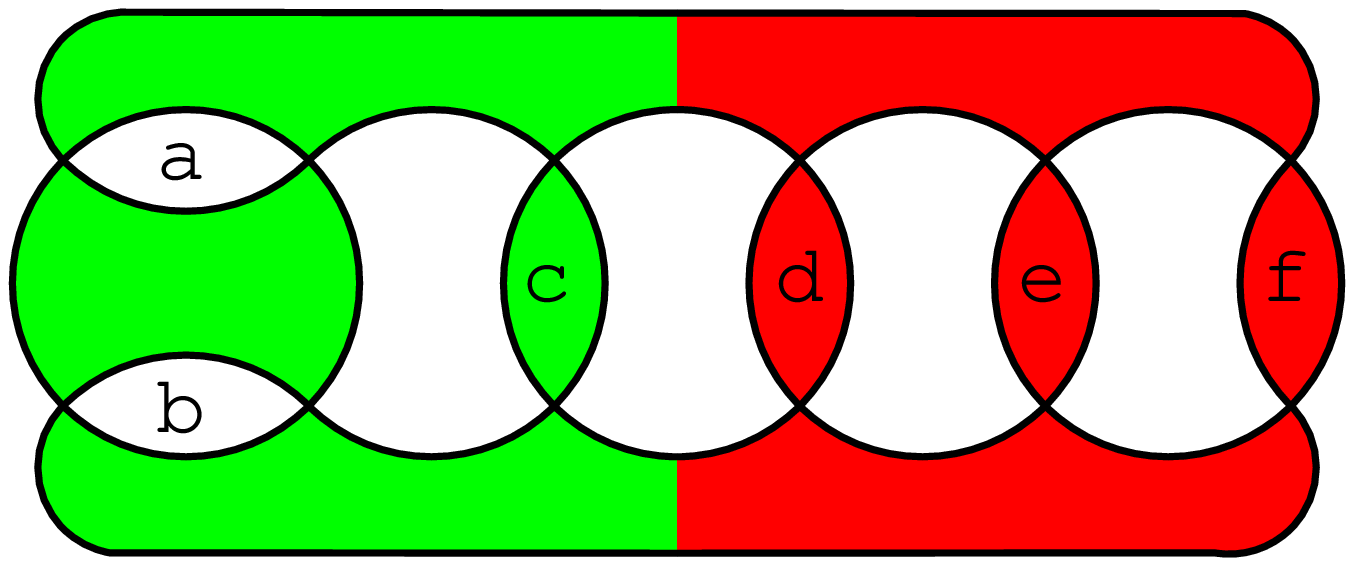}} \scalebox{0.40}{\includegraphics{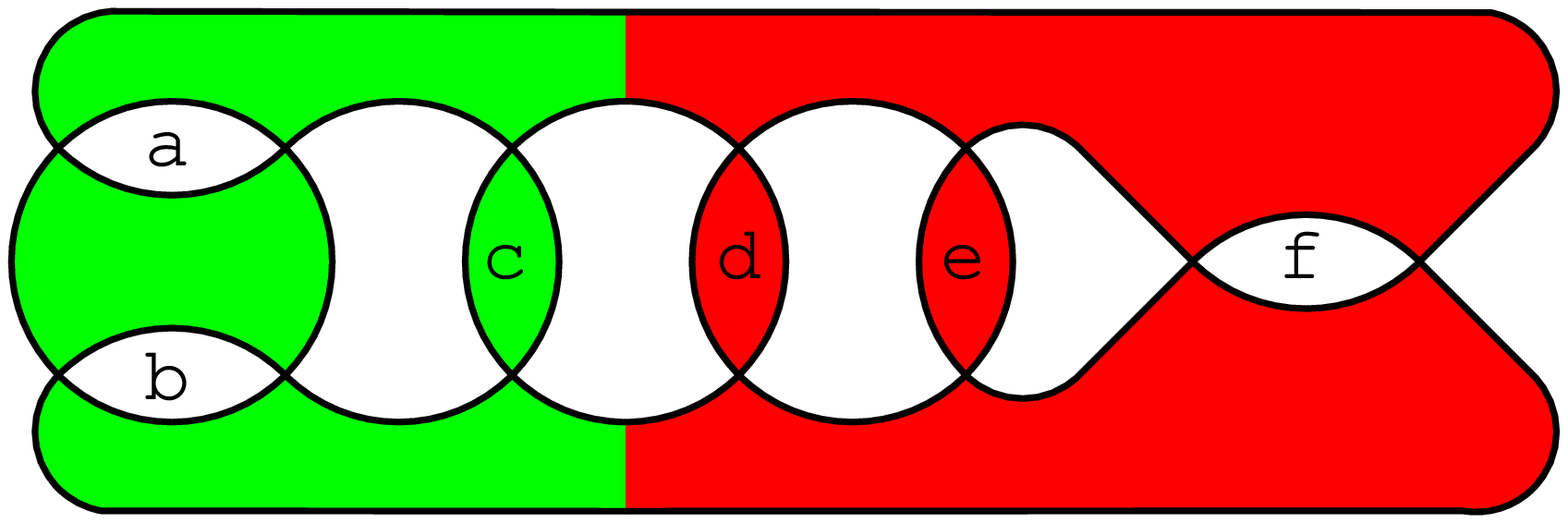}}

\caption{Families of knots formed by six conways. Its C-functions are factored in the product of two tangles formed each by three conways
$(a_1 + a_2 + a_1 a_2 a_3) a_4 a_5 a_6 + (a_1 + a_2) a_3 (a_4 a_5 + a_5 a_6 + a_6 a_4)$ and
$(a_1 + a_2 + a_1 a_2 a_3) a_4 a_5 + (a_1 + a_2) a_3 (a_4 + a_5 + a_4 a_5 a_6)$.}
\end{figure}

\begin{figure}

\psfrag{a}{\LARGE{$a_1$}}
\psfrag{b}{\LARGE{$a_2$}}
\psfrag{c}{\LARGE{$a_3$}}
\psfrag{d}{\LARGE{$a_4$}}
\psfrag{e}{\LARGE{$a_5$}}
\psfrag{f}{\LARGE{$a_6$}}

\scalebox{0.40}{\includegraphics{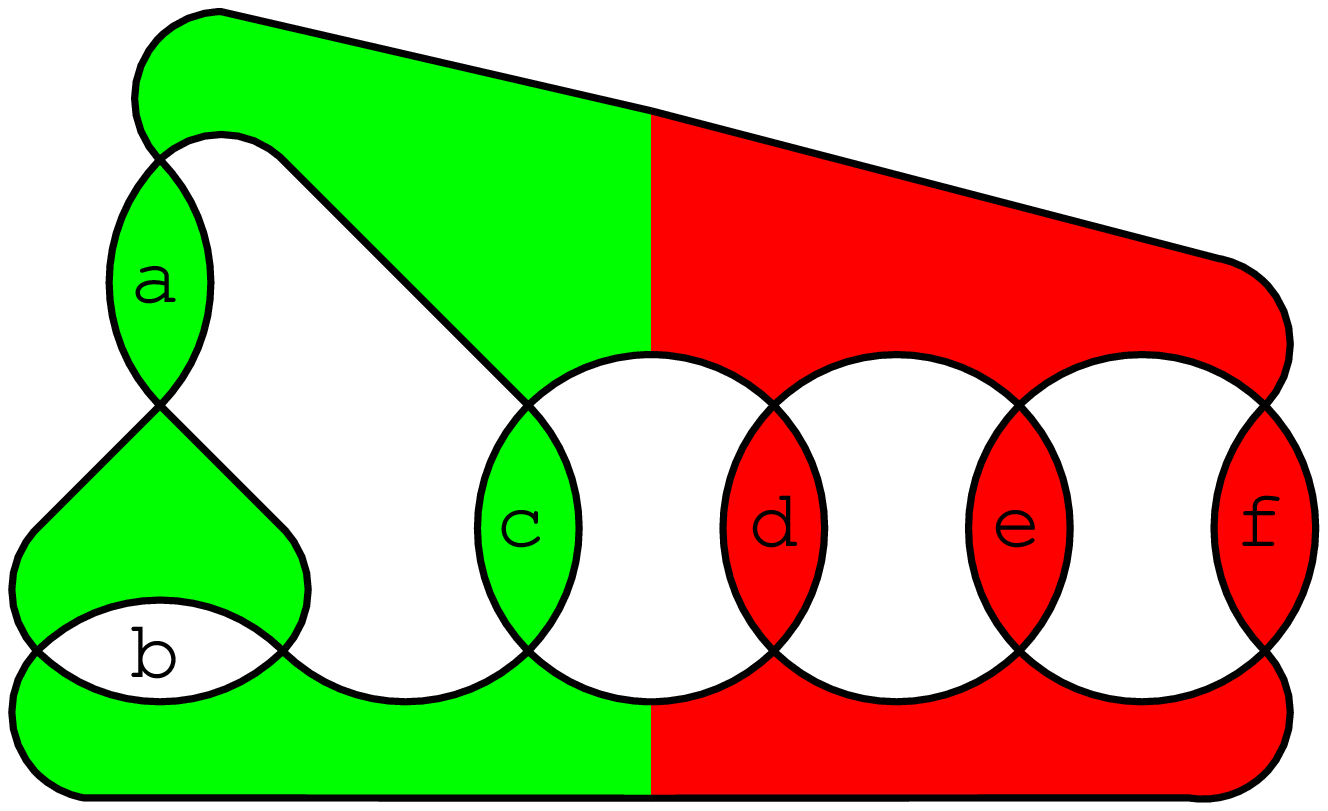}} \scalebox{0.40}{\includegraphics{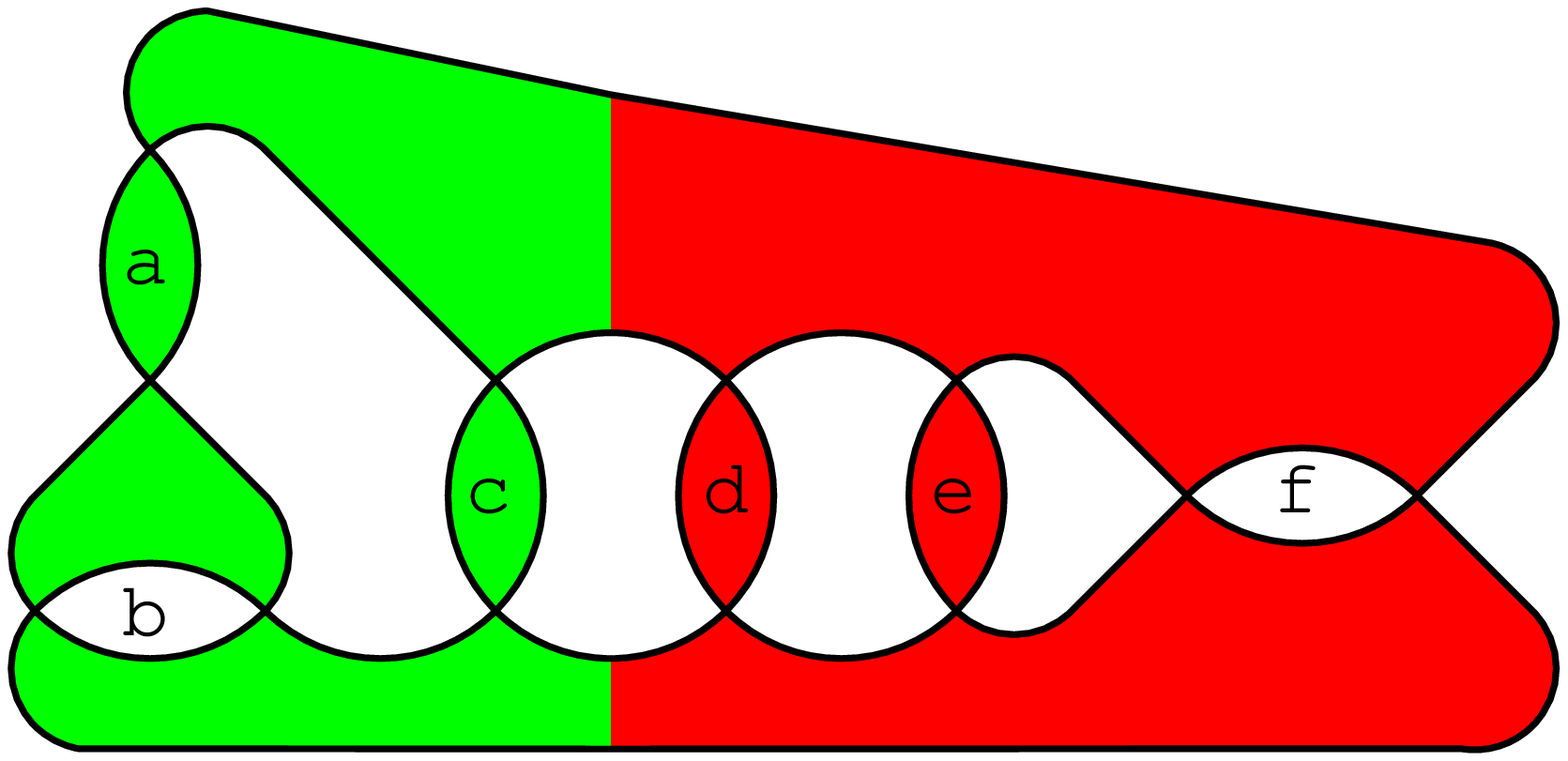}}

\caption{Two families of knots formed by six conways with the same seed. Its C-functions are
$(1 + a_2 (a_1 + a_3)) a_4 a_5 a_6 + (1 + a_1 a_2) a_3 (a_4 a_5 + a_5 a_6 + a_6 a_4)$ and
$(1 + a_2 (a_1 + a_3)) a_4 a_5 + (1 + a_1 a_2) a_3 (a_4 + a_5 + a_4 a_5 a_6)$.}
\end{figure}

\begin{figure}

FOUR FAMILIES OF KNOTS FORMED BY SIX CONWAYS CORRESPONDING TO THE LINK $6_1^3$ AS SEED

\psfrag{a}{\LARGE{$a_2$}}
\psfrag{b}{\LARGE{$a_1$}}
\psfrag{c}{\LARGE{$a_4$}}
\psfrag{d}{\LARGE{$a_3$}}
\psfrag{e}{\LARGE{$a_6$}}
\psfrag{f}{\LARGE{$a_5$}}

\scalebox{0.40}{\includegraphics{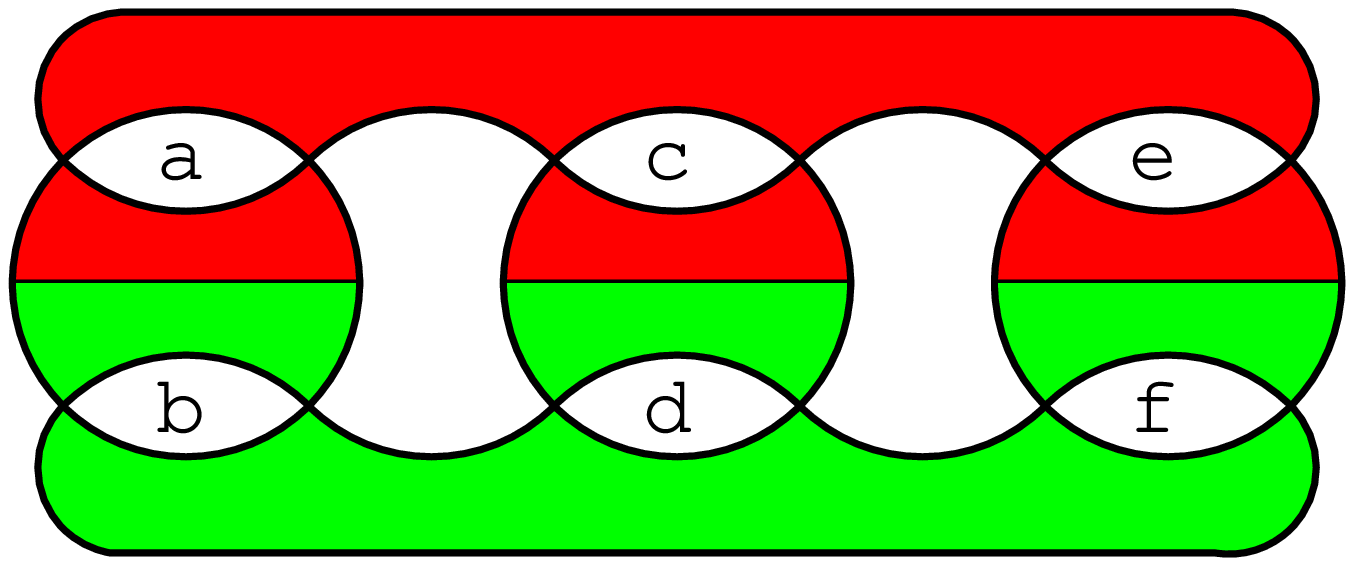}} \quad \quad \scalebox{0.40}{\includegraphics{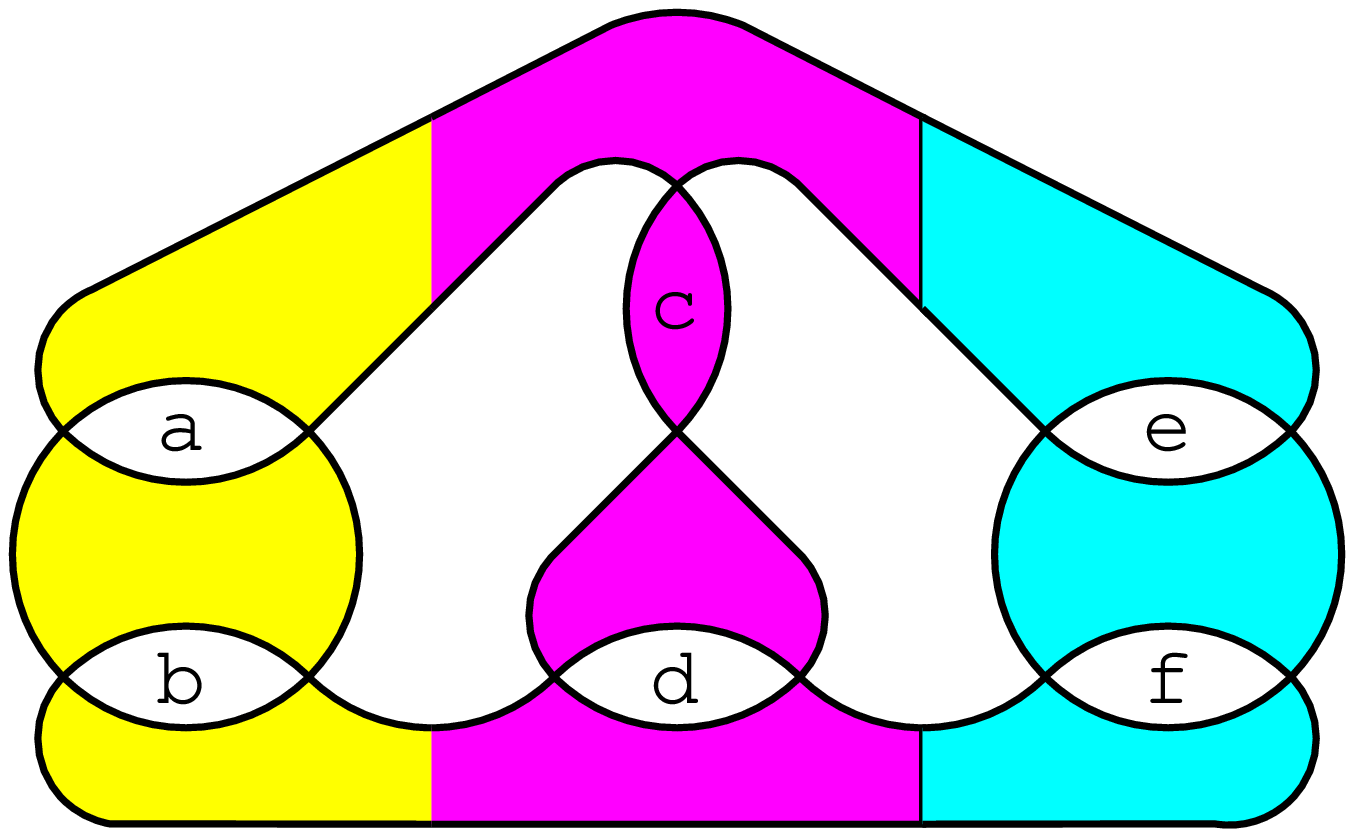}}

\caption{C-functions with twelve monomials. The C-function of the first is
$a_1 a_3 a_5 ( a_2 + a_4 + a_6) + (a_1 + a_3 + a_5) a_2 a_4 a_6 + a_3 a_5 a_2 (a_6 + a_4) + a_5 a_1 a_4 (a_6 + a_2) + a_1 a_3 a_6 (a_2 + a_4)$
and the C-function of the second is $(a_1 + a_2) (a_3 a_4 + 1) a_5 a_6 + a_1 a_2 (a_3 a_4 + 1) (a_5 + a_6) + (a_1 + a_2) a_3 (a_5 + a_6)$}
\end{figure}

\begin{figure}

\psfrag{a}{\LARGE{$a_2$}}
\psfrag{b}{\LARGE{$a_1$}}
\psfrag{c}{\LARGE{$a_4$}}
\psfrag{d}{\LARGE{$a_3$}}
\psfrag{e}{\LARGE{$a_6$}}
\psfrag{f}{\LARGE{$a_5$}}

\scalebox{0.40}{\includegraphics{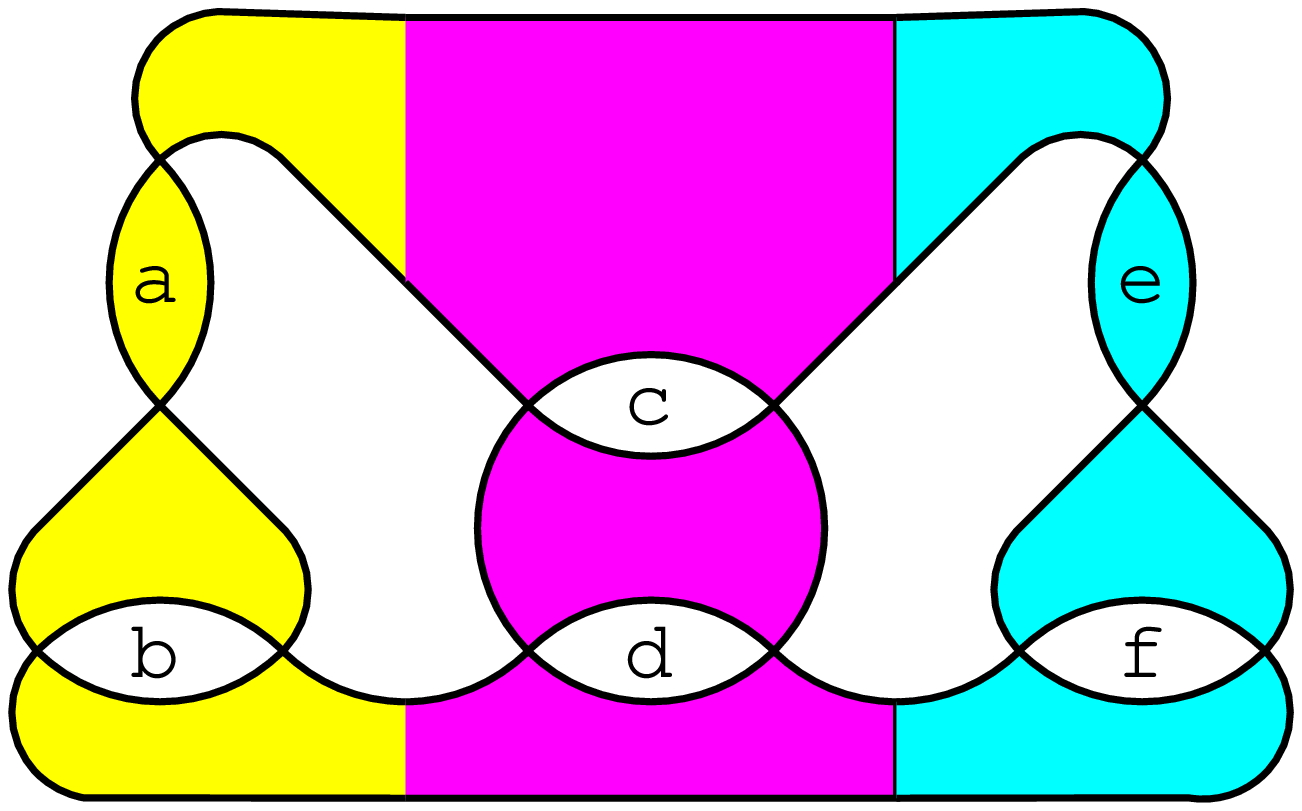}} \quad \scalebox{0.40}{\includegraphics{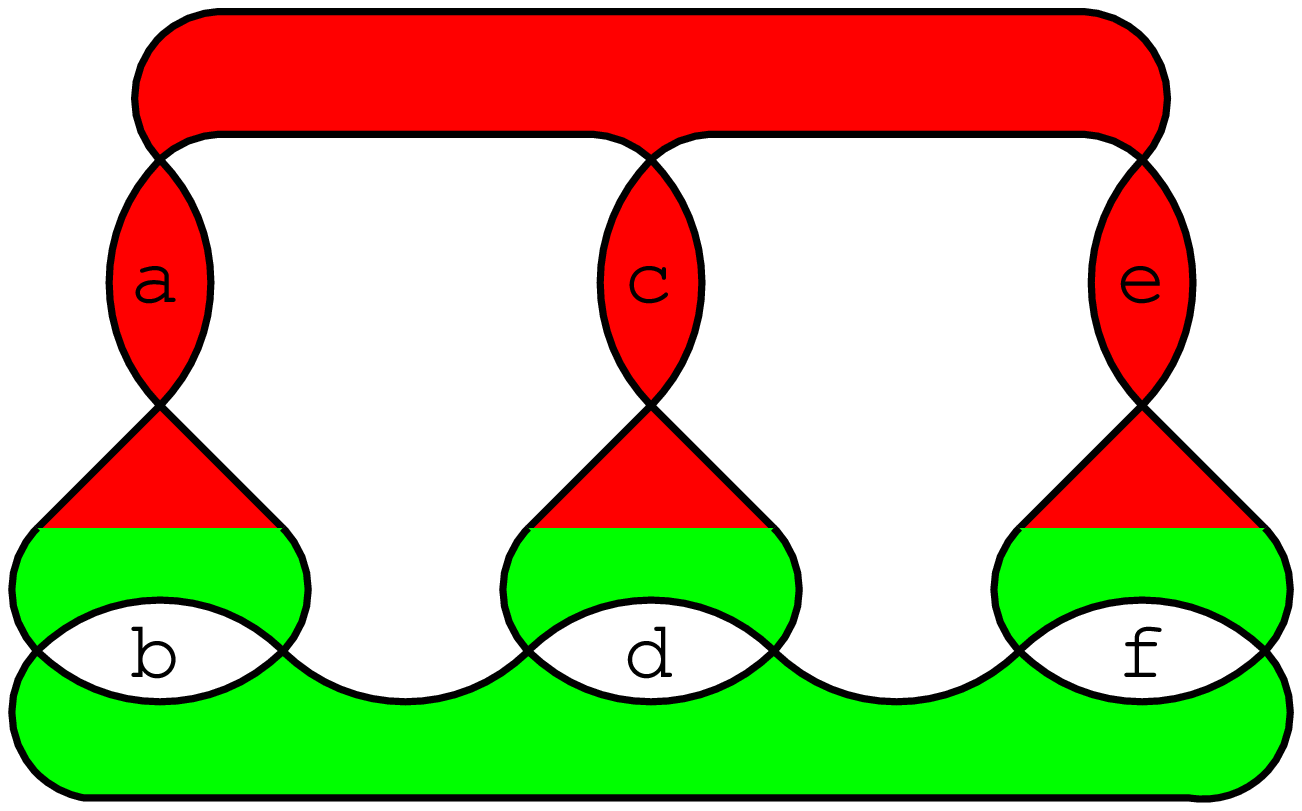}}

\caption{The other Conway's functions with twelve terms. The C-function of the first is
$(a_1 + a_2) (a_3 + a_4) a_5 a_6 + a_1 a_2 (a_3 + a_4) (a_5 + a_6) + (a_1 + a_2) a_3 a_4 (a_5 + a_6)$ and the C-function of the second is
$a_1 a_3 a_5 (a_2 a_4 + a_4 a_6 + a_6 a_2) + a_3 a_5 (a_4 + a_6) + a_5 a_1 (a_6 + a_2) + a_1 a_3 (a_2 + a_4) + a_1 + a_3 + a_5$.}

\end{figure}

\begin{figure}
FAMILIES FORMED BY SIX CONWAYS ASSOCIATED TO THE LINK $6_2^2$ AS SEED, THREE CASES
\psfrag{a}{\LARGE{$a_1$}}
\psfrag{b}{\LARGE{$a_2$}}
\psfrag{c}{\LARGE{$a_3$}}
\psfrag{d}{\LARGE{$a_4$}}
\psfrag{e}{\LARGE{$a_5$}}
\psfrag{f}{\LARGE{$a_6$}}

\hfil \scalebox{0.40}{\includegraphics{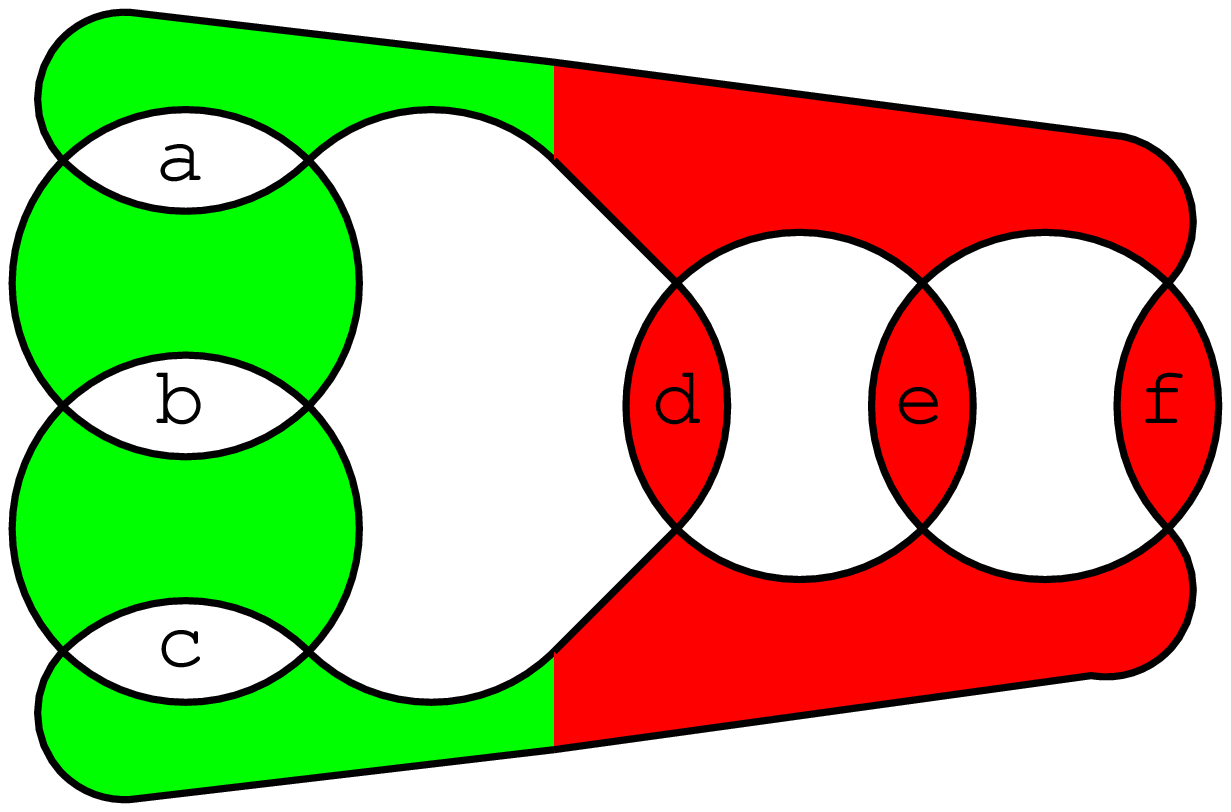}} \quad \scalebox{0.40}{\includegraphics{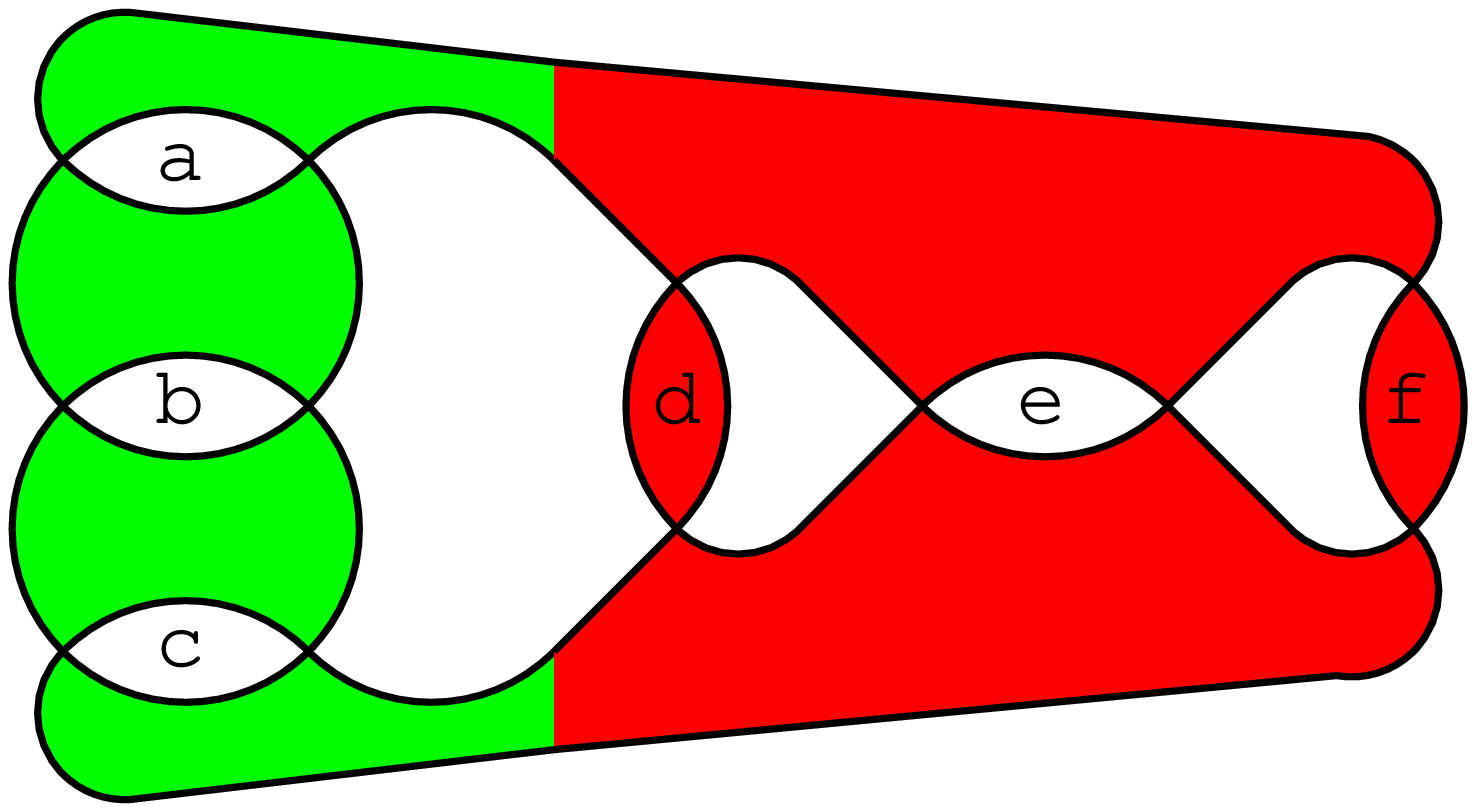}}\hfil

\caption{Families formed by six conways. Its C-functions with ten terms are $(a_1 a_2 + a_2 a_3 + a_3 a_1)(a_4 a_5 + a_5 a_6 + a_6 a_4) + a_1 a_2 a_3 a_4 a_5 a_6$ and $(a_1 a_2 + a_2 a_3 + a_3 a_1) (a_4 a_5 a_6 + a_4 + a_6) + a_1 a_2 a_3 a_4 a_6$.}

\end{figure}

\begin{figure}

\psfrag{a}{\LARGE{$a_1$}}
\psfrag{b}{\LARGE{$a_2$}}
\psfrag{c}{\LARGE{$a_3$}}
\psfrag{d}{\LARGE{$a_4$}}
\psfrag{e}{\LARGE{$a_5$}}
\psfrag{f}{\LARGE{$a_6$}}

\hfil \scalebox{0.40}{\includegraphics{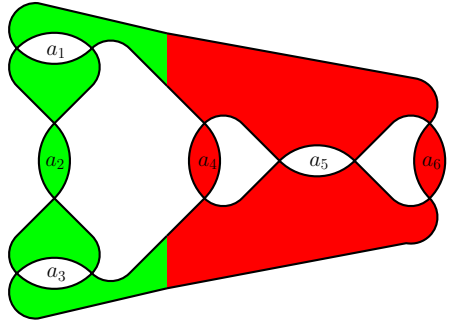}} \hfil

\caption{The third family formed by six conways with C-function of ten terms $(a_1 a_2 a_3 + a_1 + a_3) (a_4 a_5 a_6 + a_4 + a_6) + a_1 a_3 a_4 a_6$.}

\end{figure}

\begin{figure}
FAMILIES FORMED BY SIX CONWAYS WITH THE LINK $6_3^2$ AS SEED, SIX CASES
\psfrag{a}{\LARGE{$a_1$}}
\psfrag{b}{\LARGE{$a_2$}}
\psfrag{c}{\LARGE{$a_3$}}
\psfrag{d}{\LARGE{$a_4$}}
\psfrag{e}{\LARGE{$a_5$}}
\psfrag{f}{\LARGE{$a_6$}}

\hfil \scalebox{0.40}{\includegraphics{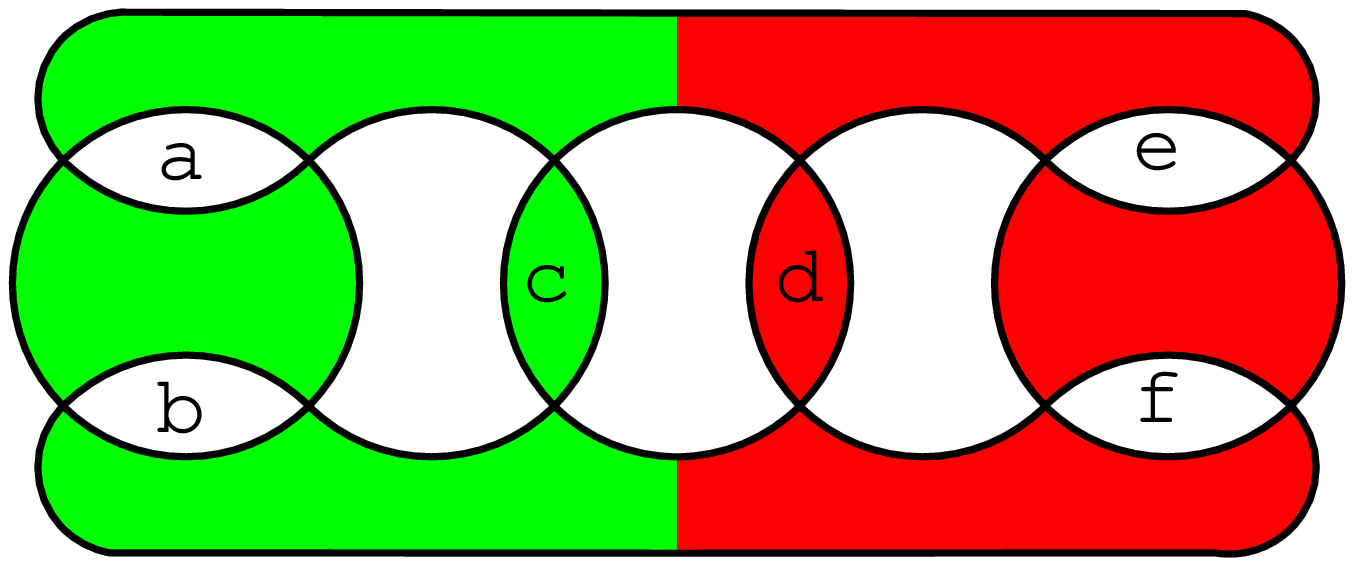}} \quad \scalebox{0.40}{\includegraphics{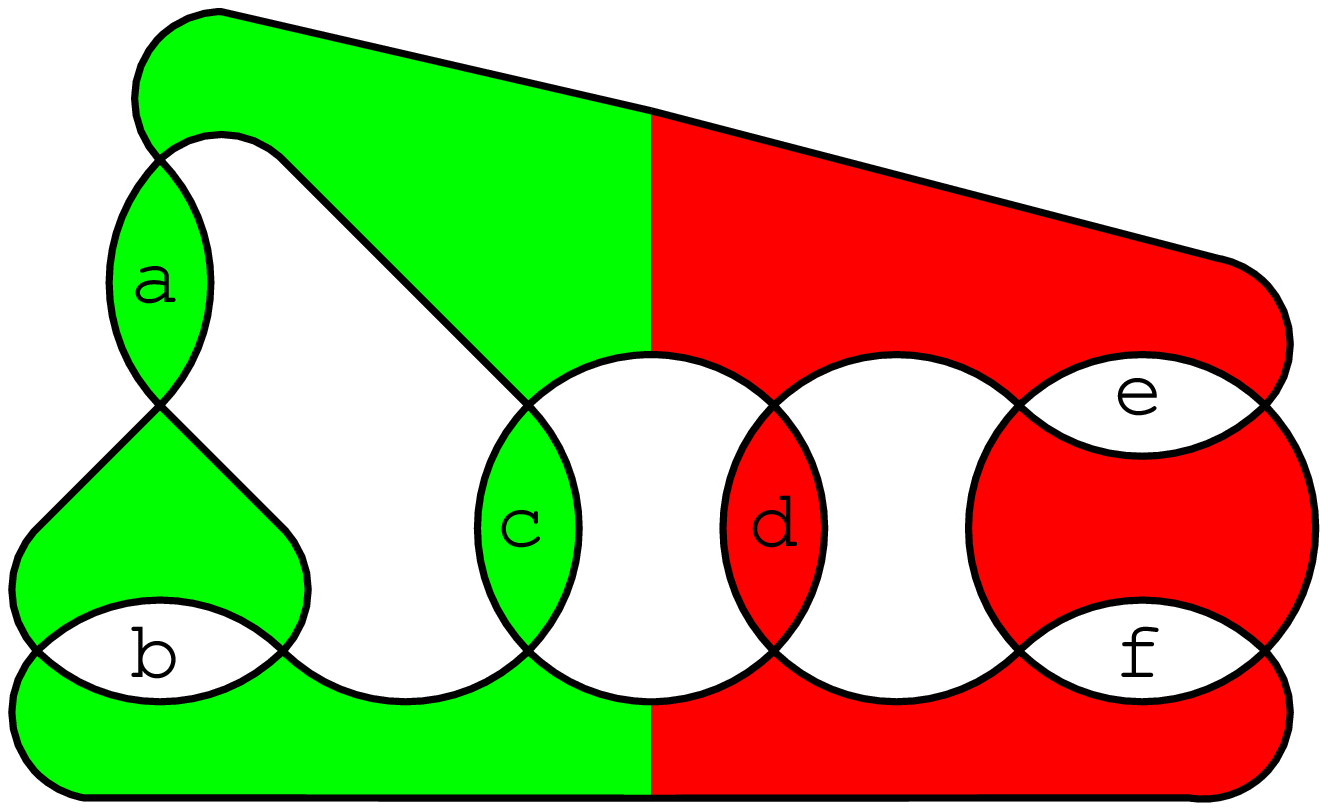}} \hfil

\caption{First two families formed by six conways. Its C-functions with twelve terms are $a_3 (a_1 + a_2)(a_4 a_5 a_6 + a_5 + a_6) + (a_1 a_2 a_3 + a_1 + a_2) a_4 (a_5 + a_6)$
and $a_3 (1 + a_1 a_2)(a_4 a_5 a_6 + a_5 + a_6) + ¨(1 + (a_1 + a_3) a_2) a_4 (a_5 + a_6)$.}

\end{figure}

\begin{figure}

\psfrag{a}{\LARGE{$a_1$}}
\psfrag{b}{\LARGE{$a_2$}}
\psfrag{c}{\LARGE{$a_3$}}
\psfrag{d}{\LARGE{$a_4$}}
\psfrag{e}{\LARGE{$a_5$}}
\psfrag{f}{\LARGE{$a_6$}}

\hfil \scalebox{0.40}{\includegraphics{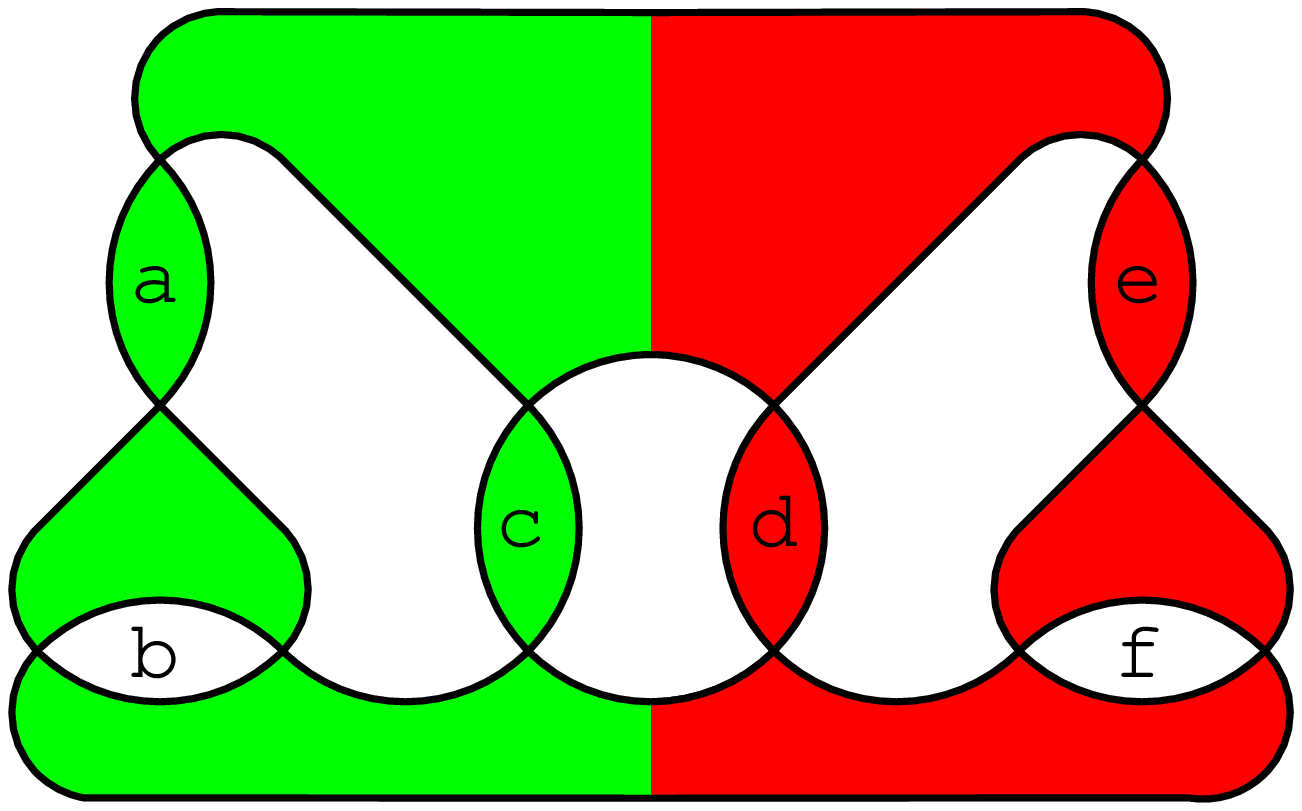}} \quad \scalebox{0.40}{\includegraphics{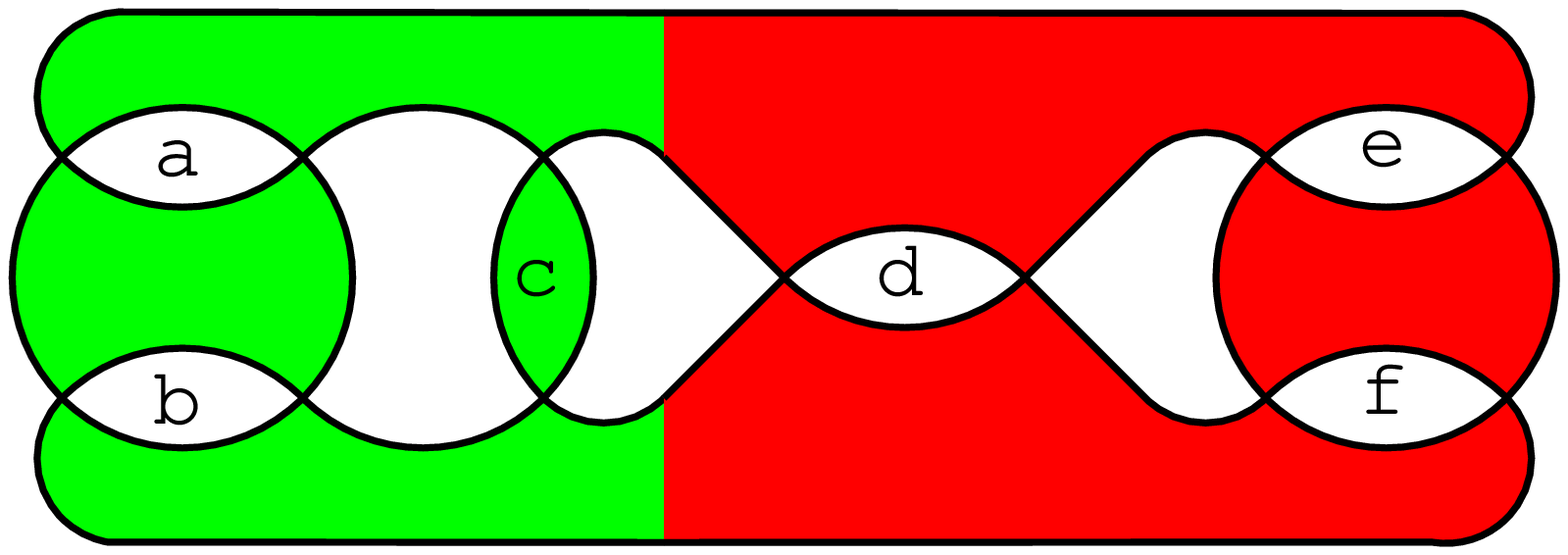}} \hfil

\caption{Two more cases of families formed by six conways. Its C-functions with twelve terms are $a_3 (1 + a_1 a_2)[a_6 (a_4 + a_5) + 1] + [a_2 (a_1 + a_3) + 1 ] a_4 (a_5 a_6 + 1)$
and $(a_1 + a_2) a_3 (a_4 a_5 + a_5 a_6 + a_6 a_4) + (a_1 a_2 a_3 + a_1 + a_2) (a_5 + a_6)$.}

\end{figure}

\begin{figure}

\psfrag{a}{\LARGE{$a_1$}}
\psfrag{b}{\LARGE{$a_2$}}
\psfrag{c}{\LARGE{$a_3$}}
\psfrag{d}{\LARGE{$a_4$}}
\psfrag{e}{\LARGE{$a_5$}}
\psfrag{f}{\LARGE{$a_6$}}

\scalebox{0.40}{\includegraphics{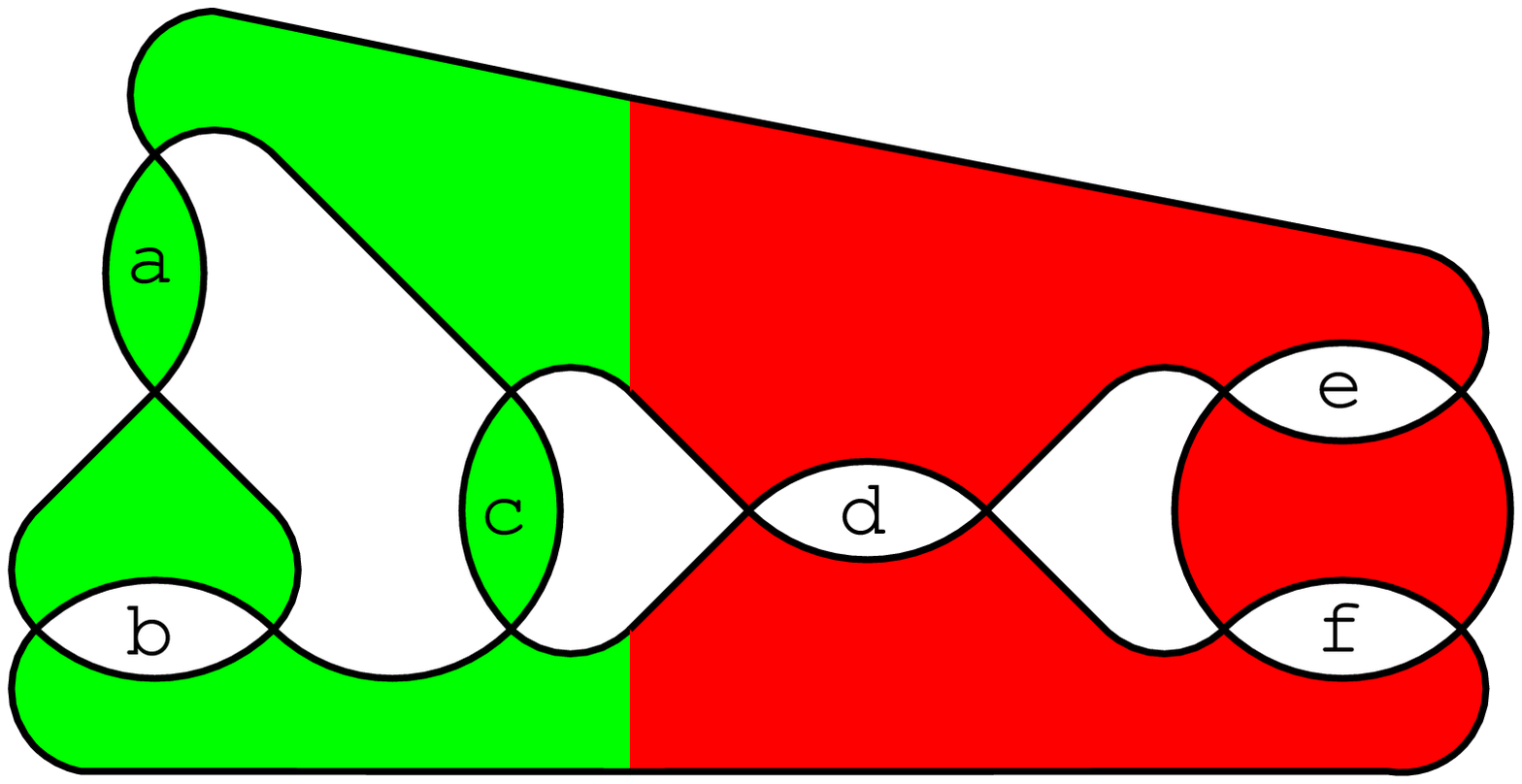}} \quad \scalebox{0.40}{\includegraphics{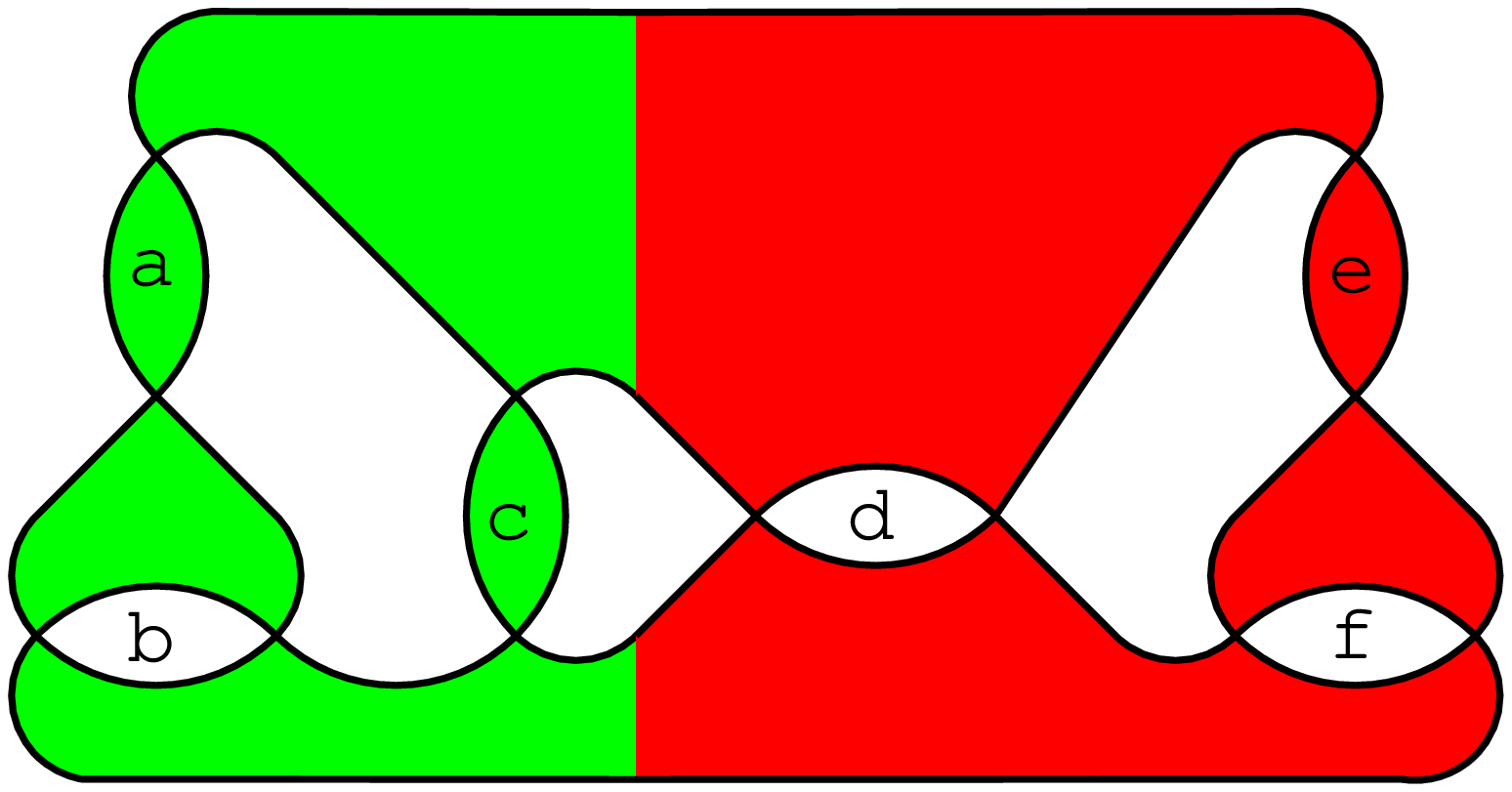}}

\caption{The two other families formed by six conways with C-functions of twelve terms $a_3 (1 + a_1 a_2) (a_4 a_5 + a_5 a_6 + a_6 a_4) + [a_2 (a_1 + a_3) + 1 ] (a_5 + a_6)$
and $(1 + a_1 a_2) a_3 (a_4 a_5 a_6 + a_4 + a_6) + [1 + (a_1 + a_3) a_2] (a_5 a_6 + 1)$.}

\end{figure} 

\begin{figure}
FAMILIES FORMED BY SIX CONWAYS WITH THE KNOT $6_2$ AS SEED, EIGHT CASES
\psfrag{a}{\LARGE{$a_1$}}
\psfrag{b}{\LARGE{$a_2$}}
\psfrag{c}{\LARGE{$a_3$}}
\psfrag{d}{\LARGE{$a_4$}}
\psfrag{e}{\LARGE{$a_5$}}
\psfrag{f}{\LARGE{$a_6$}}

\scalebox{0.40}{\includegraphics{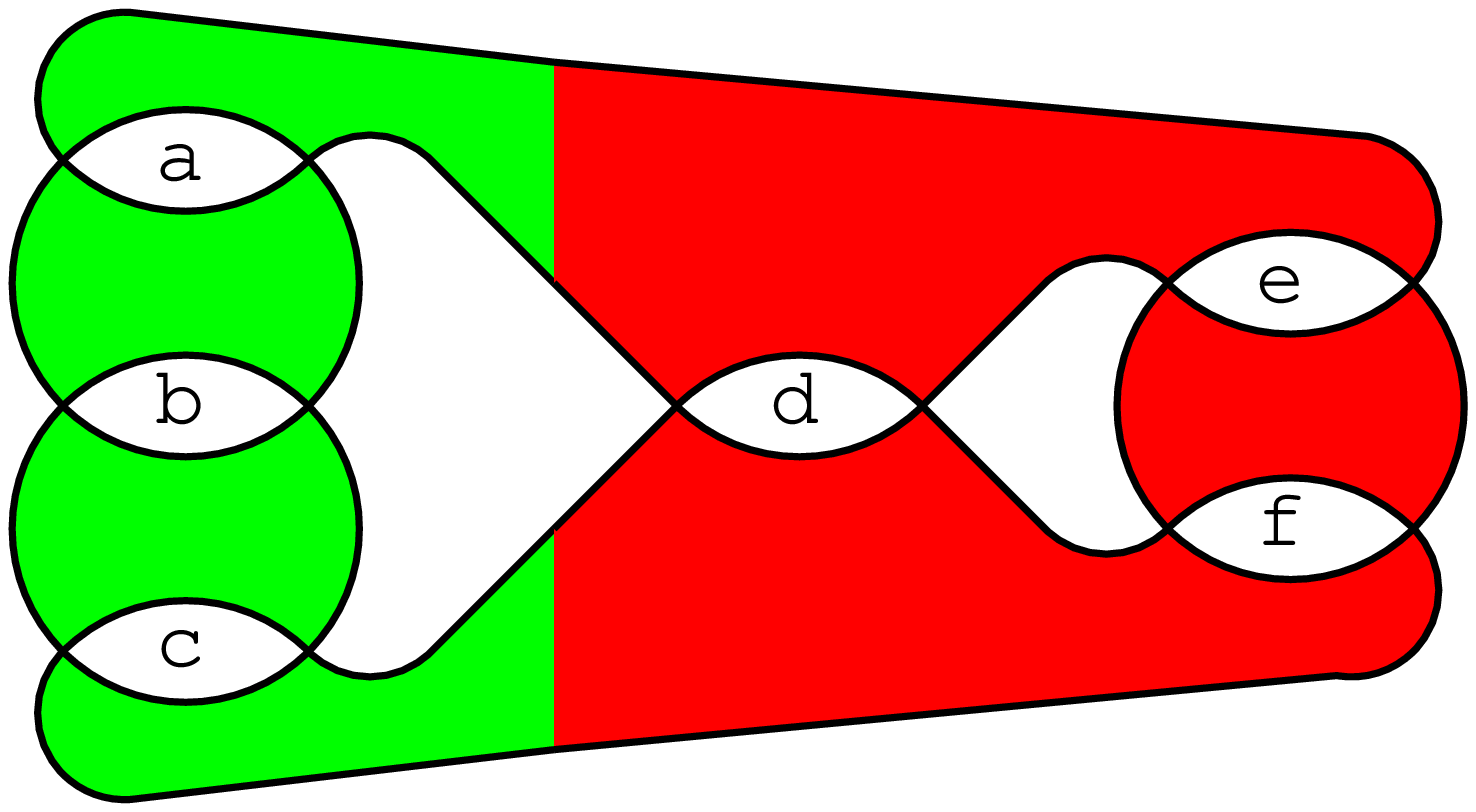}} \quad \scalebox{0.40}{\includegraphics{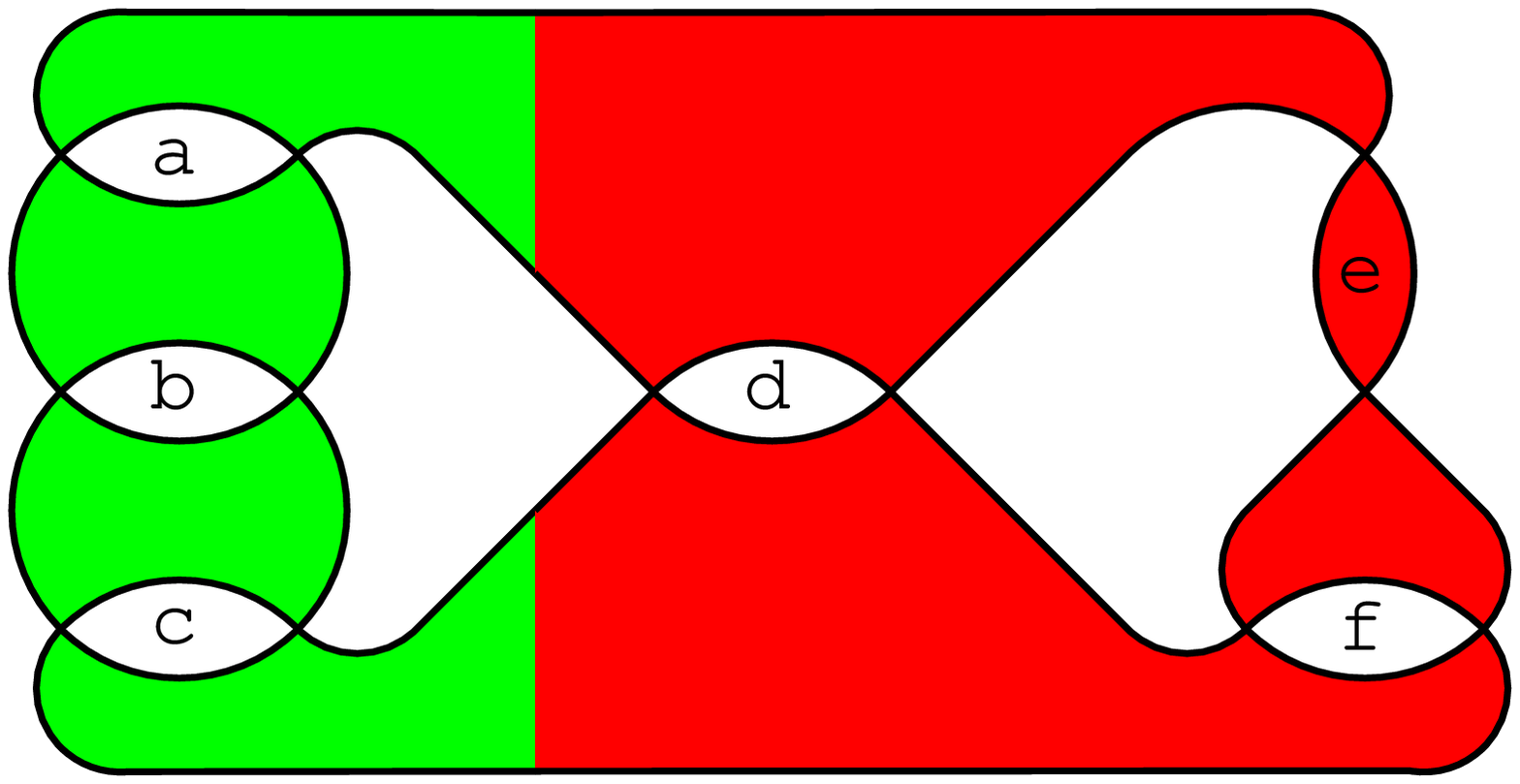}}

\caption{Families formed by six conways with C-functions with eleven terms $(a_1 a_2 + a_2 a_3 + a_3 a_1) (a_4 a_5 + a_5 a_6 + a_6 a_4) + a_1 a_2 a_3 (a_5 + a_6)$
and $(a_1 a_2 + a_2 a_3 + a_3 a_1) (a_4 a_5 a_6 + a_4 + a_6) + a_1 a_2 a_3 (a_5 a_6 + 1)$}

\end{figure}

\begin{figure}

\psfrag{a}{\LARGE{$a_1$}}
\psfrag{b}{\LARGE{$a_2$}}
\psfrag{c}{\LARGE{$a_3$}}
\psfrag{d}{\LARGE{$a_4$}}
\psfrag{e}{\LARGE{$a_5$}}
\psfrag{f}{\LARGE{$a_6$}}

\scalebox{0.40}{\includegraphics{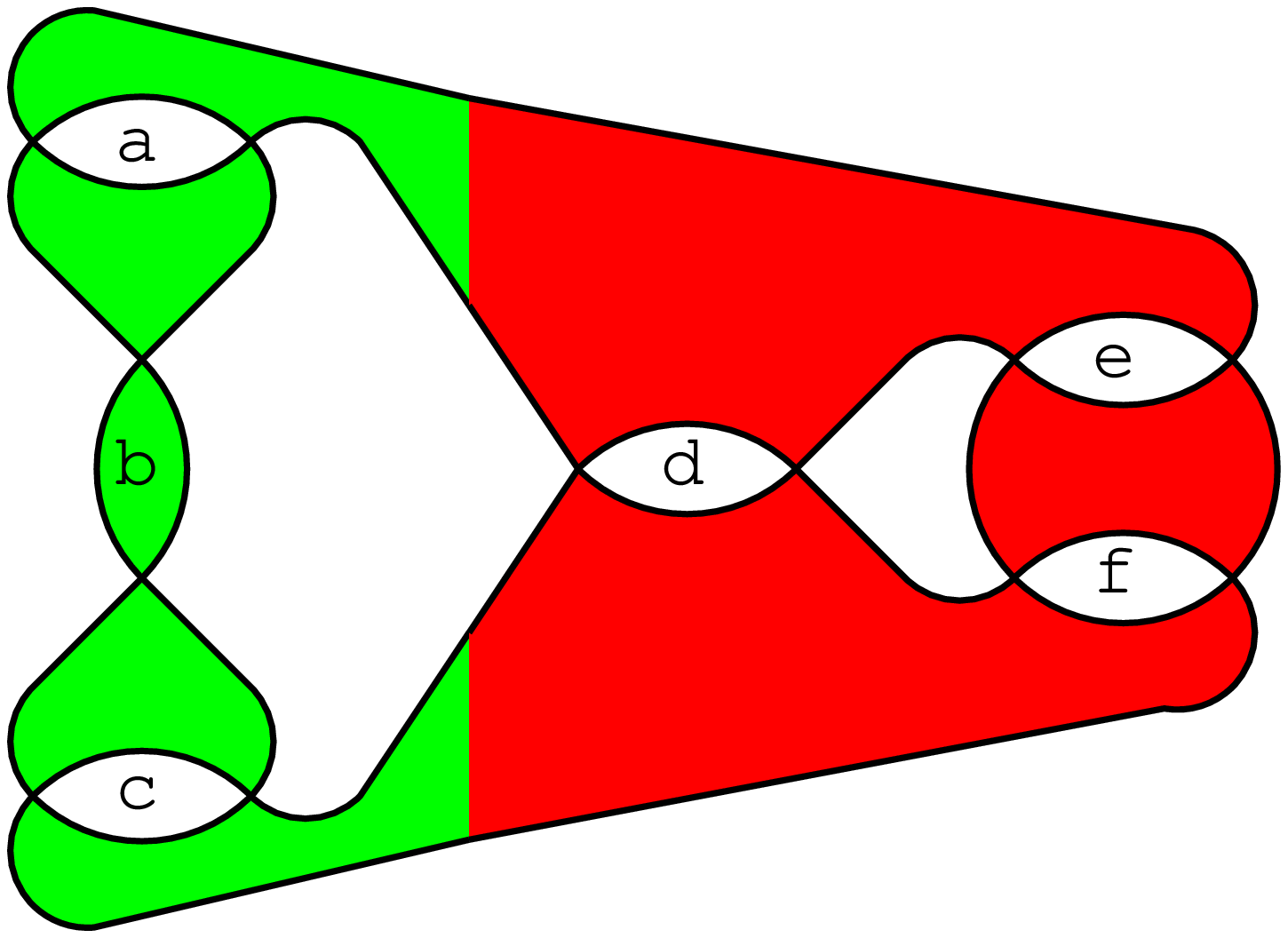}} \quad \scalebox{0.40}{\includegraphics{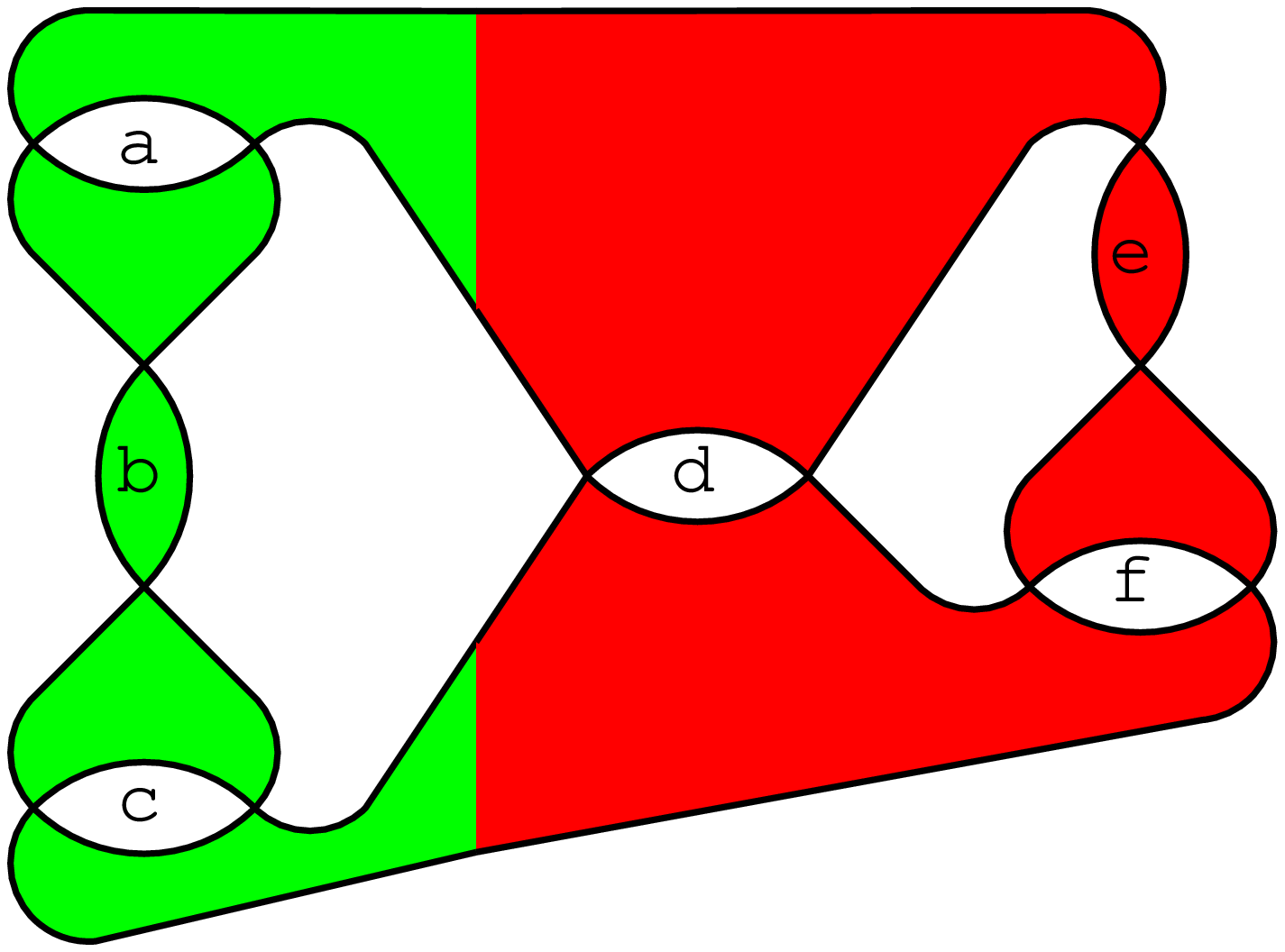}}

\caption{Families formed by six conways with C-functions of eleven terms $(a_1 a_2 a_3 + a_1 + a_3) (a_4 a_5 + a_5 a_6 + a_6 a_4) + a_1 a_3 (a_5 + a_6)$ and $(a_1 a_2 a_3 + a_1 + a_3) (a_4 a_5 a_6 + a_4 + a_6) + a_1 a_3 (a_5 a_6 + 1)$}

\end{figure}

\begin{figure}

\psfrag{a}{\LARGE{$a_1$}}
\psfrag{b}{\LARGE{$a_2$}}
\psfrag{c}{\LARGE{$a_3$}}
\psfrag{d}{\LARGE{$a_4$}}
\psfrag{e}{\LARGE{$a_5$}}
\psfrag{f}{\LARGE{$a_6$}}

\scalebox{0.40}{\includegraphics{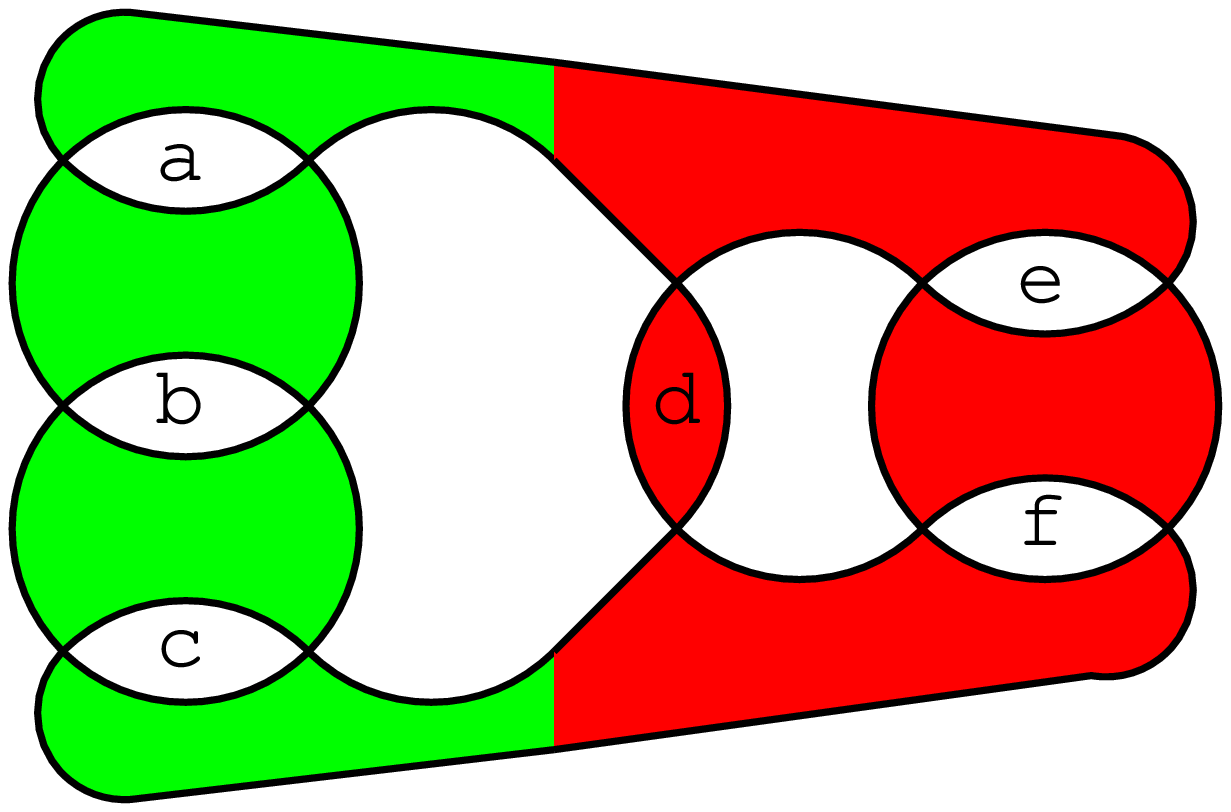}} \quad \scalebox{0.40}{\includegraphics{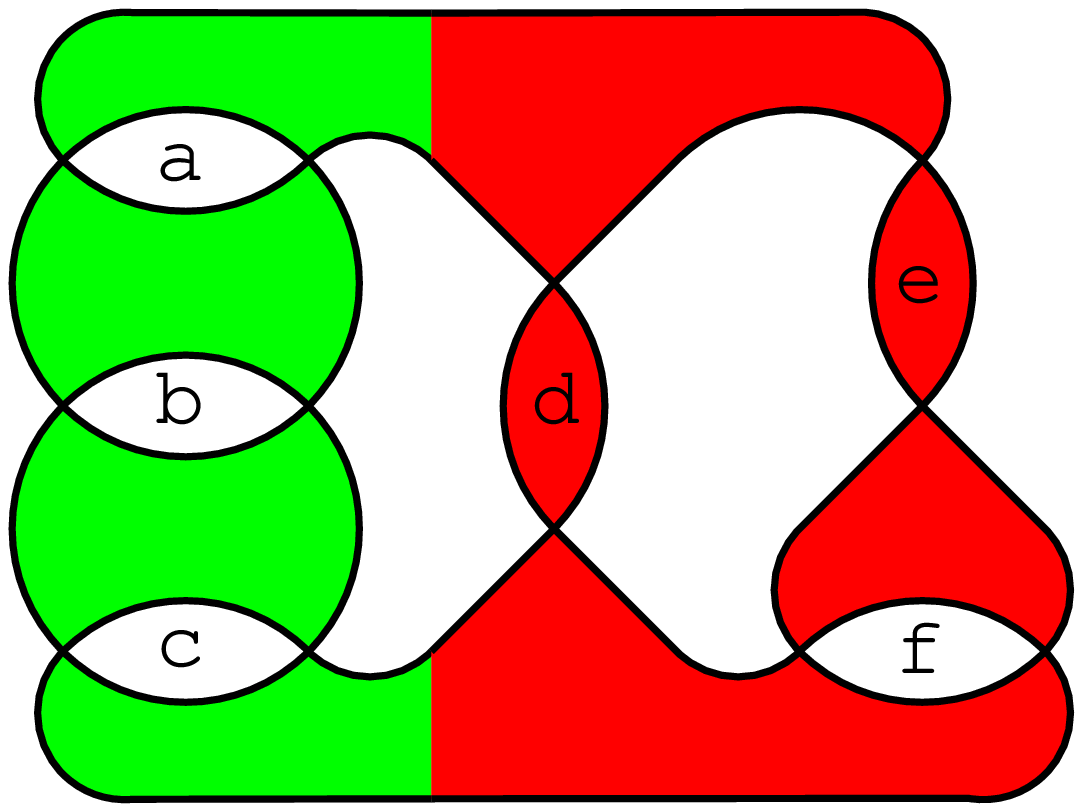}}

\caption{Families formed by six conways with C-functions with eleven terms $(a_1 a_2 + a_2 a_3 + a_3 a_1) (a_4 a_5 a_6 + a_5 + a_6) + a_1 a_2 a_3 a_4 (a_5 + a_6)$
and $(a_1 a_2 + a_2 a_3 + a_3 a_1) ((a_4 + a_5) a_6 + 1) + a_1 a_2 a_3 a_4 (a_5 a_6 + 1)$}

\end{figure}

\begin{figure}

\psfrag{a}{\LARGE{$a_1$}}
\psfrag{b}{\LARGE{$a_2$}}
\psfrag{c}{\LARGE{$a_3$}}
\psfrag{d}{\LARGE{$a_4$}}
\psfrag{e}{\LARGE{$a_5$}}
\psfrag{f}{\LARGE{$a_6$}}

\scalebox{0.40}{\includegraphics{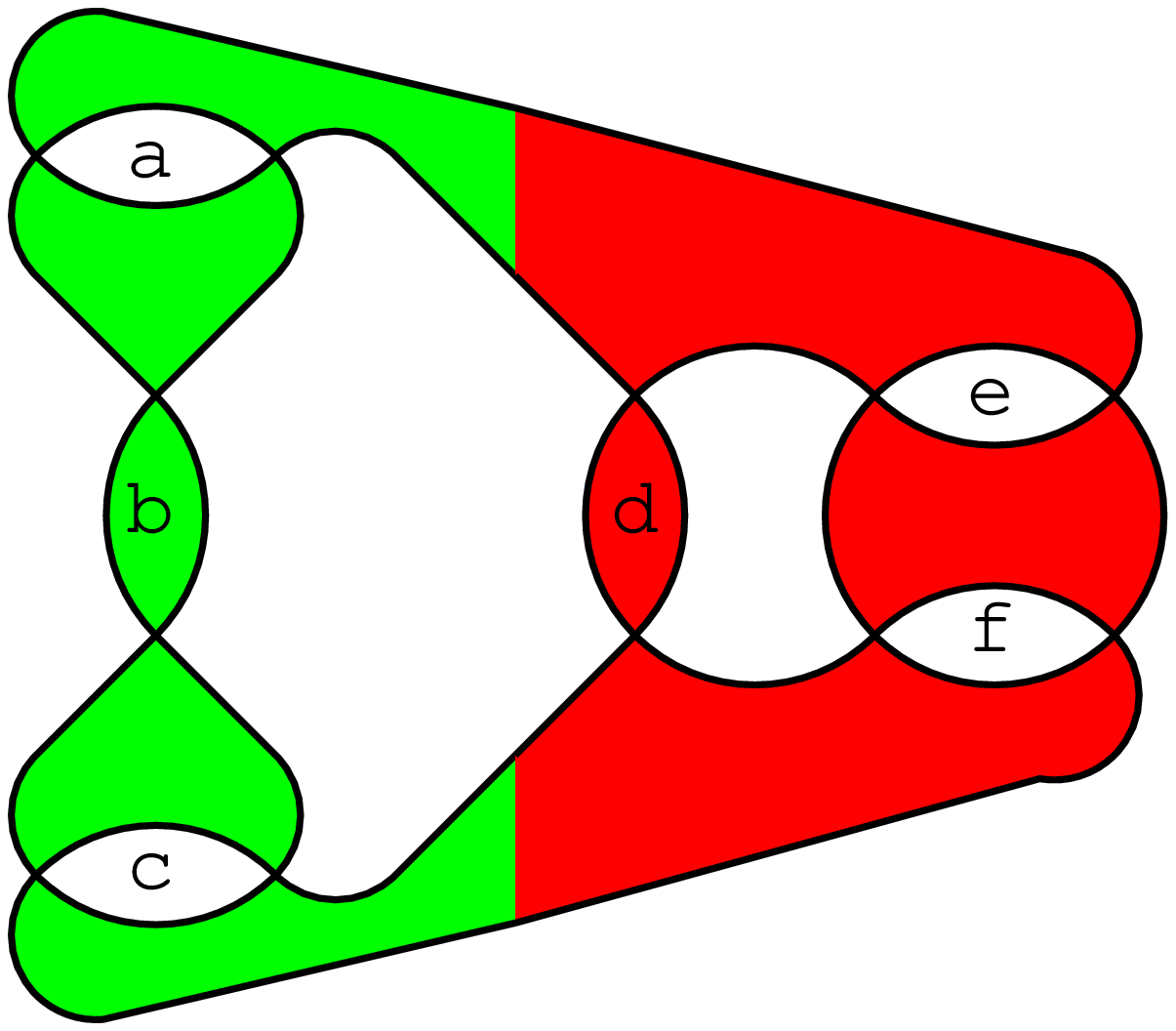}} \quad \scalebox{0.40}{\includegraphics{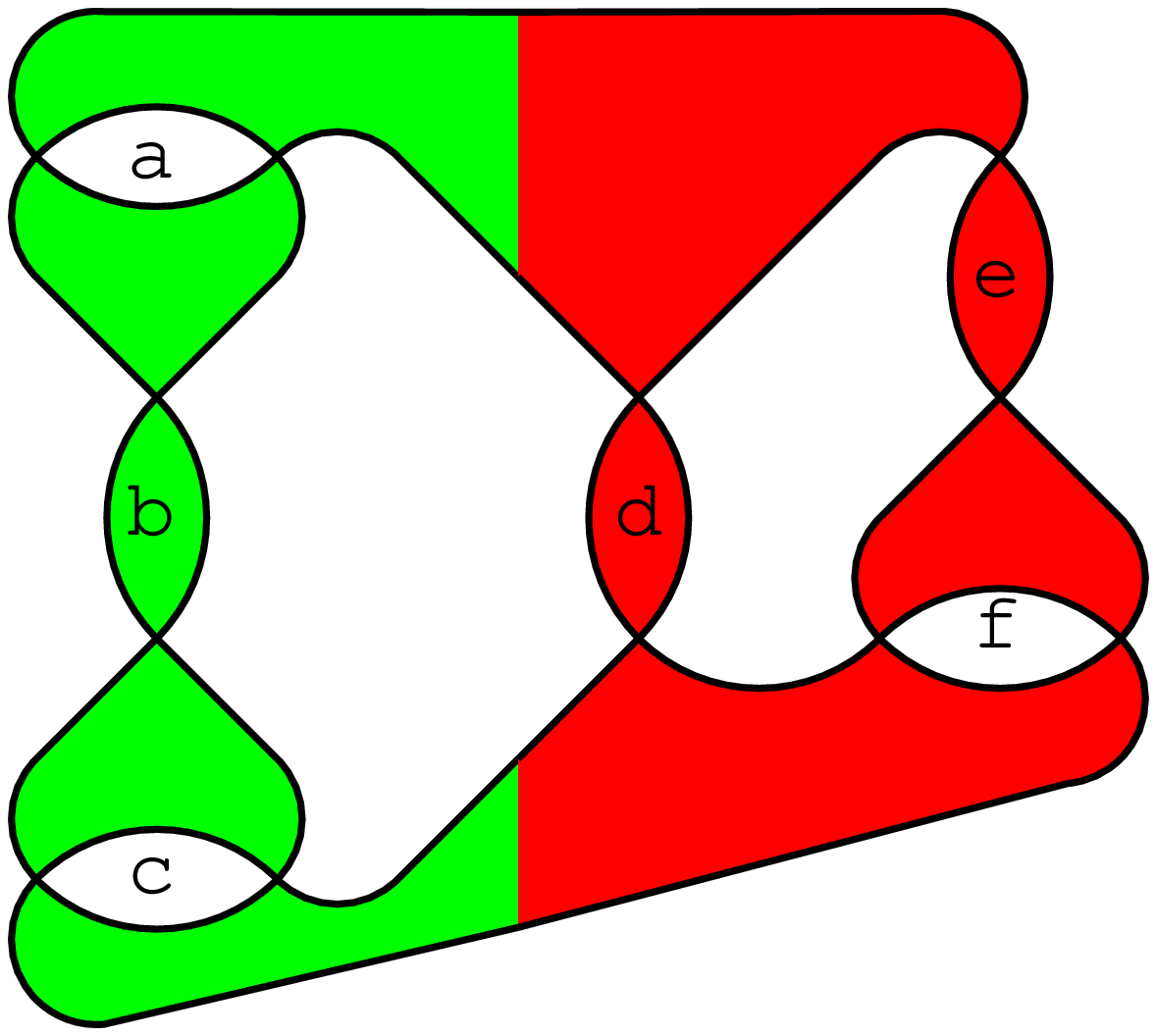}}

\caption{Families formed by six conways with C-functions of eleven terms $(a_1 a_2 a_3 + a_1 + a_3) (a_4 a_5 a_6 + a_5 + a_6) + a_1 a_3 a_4 (a_5 + a_6)$
and $(a_1 a_2 a_3 + a_1 + a_3) ((a_4 + a_5) a_6 + 1) + a_1 a_3 a_4 (a_5 a_6 + 1)$}

\end{figure}

\begin{figure}
FAMILIES FORMED BY SIX CONWAYS WITH THE KNOT $6_3$ AS SEED, TEN CASES
\psfrag{a}{\LARGE{$a_1$}}
\psfrag{b}{\LARGE{$a_2$}}
\psfrag{c}{\LARGE{$a_3$}}
\psfrag{d}{\LARGE{$a_4$}}
\psfrag{e}{\LARGE{$a_5$}}
\psfrag{f}{\LARGE{$a_6$}}

\scalebox{0.40}{\includegraphics{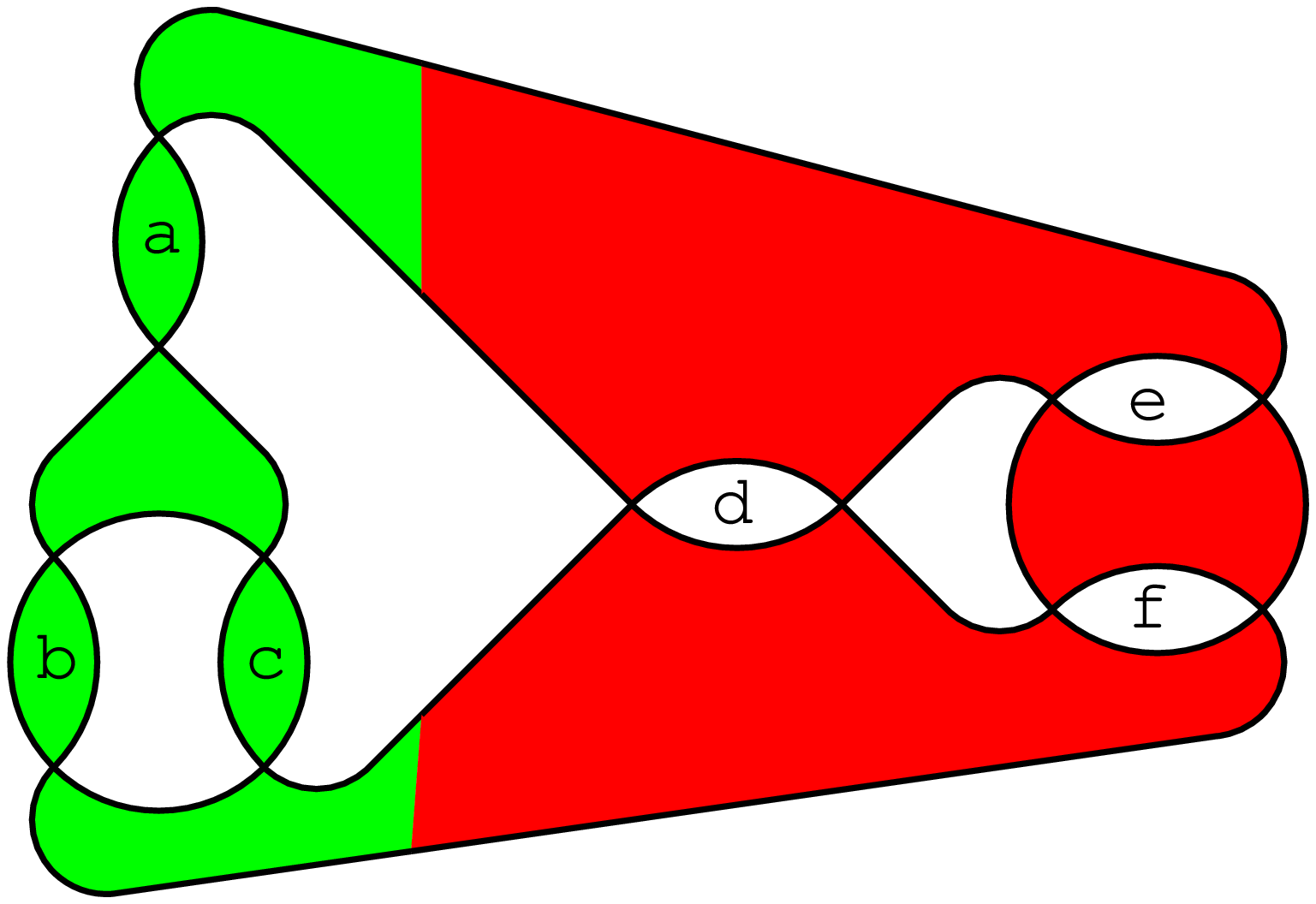}} \quad \scalebox{0.40}{\includegraphics{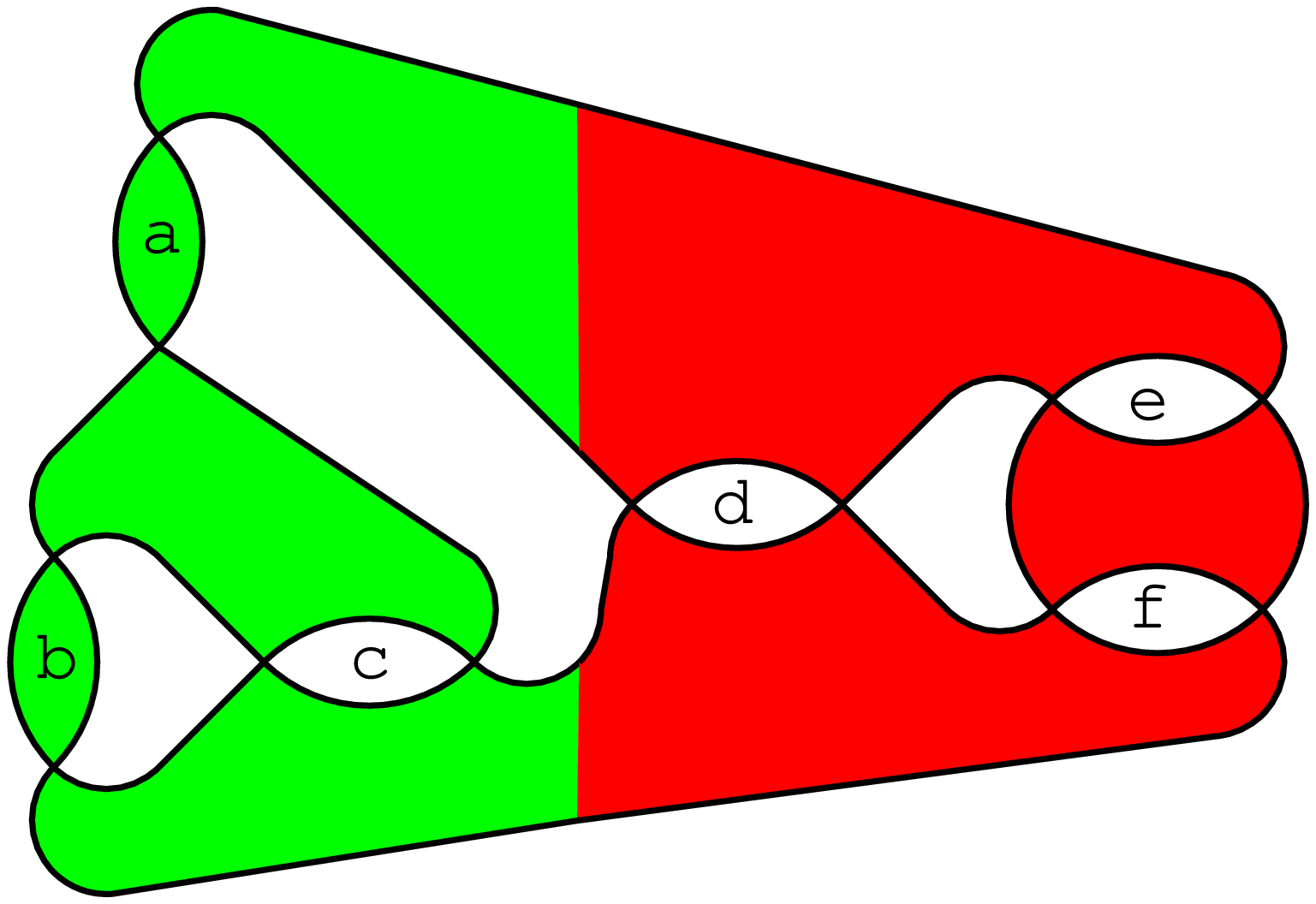}}

\caption{Families formed by six conways with C-functions formed by thirteen terms $(a_1 a_2 + a_2 a_3 + a_3 a_1) (a_4 a_5 + a_5 a_6 + a_6 a_4) +
(a_2 + a_3) (a_5 + a_6)$ and $(a_1 a_2 a_3 + a_1 + a_2) (a_4 a_5 + a_5 a_6 + a_6 a_4) + (a_2 a_3 + 1) (a_5 + a_6)$}

\end{figure}

\begin{figure}

\psfrag{a}{\LARGE{$a_1$}}
\psfrag{b}{\LARGE{$a_2$}}
\psfrag{c}{\LARGE{$a_3$}}
\psfrag{d}{\LARGE{$a_4$}}
\psfrag{e}{\LARGE{$a_5$}}
\psfrag{f}{\LARGE{$a_6$}}

\scalebox{0.40}{\includegraphics{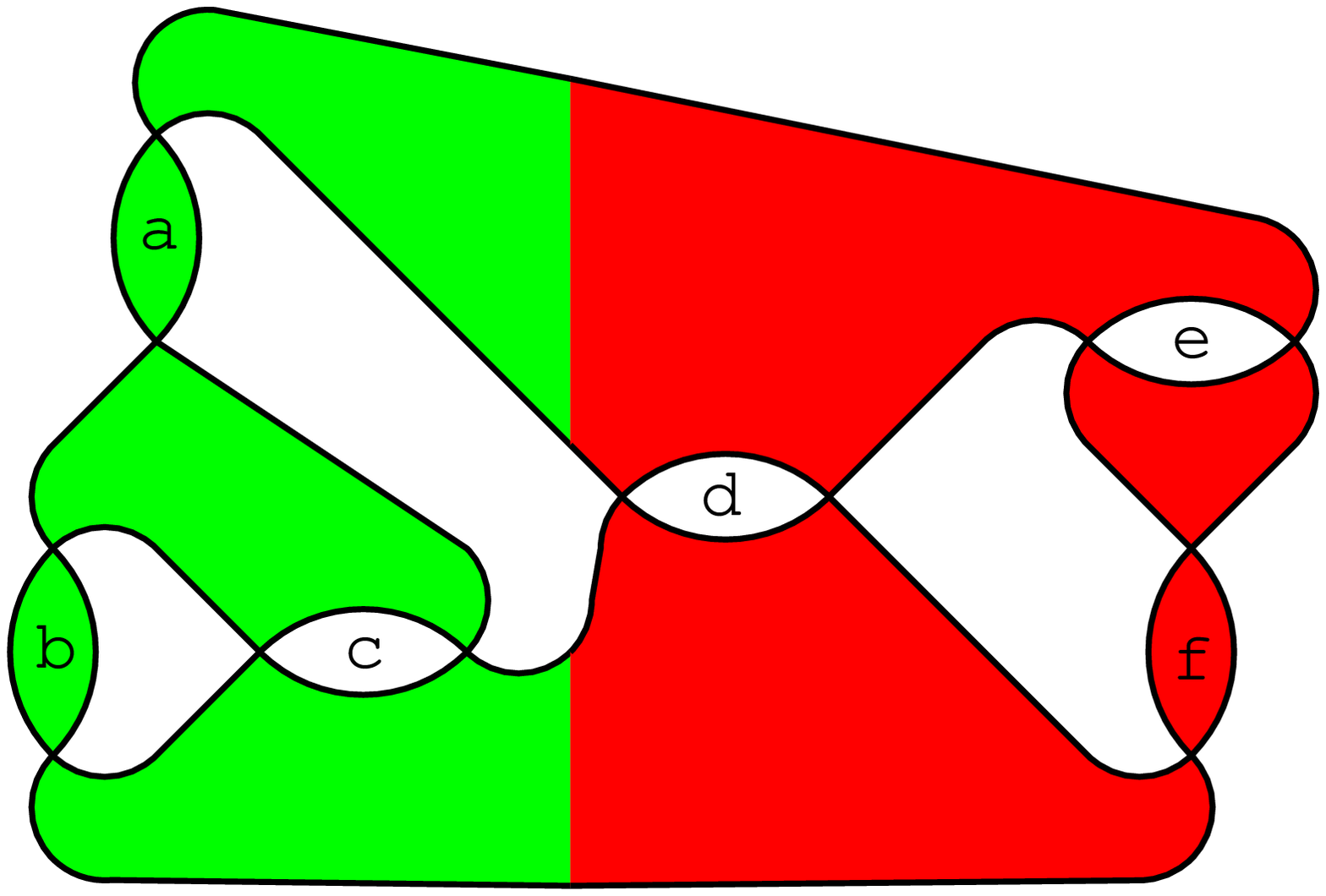}} \quad \scalebox{0.40}{\includegraphics{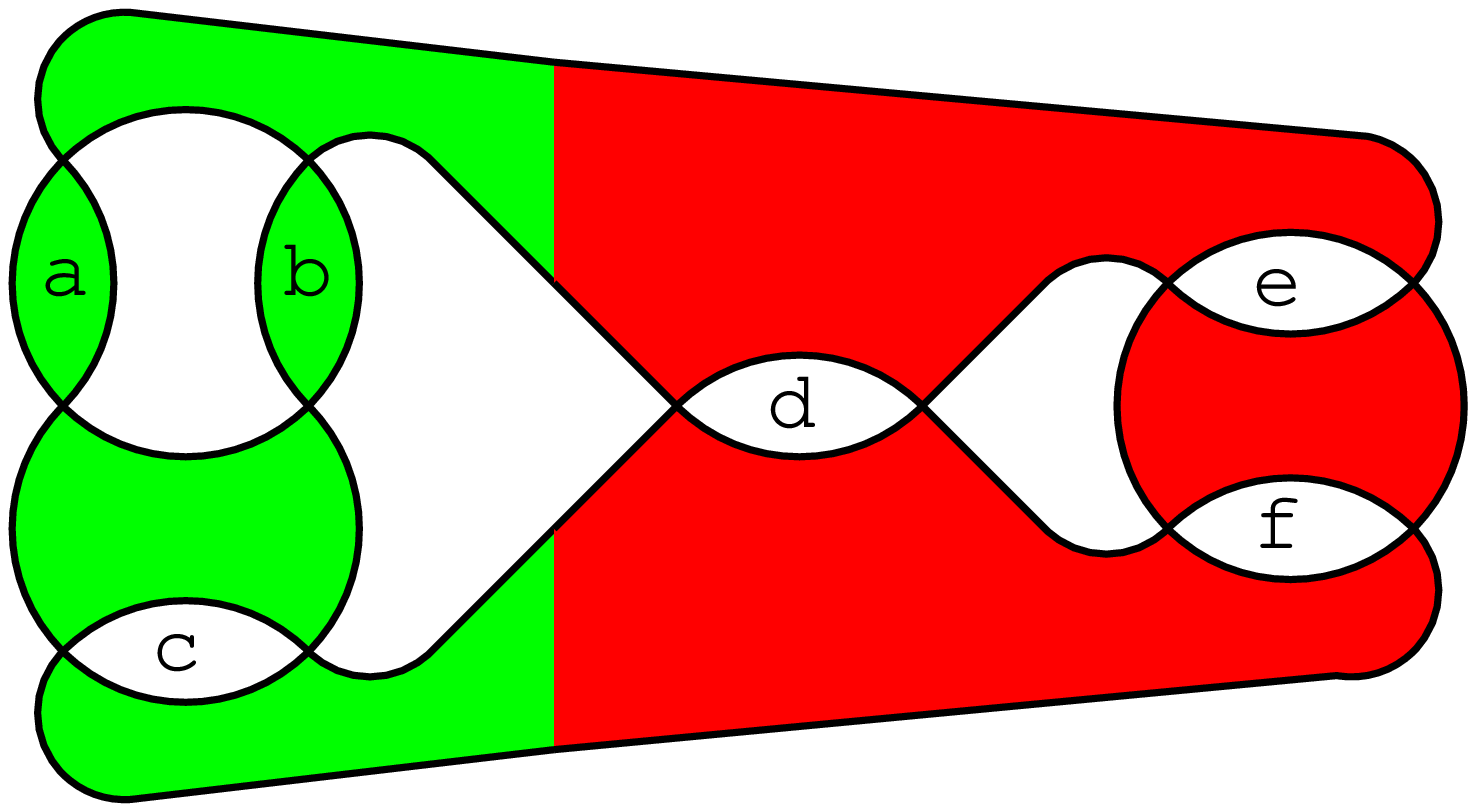}}

\caption{Families formed by six conways with its C-functions of thirteen terms $(a_1 a_2 a_3 + a_1 + a_2) (a_4 a_5 a_6 + a_4 + a_5) +
(a_2 a_3 + 1) (a_5 a_6 + 1)$
and $(a_1 a_2 a_3 + a_1 + a_2) (a_4 a_5 + a_5 a_6 + a_6 a_4) + a_3 (a_1 + a_2) (a_5 + a_6)$}

\end{figure}

\begin{figure}

\psfrag{a}{\LARGE{$a_1$}}
\psfrag{b}{\LARGE{$a_2$}}
\psfrag{c}{\LARGE{$a_3$}}
\psfrag{d}{\LARGE{$a_4$}}
\psfrag{e}{\LARGE{$a_5$}}
\psfrag{f}{\LARGE{$a_6$}}

\scalebox{0.40}{\includegraphics{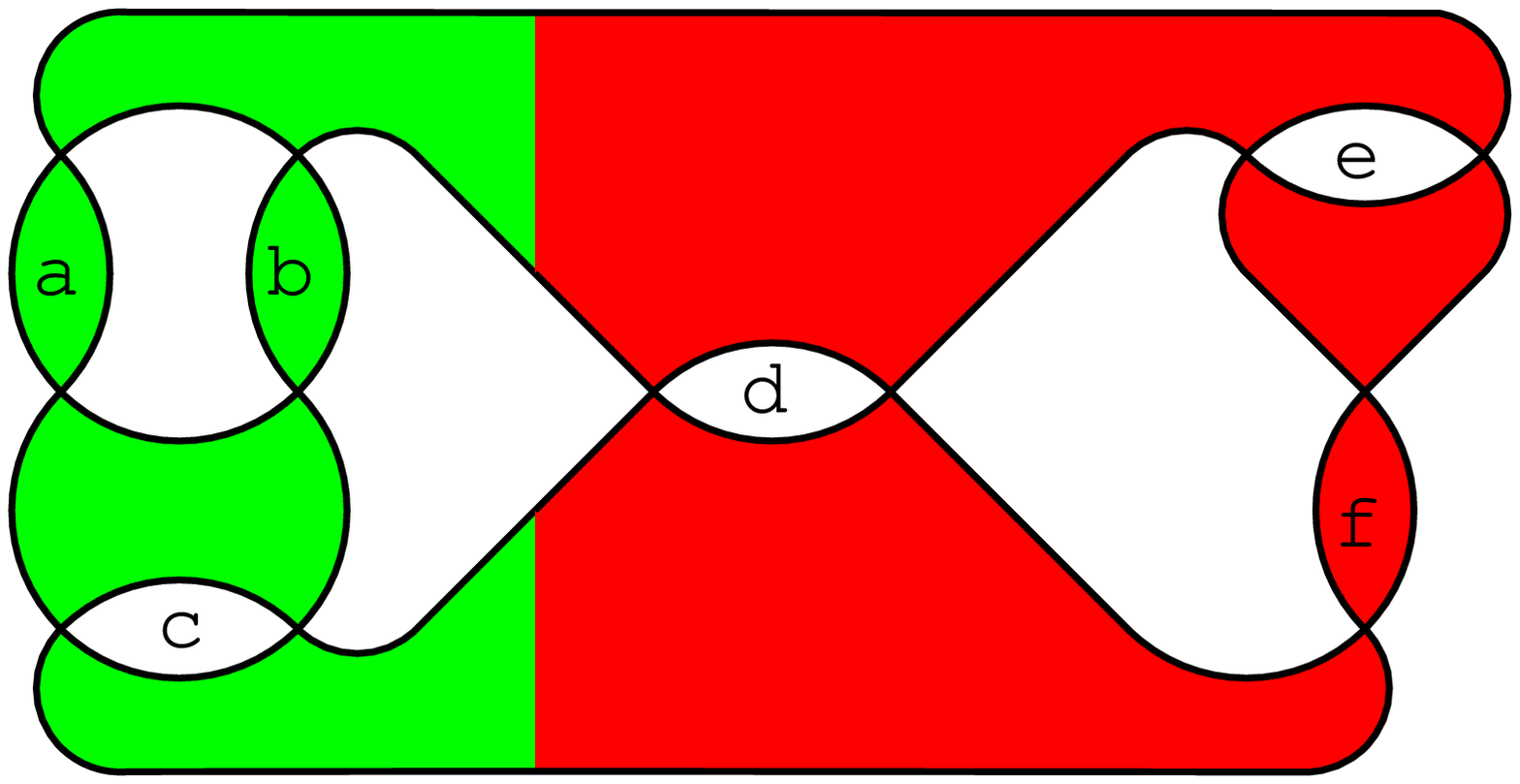}} \quad \scalebox{0.40}{\includegraphics{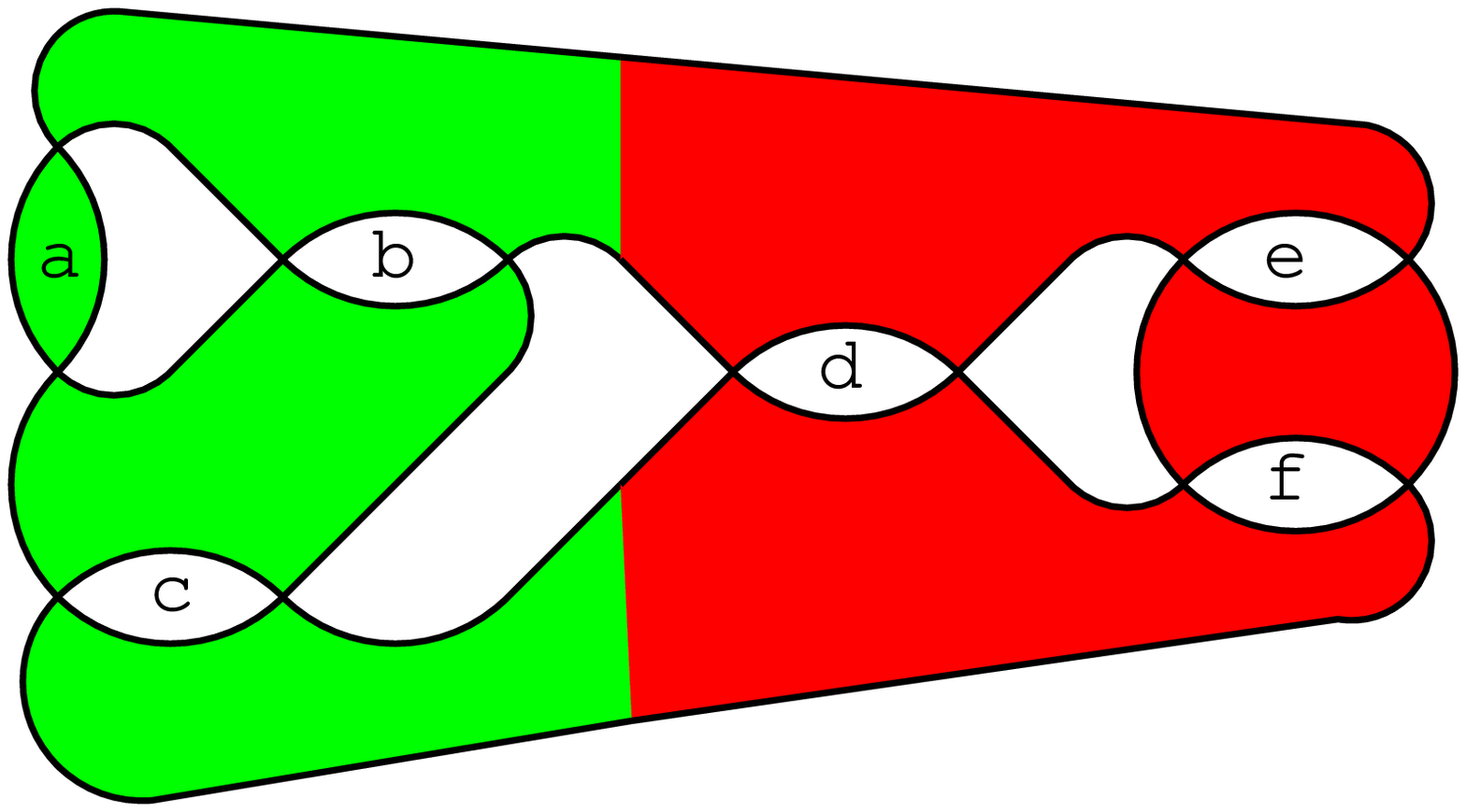}}

\caption{Families formed by six conways with C-functions with thirteen terms $(a_1 a_2 a_3 + a_1 + a_2) (a_4 a_5 a_6 + a_4 + a_5) +
a_3 (a_1 + a_2) (a_5 a_6 + 1)$ and $(a_1 (a_2 + a_3) + 1) (a_4 a_5 + a_5 a_6 + a_6 a_4) + a_3 (a_1 a_2 + 1) (a_5 + a_6)$}

\end{figure}

\begin{figure}

\psfrag{a}{\LARGE{$a_1$}}
\psfrag{b}{\LARGE{$a_2$}}
\psfrag{c}{\LARGE{$a_3$}}
\psfrag{d}{\LARGE{$a_4$}}
\psfrag{e}{\LARGE{$a_5$}}
\psfrag{f}{\LARGE{$a_6$}}

\scalebox{0.40}{\includegraphics{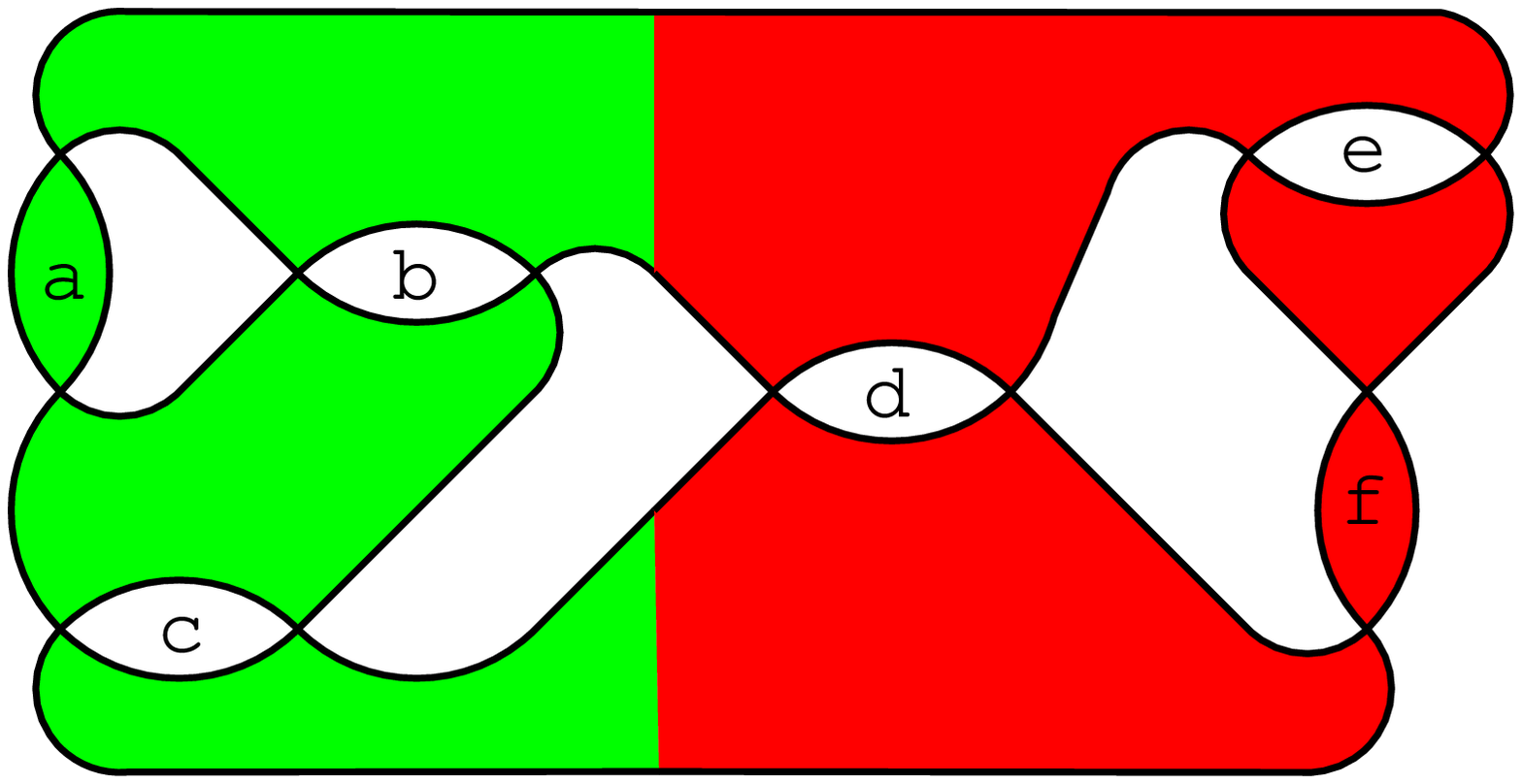}} \quad \scalebox{0.40}{\includegraphics{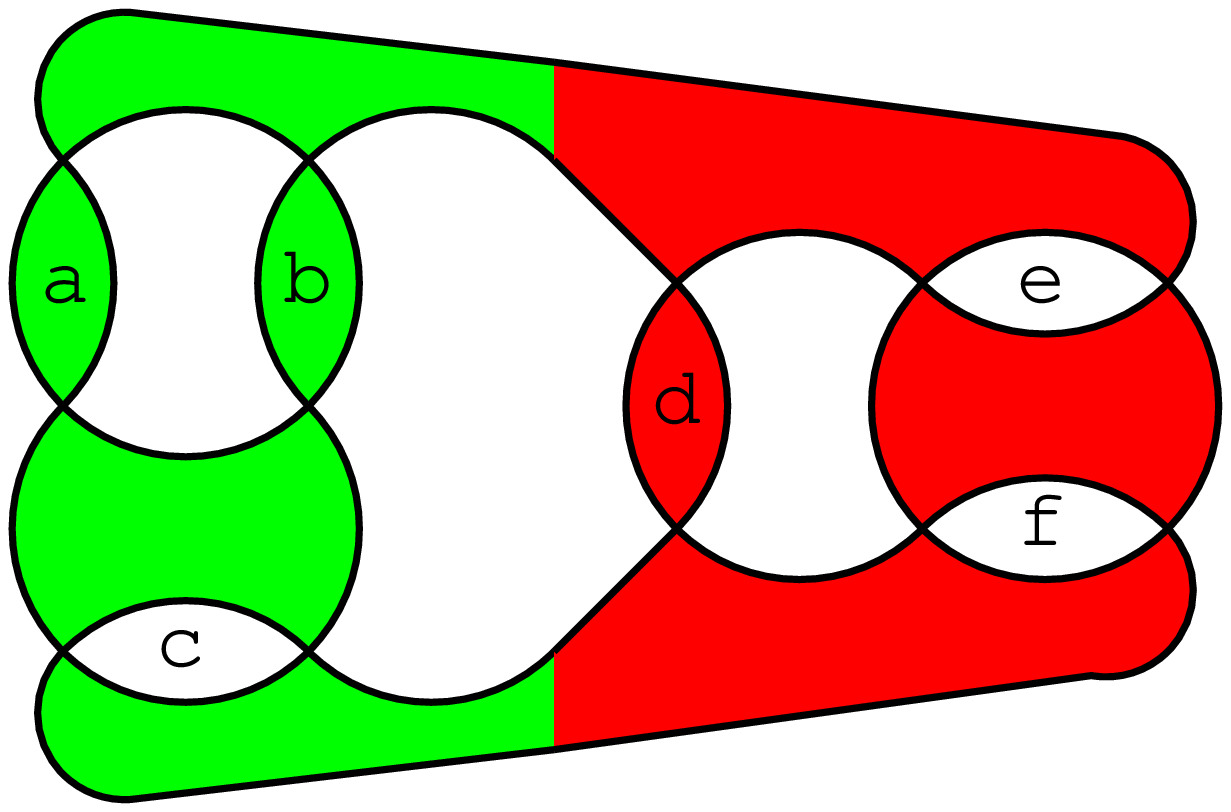}}

\caption{Families formed by six conways with C-functions of thirteen terms $(a_1 (a_2 + a_3) + 1) (a_4 a_5 a_6 + a_4 + a_5) + a_3 (a_1 a_2 + 1) (a_5 a_6 + 1)$ and $(a_1 a_2 a_3 + a_1 + a_2) (a_4 a_5 a_6 + a_5 + a_6) + a_3 (a_1 + a_2) a_4 (a_5 + a_6)$}

\end{figure}

\begin{figure}

\psfrag{a}{\LARGE{$a_1$}}
\psfrag{b}{\LARGE{$a_2$}}
\psfrag{c}{\LARGE{$a_3$}}
\psfrag{d}{\LARGE{$a_4$}}
\psfrag{e}{\LARGE{$a_5$}}
\psfrag{f}{\LARGE{$a_6$}}

\scalebox{0.40}{\includegraphics{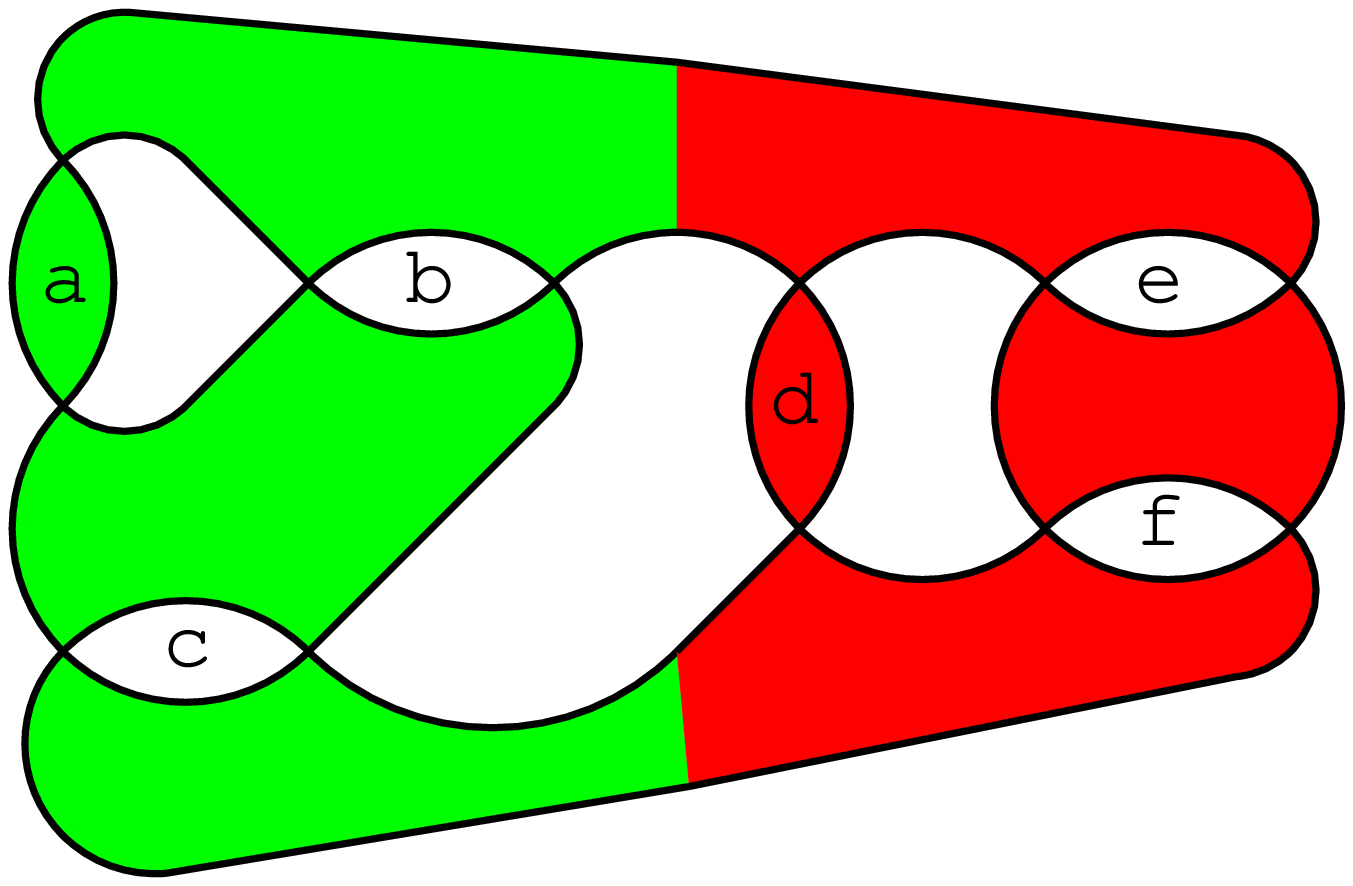}} \quad \scalebox{0.40}{\includegraphics{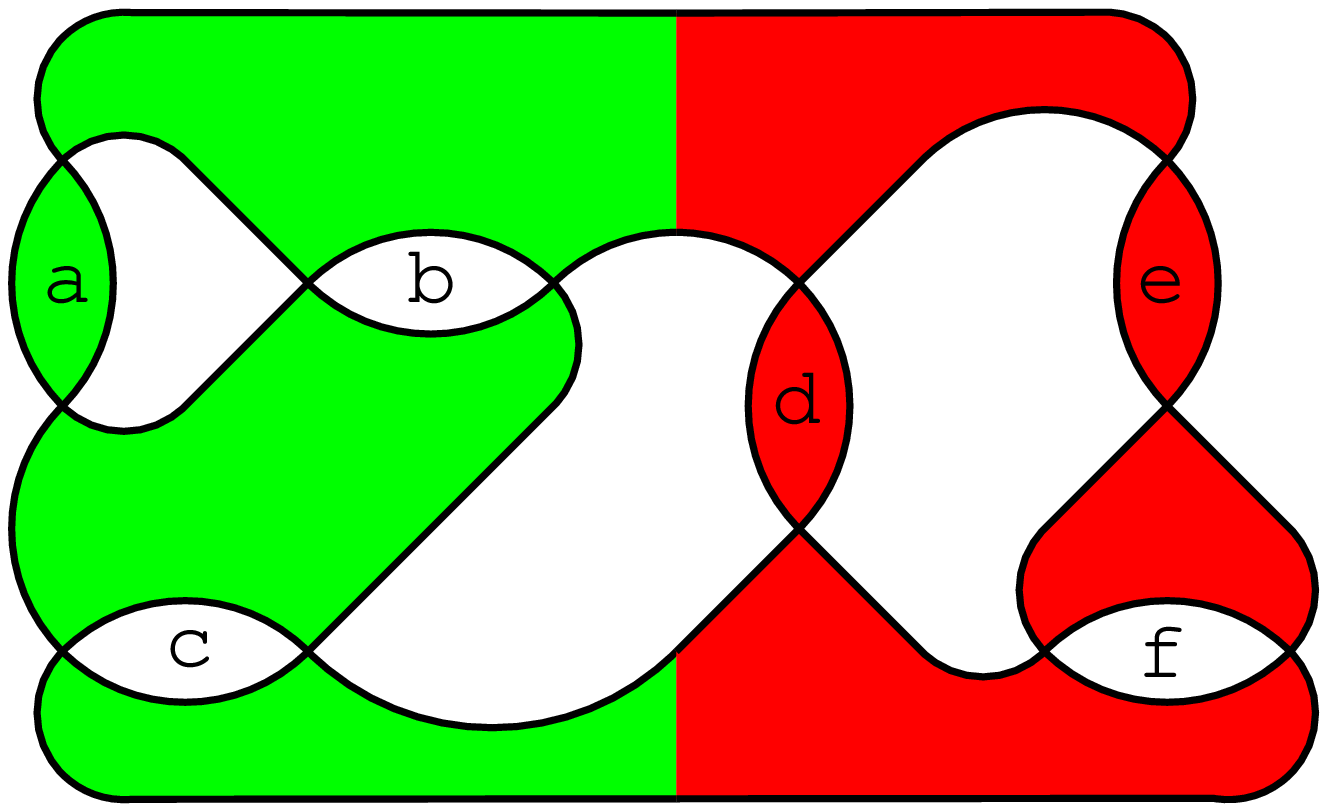}}

\caption{Families formed by six conways. Its C-functions of thirteen terms are $(a_1 (a_2 + a_3) + 1) (a_4 a_5 a_6 + a_5 + a_6) + a_3 (a_1 a_2 + 1) a_4 (a_5 + a_6)$ and $(a_1 (a_2 + a_3) + 1) ((a_4 + a_5) a_6 + 1) + a_3 (a_1 a_2 + 1) a_4 (a_5 a_6 + 1)$}

\end{figure}

\begin{figure}

FAMILIES FORMED BY SIX CONWAYS WITH THE BORROMEAN LINK $C_2^3$ AS SEED. SEVEN CASES.

\psfrag{a}{\LARGE{$a_5$}}
\psfrag{b}{\LARGE{$a_2$}}
\psfrag{c}{\LARGE{$a_3$}}
\psfrag{d}{\LARGE{$a_4$}}
\psfrag{e}{\LARGE{$a_6$}}
\psfrag{f}{\LARGE{$a_1$}}

\hfil \scalebox{0.6}{\includegraphics{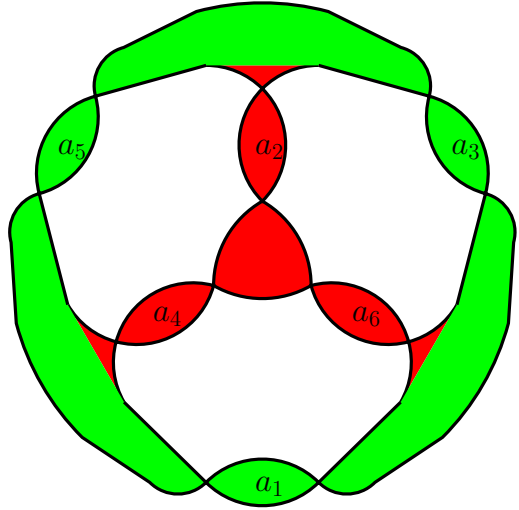}} \hfil

\caption{Family formed by six conways. The C-function has sixteen terms
$(a_1 + a_3 + a_5) (a_2 a_4 + a_4 a_6 + a_6 a_2) + a_3 a_5 (a_4 + a_6) + a_5 a_1 (a_2 + a_6) + a_1 a_3 (a_2 + a_4) + a_1 a_3 a_5$.}

\end{figure}

\begin{figure}

\centering

\psfrag{a}{\LARGE{$a_5$}}
\psfrag{b}{\LARGE{$a_2$}}
\psfrag{c}{\LARGE{$a_3$}}
\psfrag{d}{\LARGE{$a_4$}}
\psfrag{e}{\LARGE{$a_6$}}
\psfrag{f}{\LARGE{$a_1$}}

\scalebox{0.4}{\includegraphics{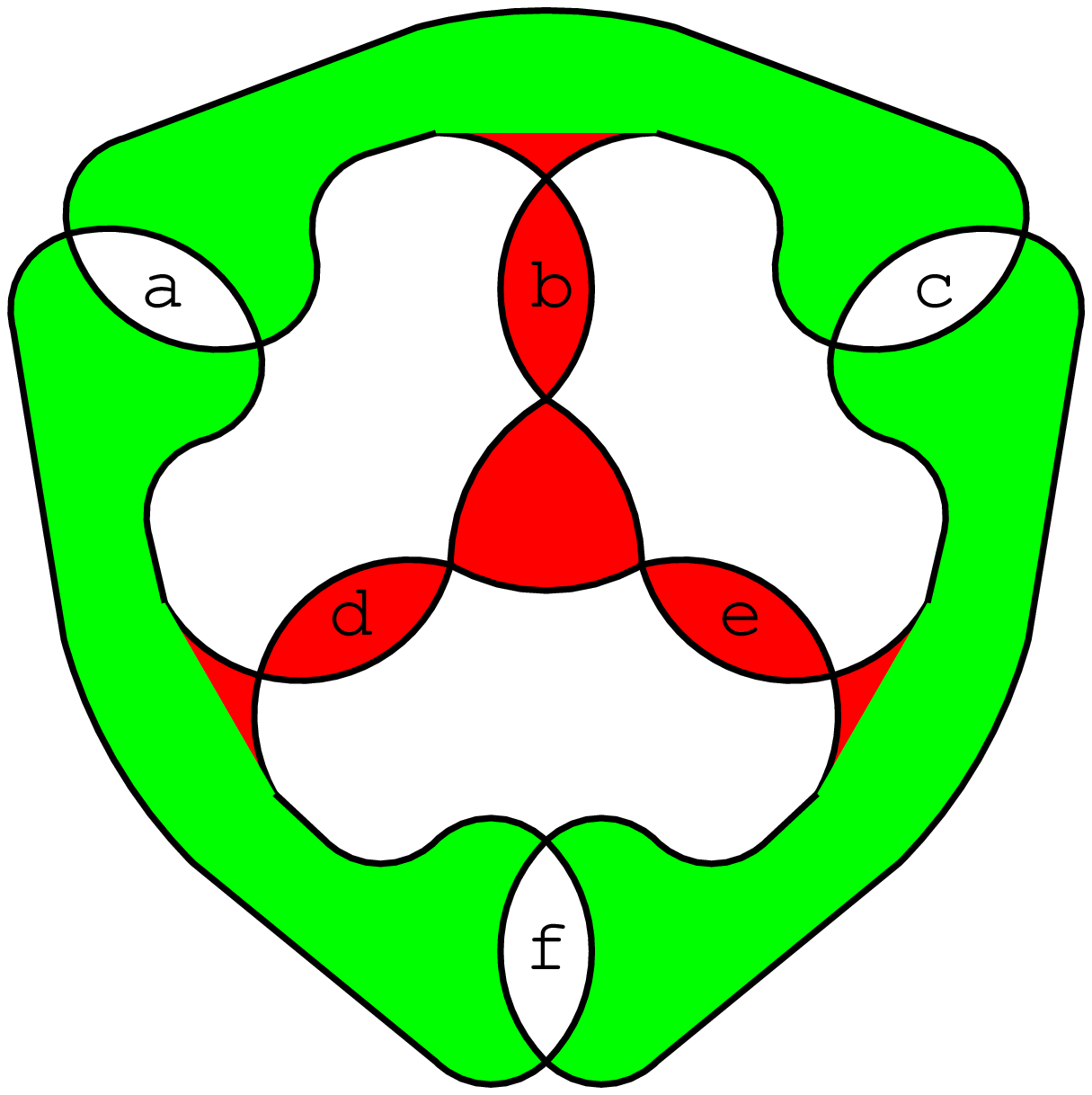}}
\caption{Family formed by six conways. The C-function has 16 monomials. The C-function is $(a_1 a_3 + a_3 a_5 + a_5 a_1) (a_2 a_4 + a_4 a_6 + a_6 a_2) + a_1 (a_4 + a_6) + a_3 (a_2 + a_6) + a_5 (a_2 + a_4) + 1$.}

\end{figure}

\begin{figure}

\centering

\psfrag{a}{\LARGE{$a_5$}}
\psfrag{b}{\LARGE{$a_2$}}
\psfrag{c}{\LARGE{$a_3$}}
\psfrag{d}{\LARGE{$a_4$}}
\psfrag{e}{\LARGE{$a_6$}}
\psfrag{f}{\LARGE{$a_1$}}

\scalebox{0.4}{\includegraphics{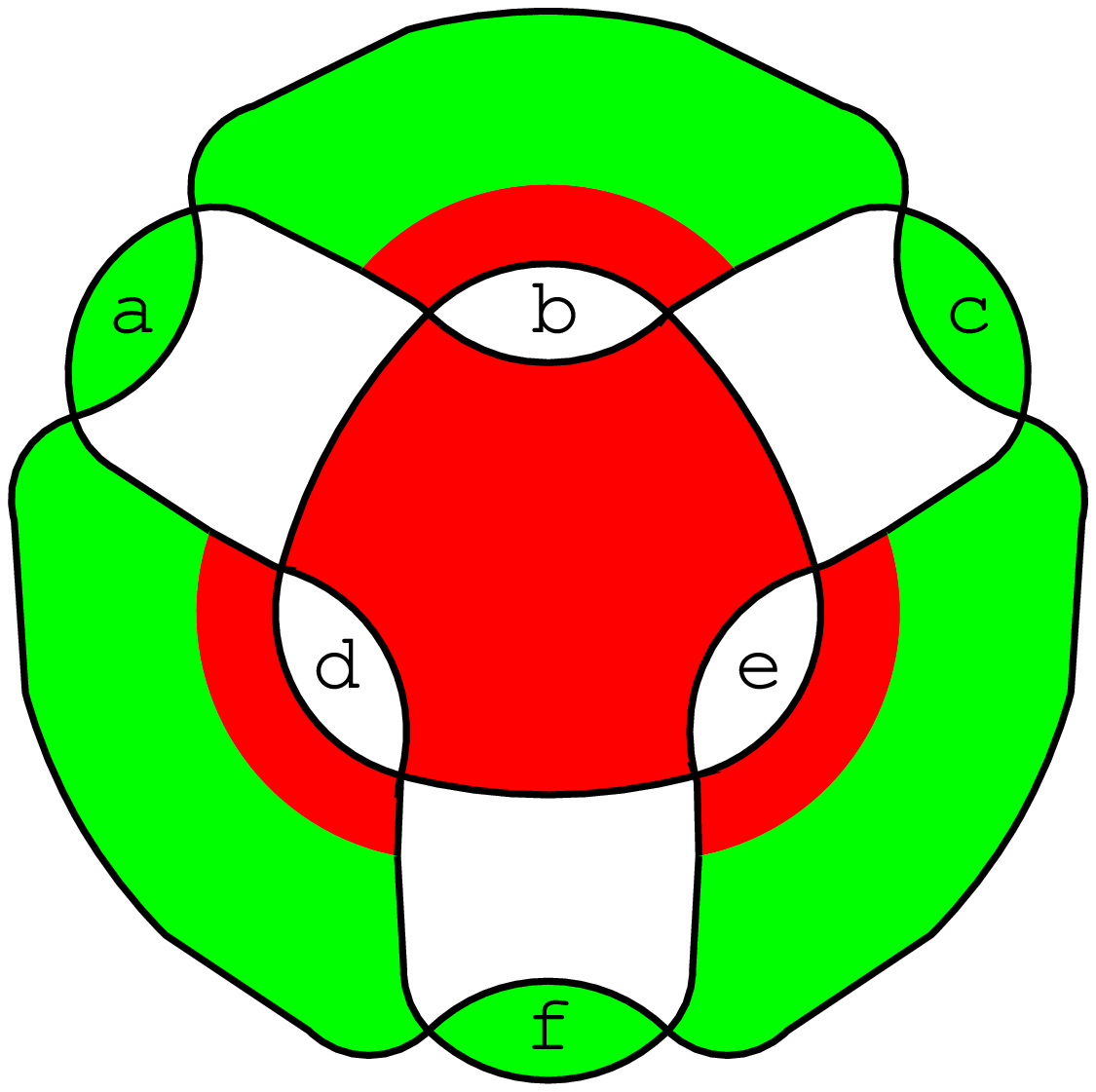}}
\caption{Family formed by six conways. The C-function has 16 monomials. The C-function is $(a_1 + a_3 + a_5) (a_2 + a_4 + a_6) + a_3 a_5 a_2 (a_6 + a_4) + a_5 a_1 a_4 (a_6 + a_2 ) + a_1 a_3 a_6(a_2 + a_4) + a_1 a_3 a_5 a_2 a_4 a_6.$}

\end{figure}

\begin{figure}

\centering

\psfrag{a}{\LARGE{$a_5$}}
\psfrag{b}{\LARGE{$a_2$}}
\psfrag{c}{\LARGE{$a_3$}}
\psfrag{d}{\LARGE{$a_4$}}
\psfrag{e}{\LARGE{$a_6$}}
\psfrag{f}{\LARGE{$a_1$}}

\scalebox{0.4}{\includegraphics{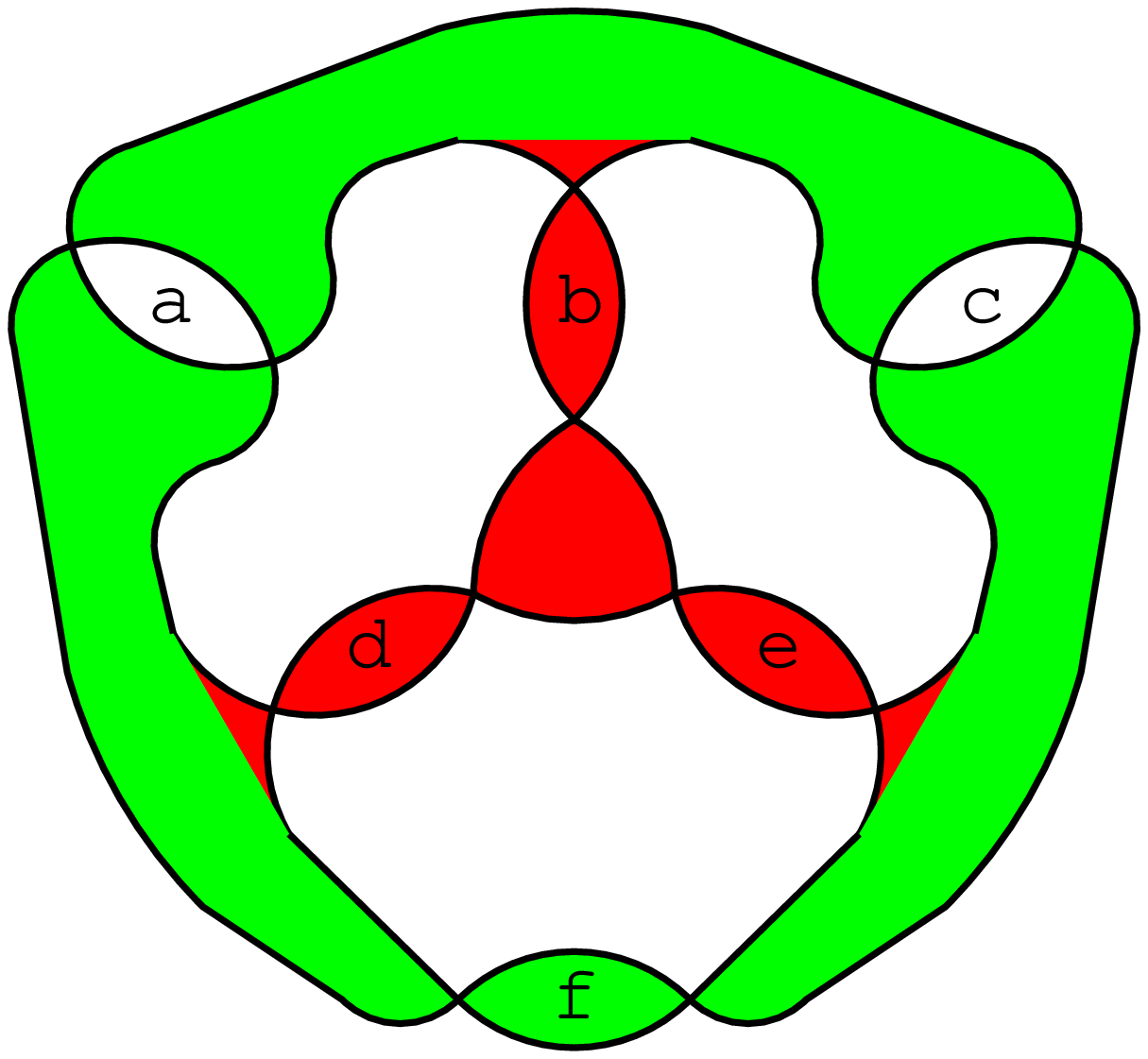}}
\caption{Family formed by six conways. The C-function has 16 monomials. The C-function is $(a_1 a_3 a_5 + a_3 + a_5) (a_2 a_4 + a_4 a_6 + a_6 a_2) + a_1 a_3 (a_2 + a_6) + a_1 a_5 (a_4 + a_6) + a_1 + a_4 + a_6$. }
\end{figure}

\begin{figure}

\psfrag{a}{\LARGE{$a_5$}}
\psfrag{b}{\LARGE{$a_2$}}
\psfrag{c}{\LARGE{$a_3$}}
\psfrag{d}{\LARGE{$a_4$}}
\psfrag{e}{\LARGE{$a_6$}}
\psfrag{f}{\LARGE{$a_1$}}

\hfil \scalebox{0.6}{\includegraphics{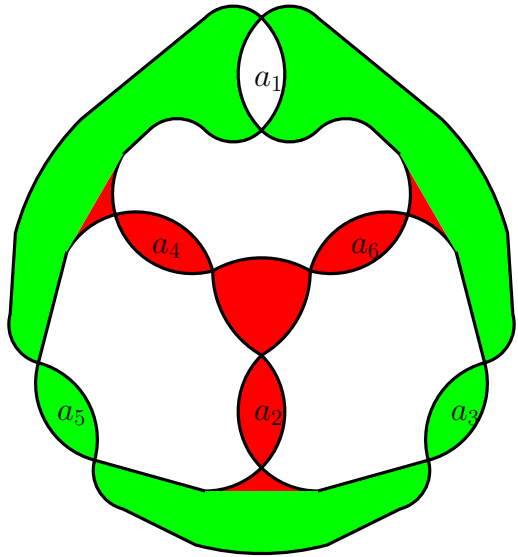}} \hfil

\caption{Family formed by six conways. The C-function with sixteen terms is
$(a_2 a_4 + a_4 a_6 + a_6 a_2) + a_5 (a_2 + a_6) + a_3 (a_2 + a_4) + a_3 a_5 + a_1 [( a_3 + a_5) (a_2 a_4 + a_4 a_6 + a_6 a_2) + a_3 a_5 (a_4 + a_6)]$.}

\end{figure}

\begin{figure}

\psfrag{a}{\LARGE{$a_5$}}
\psfrag{b}{\LARGE{$a_2$}}
\psfrag{c}{\LARGE{$a_3$}}
\psfrag{d}{\LARGE{$a_4$}}
\psfrag{e}{\LARGE{$a_6$}}
\psfrag{f}{\LARGE{$a_1$}}

\hfil \scalebox{0.6}{\includegraphics{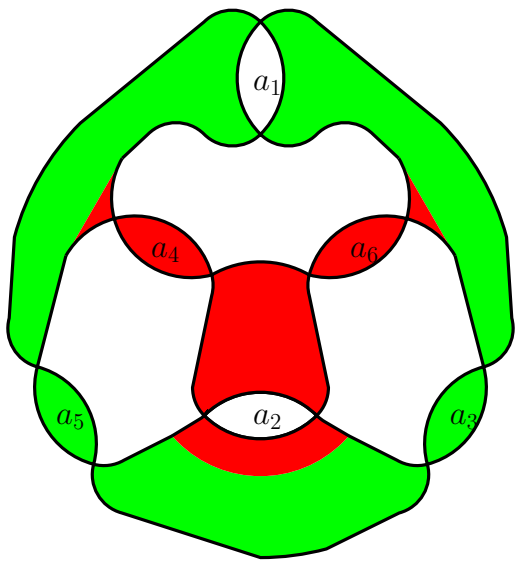}} \hfil

\caption{Family formed by six conways. The C-function has sixteen terms. The C-function is
$(a_2 + a_3) [1 + a_1 a_4 a_5 a_6] + [1 + a_1 a_2 a_3 a_4 ] (a_6 + a_5) + (a_2 + a_3) a_4 (a_6 + a_5) + (a_2 + a_6) a_1 (a_3 + a_5)$.}

\end{figure}

\begin{figure}

\centering

\psfrag{a}{\LARGE{$a_1$}}
\psfrag{b}{\LARGE{$a_2$}}
\psfrag{c}{\LARGE{$a_3$}}
\psfrag{d}{\LARGE{$a_4$}}
\psfrag{e}{\LARGE{$a_5$}}
\psfrag{f}{\LARGE{$a_6$}}

\scalebox{0.6}{\includegraphics{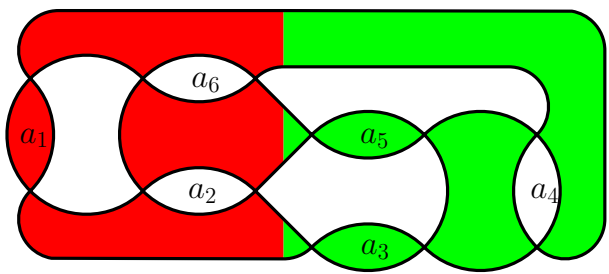}}

\caption{Family formed by six conways. The C-function has sixteen monomials, $\left( a_1 a_2 a_6 + a_2 + a_6 \right) \left(a_3 a_4 a_5 + a_3 + a_5 \right) + a_1 a_4  + (a_1 a_2 + a_1 a_6) a_4 a_5 +
(1 + a_1 a_2) a_3 a_4  + 1 + a_1 a_6 $}

\end{figure}

\

{\large \bf ACKNOWLEDGMENTS }

I thank Prof. Eduardo Virue\~na from ESFM, IPN his help to learn the drawing techniques used in this work. The first figure is the resulting projection
of his program for drawing (moving) knots in three dimensions, given just the parametric equations of the knot in three dimensions. My thanks to
Prof. Pablo Lonngi
from UAM-I for correcting the English language of this paper to make it comprehensible to almost everyone with some interest in knots.

\

\end{document}